\documentclass[11pt,leqno]{article} 

\usepackage{fullpage}

\usepackage{amssymb} 
\usepackage{amsmath} 
\date{}
\input xy 
\xyoption{all}
\title{Classification of linearly compact simple rigid 
 superalgebras} 
\author{{\sc Nicoletta Cantarini}\thanks{Dipartimento di Matematica
Pura ed Applicata, Universit\`a di Padova, Padova, Italy
- Partially supported by Progetto di ateneo CPDA071244}
\and\setcounter{footnote}{6}
{\sc Victor G.\ Kac}\thanks{Department of Mathematics, MIT, Cambridge,
Massachusetts 02139, USA
- Partially supported by an NSF grant}}

\newtheorem{theorem}{Theorem}[section] 
\newtheorem{lemma}[theorem]{Lemma} 
\newtheorem{corollary}[theorem]{Corollary} 
\newtheorem{proposition}[theorem]{Proposition} 
\newtheorem{definition}[theorem]{Definition} 
\newtheorem{remark}[theorem]{Remark}

\newtheorem{example}[theorem]{Example}

\def\Z{\mathbb{Z}} 

\def\g{\mathfrak{g}}
\def\a{\mathfrak{a}}

\def\fg{\mathfrak{g}}
\def\h{\mathfrak{h}} 
\def\s{\mathfrak{s}}

\def\F{\mathbb{F}}
\def\FF{\mathbb{F}}

\def\0{\bar{0}}
\def\1{\bar{1}}
\def\O{{\mathcal O}}
\def\growth{~\mbox{growth}~}
\def\size{~\mbox{size}~}

\numberwithin{equation}{section}

\makeatletter

\def\enumerate{%
  \ifnum \@enumdepth >\thr@@\@toodeep\else
    \advance\@enumdepth\@ne
    \edef\@enumctr{enum\romannumeral\the\@enumdepth}%
      \list
        {\csname label\@enumctr\endcsname}%
        {\usecounter\@enumctr
          \addtolength{\leftmargin}{-\leftmargin}
          \settowidth{\labelwidth}{(99)}
          \itemindent = \labelwidth
          \addtolength{\itemindent}{\labelsep}
        \listparindent=1em      
          \def\makelabel##1{{##1}\hfill}
          }%
  \fi}

\makeatother

\begin{document} 
\maketitle 
\begin{abstract}
The notion of an anti-commutative (resp.\ commutative) rigid superalgebra 
is a natural generalisation of the notion of a Lie (resp.\ Jordan) 
superalgebra. Intuitively rigidity means that small deformations of the 
product under the structural group produce an isomorphic algebra.
In this paper we classify all linearly compact simple 
anti-commutative (resp.\ commutative) rigid superalgebras.
Beyond Lie (resp.\ Jordan) superalgebras the complete list
includes four series and twenty two exceptional superalgebras (resp.\
ten exceptional superalgebras).   
\end{abstract} 
\section*{Introduction}
The Killing-Cartan classification of simple finite-dimensional Lie algebras
over an algebraically closed field $\F$ of characteristic 0 and Cartan's
classification of simple infinite-dimensional Lie algebras of vector fields
can be viewed as two parts of the classification of simple linearly compact Lie
algebras. (Recall that a topological algebra is called linearly compact
if the underlying topological vector space is linearly compact, i.e., it is
isomorphic to the topological direct product of finite-dimensional vector 
spaces with the discrete topology). Thus, a complete list of simple linearly 
compact Lie algebras consists of classical finite-dimensional series: $sl_m$,
$so_m$, $sp_m$ (m even), five exceptional (finite-dimensional) Lie algebras: 
$E_6$,
$E_7$, $E_8$, $F_4$, $G_2$, and four infinite-dimensional series of Lie 
algebras of formal vector fields: $W_m$, $S_m$, $H_m$ ($m$ even) and $K_m$
($m$ odd).

This classification was extended to the Lie superalgebra case in \cite{K2}
and \cite{K3}. A complete list of simple linearly compact Lie superalgebras
consists of classical finite-dimensional series $sl(m,n)/\delta_{m,n}I$,
$osp(m,n)$ (which include the classical Lie algebra series when $m$ or $n$
are 0), two ``strange'' series $p(n)$, $q(n)$, three exceptional finite-dimensional
Lie superalgebras $D(2,1;\alpha)$, $F(4)$, $G(3)$, in addition to the five
exceptional Lie algebras, eleven series of Lie superalgebras of vector fields
(which are finite-dimensional if $m=0$): $W(m,n)$, $S(m,n)$, $H(m,n)$ ($m$ even),
$K(m,n)$ ($m$ odd), $HO(m,m)$ ($m\geq 2$), $SHO(m,m)$ ($m\geq 3$),
$KO(m,m+1)$ ($m\geq 1$), $SKO(m,m+1;\beta)$ ($m\geq 2$), $\tilde{S}(0,n)$
($n\geq 4$, even), $SHO^\sim(m,m)$ ($m\geq 2$, even), $SKO^\sim(m,m+1)$
($m\geq 3$, odd) (which include the Lie algebra series $W_m=W(m,0)$,
$S_m=S(m,0)$, $H_m=H(m,0)$, $K_m=K(m,0)$), and five exceptional
infinite-dimensional Lie superalgebras: $E(1,6)$, $E(3,6)$, $E(3,8)$,
$E(4,4)$, $E(5,10)$.

Closely related, via the Tits-Kantor-Koecher (TKK) construction, is the
classification of simple linearly compact Jordan superalgebras, obtained in
\cite{CantaK2}. It is based on the observation that isomorphism classes
of such Jordan superalgebras correspond bijectively to conjugacy classes
of short $sl_2$ subalgebras in the Lie superalgebra $Der L$ of continuous
derivations of a simple linearly compact Lie superalgebra $L$. Recall that
an $sl_2$ subalgebra of a Lie superalgebra is called short if it has
a basis $e, h, f$, such that $[h,e]=-e$, $[h,f]=f$, $[e,f]=h$ and the
only eigenvalues of $ad\, h$ on $L$ are 0,1 or $-1$. 
Given such a subalgebra in $Der L$, the
corresponding simple Jordan superalgebra is the $-1$-eigenspace $J$ of $ad\, h$
in $L$, endowed with the product
\begin{equation}
a\circ b=[[f,a],b], ~~a, b\in J.
\label{jp}
\end{equation}

The complete list of simple linearly compact Jordan superalgebras
consists of classical finite-dimensional series $gl(m,n)_+$, $osp(m,n)_+$
($n$ even), $(m,n)_+$ ($n$ even), two ``strange'' series $p(n)_+$,
$q(n)_+$, the series $JP(m,n)$ (which is finite-dimensional if $m=0$), four
 finite-dimensional exceptional superalgebras $E$, $F$, $D_t$, $K$, and two 
exceptional infinite-dimensional superalgebras, $JCK$ and $JS$. 
The finite-dimensional part of this classification goes back to \cite{A}
and  \cite{K1} 
in the non-super and super cases respectively.

The constructions of finite-dimensional classical and ``strange'' series
of Lie and Jordan superalgebras are parallel and very simple. Recall that
all finite-dimensional associative simple superalgebras consist of two series:
$mat(m,n)$ (the superalgebra of endomorphisms of the superspace $\F^{m|n}$)
and, in the case $m=n$, its subalgebra $qmat(n,n)$, consisting of
matrices commuting with the matrix 
$\left(\begin{array}{cc}
0 & 1\\
1 & 0
\end{array}
\right)$. Given an associative superalgebra $A$, one denotes by $A_{-}$
(resp.\ $A_+$) the space $A$ with the product 
\begin{equation}
[a, b]=ab-(-1)^{p(a)p(b)}ba ~~~~
({\mbox{resp.}}~ a\circ b=ab+(-1)^{p(a)p(b)}ba).
\label{lie}
\end{equation}
Here and further, $p(a) \in \Z/2\Z$ denotes the parity of $a \in A$.
Then $sl(m,n)/\delta_{m,n}I$ (resp.\ $q(n)$) is the derived algebra divided
by the center for $mat(m,n)_{-}$ (resp.\ $qmat(n,n)_-$) and $gl(m,n)_+=
mat(m,n)_+$.
The Lie (resp.\ Jordan) superalgebra $osp(m,n)$ (resp.\ $osp(m,n)_+$) is the
subalgebra of skewsymmetric (resp.\ symmetric) endomorphisms in $mat(m,n)_-$
(resp.\ $mat(m,n)_+$) with respect to a non-degenerate even supersymmetric 
bilinear form on $\F^{m|n}$. The Lie (resp.\ Jordan) superalgebra $p(n)$ 
(resp.\ $p(n)_+$) is the derived algebra divided by the center for the
subalgebra of skew-symmetric endomorphisms in $mat(n,n)$ (resp.\
the subalgebra of symmetric endomorphisms) with respect to an odd
non-degenerate bilinear form on $\F^{n|n}$. Finally, $(m,n)_+=\F e\oplus
\F^{m|n}$, where $e$ is the identity element and $a\circ b=(a|b)e$
for $a,b\in\F^{m|n}$, where $(\cdot|\cdot)$ is an even non-degenerate
supersymmetric bilinear form on $\F^{m|n}$.

The series $JP(m,n)$ is obtained from the generalized Poisson superalgebras
$P(m,n)$ via the following Kantor-King-McCrimmon (KKM) construction
\cite{Ka2}, \cite{KMCC}. Recall \cite{CantaK2}
that a generalized Poisson superalgebra $P$ is
a unital commutative associative superalgebra, endowed with a Lie
superalgebra bracket $\{\cdot,\cdot\}$, satisfying the generalized
Leibniz rule
\begin{equation}
\{a,bc\}=\{a,b\}c+(-1)^{p(a)p(b)}b\{a,c\}+D(a)bc,
\label{leibniz}
\end{equation} 
where $D(a)=\{1,a\}$. The associated Jordan superalgebra is 
$JP=P\oplus \bar{P}$, where $\bar{P}$ is a copy of $P$ with reversed parity, 
endowed with the following Jordan product, where $a, b\in P$, $\bar{a},
\bar{b}\in\bar{P}$:
\begin{equation}
a\circ b=ab, ~~\bar{a}\circ b=\overline{ab}, 
~~a\circ\bar{b}=(-1)^{p(a)}\overline{ab},
~~\bar{a}\circ\bar{b}=(-1)^{p(a)}\{a,b\}_D,
\label{jordan}
\end{equation}
where $\{a,b\}_D=\{a,b\}-\frac{1}{2}(aD(b)-D(a)b)$.

The KKM construction allows one to derive a classification of simple linearly
compact generalized Poisson superalgebras from that of simple linearly
compact Jordan superalgebras \cite{CantaK2}. The result is that, up to gauge
equivalence, any such generalized Poisson superalgebra is one of the $P(m,n)$.
 Recall that a bracket $\{a,b\}^{\varphi}=\varphi^{-1}\{
\varphi a,\varphi b\}$, where $\varphi$ is an even invertible element,
is called gauge equivalent to the bracket $\{\cdot,\cdot\}$. Also,
the generalized Poisson superalgebra
$P(m,n)$ is the associative commutative superalgebra of formal
power series in the even indeterminates $p_1, \dots,p_k, q_1, \dots, q_k$
(resp.\ $p_1, \dots,p_k, q_1, \dots, q_k, t$) if $m=2k$ (resp.\ $m=2k+1$)
and in the $n$ odd indeterminates $\xi_1, \dots, \xi_n$,
endowed with the bracket:
\begin{equation}
\{a,b\}=\sum_{i=1}^k(\frac{\partial a}{\partial p_i}
\frac{\partial b}{\partial q_i}-\frac{\partial a}{\partial q_i}
\frac{\partial b}{\partial p_i})+(-1)^{p(a)}
\sum_{i=1}^n\frac{\partial a}{\partial\xi_i}\frac{\partial b}{\partial\xi_i}
+(2-E)(a)\frac{\partial b}{\partial t}-\frac{\partial a}{\partial t}
(2-E)(b),
\label{poisson}
\end{equation}
where $E=\sum_{i=1}^k(p_i\frac{\partial}{\partial p_i}+q_i\frac{\partial}
{\partial q_i})+\sum_{i=1}^n\xi_i\frac{\partial}{\partial\xi_i}$
(the last two terms in (\ref{poisson}) vanish if $m$ is even).

Our first basic idea is to generalize the notion of a Jordan
superalgebra as follows. Note that a commutative superalgebra is a vector
superspace $J$ together with an even element $\mu$ of $Hom(S^2J,J)$,
so that $a\circ b=\mu(a \otimes b)$. Define the structural Lie superalgebra
$Str(J,\mu)$ as a subalgebra of the Lie superalgebra $End(J,J)_-$
generated by all operators of left multiplication $\mu_a(b)=\mu(a \otimes b)$,
$a\in J$, and denote by $R(J,\mu)$ the minimal submodule of the
$Str(J,\mu)$-module $Hom(S^2J,J)$, containing $\mu$. One of the basic
properties of a Jordan superalgebra is \cite{K1}
\begin{equation}
R(J,\mu)=Str(J,\mu) \mu+\F\mu.
\label{Str}
\end{equation}
Note that in the case $J$ is a finite-dimensional Jordan algebra,
$R(J,\mu)$ is the span of the orbit of $\F\mu$ under the structural group, and
 property 
(\ref{Str}) 
means that this orbit is dense in $R(J,\mu)$. Thus, (\ref{Str})  means 
a certain rigidity property, namely, ``small'' deformations of the product
by the structural group produce an isomorpic algebra.

We define a {\em commutative rigid superalgebra} 
to be a commutative
superalgebra $(J,\mu)$, satisfying (\ref{Str}). The first result
of our paper is
\begin{theorem}\label{01} Any simple linearly compact commutative rigid  
superalgebra, which is not a Jordan superalgebra, 
is isomorphic either to one of the five infinite-dimensional superalgebras
$JS_{1,1}$, $JSHO_{2,2}$, $JSKO_{1,2}$, $JS^\alpha_{1,8}$ for $\alpha =0,1$,
or to one of the five finite-dimensional superalgebras $JS_{0,2}$,
 $JW_{0,4}$, $JW_{0,8}$, $JS_{0,8}$, $JS_{0,16}$. 
\end{theorem}

In order to explain the construction of the superalgebras
appearing in the statement of Theorem~0.1, introduce some
notation.  Given a commutative associative superalgebra $A$ and its
derivation~$D$, denote by $AD$ the vector superspace of
derivations of $A$ of the form $fD$, where $f \in A$.  
Given even derivations $D_i$ of $A$ and elements $f_i \in A$ $(i=1,2)$, let
\begin{equation}
\label{eq:07a}  
[f_1D_1 \, , \, f_2D_2]_\pm = f_1 D_1 (f_2) D_2 \pm    
(-1)^{p(f_1)p(f_2)} f_2D_2 (f_1)D_1. 
\end{equation}
Given odd derivations $D_i$ of $A$ and elements $f_i \in A$ $(i=1,2)$, let
\begin{equation}
\label{eq:08a}  [f_1D_1 \, , \, f_2D_2]_+ =
f_1 D_1 (f_2) D_2 + (-1)^{p(f_1)+1)(p(f_2)+1)}f_2D_2 (f_1)D_1
\end{equation}
Recall that the product $[fD, gD]_+$ defines on $AD$
a Jordan superalgebra structure for any odd derivation $D$ of $A$
\cite[Example 4.6]{CantaK2}. 

Our first example is $JS_{1,1} = \FF [[x]]\frac{d}{dx}$ with the product
\begin{equation}
  \label{eq:0.7}
  f \frac{d}{dx} \circ g  \frac{d}{dx}
    = \Big[ f  \frac{d}{dx}\, , \, 
         g  \frac{d}{dx}\Big]_+ \, ,
\end{equation}
%
which we call the Beltrami algebra. 

We shall denote by $\bar{A}$ a copy of A with reversed parity, and by 
$\bar{a}\in \bar{A}$ the element corresponding to $a\in A$.

The next example is $JSHO_{2,2} = \F[[p,q]]  \frac{\partial}{\partial p}
\oplus \overline{\F[[p,q]]}$
with the product $(f,g \in \F[[p,q]] \, , \,\bar{f}, \bar{g}
\in \overline{\F[[p,q]]})$:
\begin{equation}
\label{eq:09a}
  f  \frac{\partial}{\partial p} \circ g  
    \frac{\partial}{\partial p}
    = \Big[ f  \frac{\partial}{\partial p}, 
        g  \frac{\partial}{\partial p}\Big]_+ \, , \, f
        \frac{\partial}{\partial p} \circ \bar{g} 
    = \overline{f \frac{\partial g}{\partial p}} \, , \, \bar{f} \circ \bar{g} 
          = \Big( \frac{\partial f}{\partial p} \frac{\partial g}{\partial q}
            -  \frac{\partial f}{\partial q} 
           \frac{\partial g}{\partial p}\Big)\frac{\partial}{\partial p}\, .
\end{equation}

The next example is $JSKO_{1,2} = \FF[[x]]
\frac{d}{d x} \oplus \overline{\FF [[x]]}$ with the product
($f,g \in \FF [[x]] \, , \, \bar{f},\bar{g} \in \overline{\FF [[x]]}$)
\begin{equation}
  \label{eq:0.8}
   f \frac{d}{d x} \circ g   \frac{d}{d x}  =
     \Big[ f \frac{d}{d x} \, , \,  g   \frac{d}{d x}\Big]_+\, , \, 
       f \frac{d}{d x} \circ \bar{g} = \overline{f  \frac{dg}{d  x}}\,,\,
         \bar{f} \circ \bar{g} = 2 \Big[ f \frac{d}{d x}\, , \, 
           g   \frac{d}{d x} \Big]_- \, .  
\end{equation}

The remaining two infinite-dimensional examples are
$JS^\alpha_{1,8}$ for $\alpha = 0,1$, constructed as follows.  Let
$A =\FF [[x,\xi_1,\xi_2]]$, where $\xi_1$ and $\xi_2$ are
Grassmann indeterminates, and let
\begin{displaymath}
  D_1 = \frac{\partial}{\partial \xi_1} +\xi_1
         \frac{\partial}{\partial x} + \alpha \xi_2  
          \frac{\partial}{\partial x}\, , \, 
  D_2 =  \frac{\partial}{\partial \xi_2} + x \xi_2  
          \frac{\partial}{\partial x } -\xi_1\xi_2 
            \frac{\partial}{\partial \xi_1}\, .
\end{displaymath}
Then $JS^\alpha_{1,8} = A D_1 \oplus A D_2$, 
with the product
$(f,g \in A)$:
\begin{equation}
  \label{eq:0.9}
  f D_i \circ g D_j = [f D_i, g D_j]_+
     \, , \quad i,j = 1,2\, .
\end{equation}
The four finite-dimensional superalgebras $JW_{0,4}$, $JW_{0,8}$, $JS_{0,8}$
and $JS_{0,16}$ are constructed in a 
similar way (see Example \ref{JDD}). Finally, 
 $JS_{0,2}=\F a+\F b$ with products: $a^2=b$, $b^2=a$, 
$a\circ b=0$.

Our second basic idea is to study odd type superalgebras, aimed at the 
classification of odd generalized Poisson superalgebras. Odd Poisson
superalgebras have appeared both in mathematics and physics literature,
like the work of Gerstenhaber on Hochshild cohomology, and the work on the 
Batalin-Vilkovisky quantization.

By definition, an {\em odd type superalgebra} is a superalgebra with reversed parity. Equivalently, it is $(J,\mu)$ with an odd $\mu\in Hom(S^2J,J)$.
The advantage of this notion is that, given an anti-commutative superalgebra
$A$ with product $\bullet$, reversing the parity of $A$
and defining the new product
\begin{equation}
a\circ b=(-1)^{p(a)}a\bullet b,
\label{odd}
\end{equation}
produces a commutative odd type superalgebra $\bar{A}$, and viceversa,
reversing the parity, we pass from a commutative odd type superalgebra
to an anti-commutative superalgebra. 

An {\em odd generalized Poisson superalgebra} is defined as a unital
associative commutative superalgebra, endowed with a bracket
$\{\cdot,\cdot\}$, which is a Lie superalgebra bracket with respect to
reversed parity, such that the following generalized odd Leibniz rule
holds:
\begin{equation}
\{a,bc\}=\{a,b\}c+(-1)^{(p(a)+1)p(b)}b\{a,c\}+(-1)^{p(a)+1}D(a)bc,
\label{oddleibniz}
\end{equation}
where $D(a)=\{1,a\}$. An {\em odd Poisson superalgebra} is a special case
of this, when $D=0$.

The basic examples of odd generalized Poisson superalgebras are $PO(n,n)$
(resp.\ $PO(n,n+1)$), which are associative commutative superalgebras
of formal power series in $n$ even indeterminates $x_1, \dots, x_n$
and $n$ (resp.\ $n+1$) odd indeterminates $\xi_1, \dots, \xi_n$
(resp.\ $\xi_1,\dots, \xi_n, \tau$), endowed with the following
bracket:
\begin{equation}
\{a,b\}=\sum_{i=1}^n(\frac{\partial a}{\partial x_i}
\frac{\partial b}{\partial \xi_i}+(-1)^{p(a)}
\frac{\partial a}{\partial \xi_i}\frac{\partial b}{\partial x_i})
+(E-2)(a)\frac{\partial b}{\partial\tau}+(-1)^{p(a)}
\frac{\partial a}{\partial\tau}(E-2)(b),
\label{PO}
\end{equation}
where $E=\sum_{i=1}^n(x_i\frac{\partial}{\partial x_i}
+\xi_i\frac{\partial}{\partial\xi_i})$ (in $PO(n,n)$ case
the last two terms in (\ref{PO}) vanish). Note also that
$PO(n,n)$ is an odd Poisson superalgebra.

We define a {\em commutative rigid odd type superalgebra} to be an odd 
type commutative
superalgebra $(J,\mu)$, where $\mu$ is an odd element of
$Hom(S^2J,J)$, satisfying condition (\ref{Str}). This is equivalent, by 
reversing the parity and redefining the product as in (\ref{odd}), to the 
notion of an anti-commutative rigid superalgebra.

Note that any Lie superalgebra $\g$ 
is an anti-commutative rigid superalgebra.
Indeed, it is easy to see that $Str(\g,\mu)=\g$ and 
$R(\g,\mu)=\F\mu$. Hence any Lie superalgebra is an anti-commutative rigid
superalgebra.
It turns out that there are many more anti-commutative rigid superalgebras.
At least one of these algebras appears in the theory of PDEs under the name
Jacobi-Mayer bracket, which is given on the factorspace by $\F1$
of the space of formal power series in $x,y,z$
by the following formula (see \cite{D2}):
$$\{f,g\}=\det
\left(
\begin{array}{ccc}
\frac{\partial f}{\partial x} & \frac{\partial f}{\partial y} & \frac{\partial f}{\partial z}\\
\frac{\partial g}{\partial x} & \frac{\partial g}{\partial y} & \frac{\partial g}{\partial z}\\
0 & -x & 1
\end{array}
\right).$$
Our second main result is
\begin{theorem}
\label{02}
Any simple linearly compact anti-commutative rigid superalgebra,
which is not a Lie superalgebra, is isomorphic
either to one of the four series,
described by Examples \ref{mainexample}, \ref{SHO}, and \ref{OJSKO},  
related to odd generalized
Poisson superalgebras, or to one of the twenty two
exceptions, described by Examples \ref{LSKO'}--\ref{OJW}, 
\ref{OJS_{1,3}}--\ref{OJS}, \ref{JSKO0}--\ref{OJH},
\ref{HO(3,3)}, \ref{OJSHO(4,4)}--\ref{OJSKO(3,4;beta)}.
\end{theorem}

The proof of Theorem \ref{01} and Theorem \ref{02} goes as follows.
Given a linearly compact vector superspace $J$, consider the universal 
$\Z$-graded
Lie superalgebra $W(J)=\prod_{k=-1}^\infty W_k(J)$.
{\em Universality} means that any transitive $\Z$-graded Lie superalgebra 
$\g=\prod_{k=-1}^\infty \g_k$ with $\g_{-1}=J$ canonically embeds in $W(J)$,
and {\em transitivity} means that $[x,\g_{-1}]=0$, $x\in\g_k$, $k\geq 0$,
implies $x=0$. Then $W_k(J)=Hom(S^{k+1}(J),J)$, so that any even (resp.\
odd) element $\mu\in W_1(J)$ defines a commutative superalgebra (resp.\ 
commutative odd type 
superalgebra) structure on $J$ and viceversa. Hence to any commutative
superalgebra (resp.\ odd type superalgebra) structure on $J$
we can canonically associate the Lie superalgebra
$Lie(J,\mu)=\prod_{k=-1}^\infty \g_k$, which is the graded subalgebra
of $W(J)$, generated by $J$ and $\mu$. We thus obtain a transitive $\Z$-graded
Lie superalgebra $Lie(J,\mu)$ such that the following two properties hold:
\begin{equation}
\g_k=\g_1^k ~~{\mbox for}~ k\geq 1,
\label{local}
\end{equation}
\begin{equation}
[\g_{-1},\mu] ~~{\mbox{generates\, the\, Lie\, superalgebra}}~ \g_0.
\label{subalgebra}
\end{equation}
As a result, the classification of algebra structures on $J$ is equivalent to 
the classification of $\Z$-graded transitive Lie superalgebras 
$\g=\prod_{k=-1}^\infty\g_k$ together with an element $\mu \in \fg_1$, 
satisfying (\ref{local}) and
(\ref{subalgebra}). We call this the {\em generalized TKK construction}.
The simplicity of the superalgebra $J$ is equivalent
to the following {\em irreducibility} of the grading:
\begin{equation}
{\mbox{the}}~ \g_0{-\mbox{module}}~ \g_{-1} ~{\mbox{is\, irreducible}}.
\label{irreducibility}
\end{equation}
Thus, the classification of simple algebra structures on $J$ is equivalent
to the classification of all transitive irreducible $\Z$-graded Lie
superalgebras (of depth 1) together with $\mu \in \g_1$, satisfying 
(\ref{local})
and (\ref{subalgebra}). Commutative (resp anti-commutative) rigid superalgebras
correspond to such $\Z$-graded superalgebras
with $\mu \in \g_1$, satisfying property (\ref{Str}), where $\mu$ is even
(resp.\ odd).
This reduces the classification
of simple linearly compact commutative (resp.\ anti-commutative)
rigid superalgebras to the classification
of linearly compact transitive $\Z$-graded Lie superalgebras 
$\g=\prod_{k=-1}^\infty \g_k$, satisfying (\ref{local}), (\ref{subalgebra}), 
(\ref{irreducibility}) and (\ref{Str}) for some even (resp.\ odd)
element $\mu \in \g_1$. Such a $\Z$-grading is called
{\em admissible}.

We obtain a classification of the latter in two steps. First, showing that
$\g$ is close to being simple and using the classification of simple
and semisimple linearly compact Lie superalgebras obtained in 
\cite{K3} and \cite{CantaK}, respectively, we obtain a classification
of all linearly compact Lie superalgebras which may have an admissible
grading.
Next, from the classification of all gradings, obtained in \cite{CantaK},
we extract the classification of all admissible gradings.

However this still does not complete the classification of linearly 
compact simple commutative (resp.\ anti-commutative) rigid
superalgebras since we still have to classify, for each
admissible grading, all even (resp.\ odd) $\mu\in\g_1$, 
satisfying (\ref{subalgebra}) and
(\ref{Str}). If $\dim \g_1<\infty$, property (\ref{Str})
implies that for the algebraic group $G$, whose Lie algebra is the even
part of $\g_0$, we have, if $\mu$ is even (resp.\ odd):
$G\cdot\mu$ is an open orbit in the even (resp.\ odd) part of $\g_1$.
Hence, up to isomorphism, there exists a unique $\mu$ in question.

If $\dim \g_1=\infty$, the situation is more complicated, since only a 
subalgebra of finite codimension $n$ of the even part of $\g_0$ can be 
integrated
to a group, which we again denote by $G$. Hence we have to classify the
orbits of $G$ in the even (resp.\ odd) part of $\g_1$ of codimension at
most $n$ (satisfying the additional property (\ref{subalgebra})).

The somewhat awkward notation for a rigid superalgebra $(J,\mu)$ is 
derived from that of the associated Lie superalgebra $\g=Lie(J,\mu)$ by
adding $J$ (resp.\ $L$) in the commutative (resp.\ anti-commutative) case 
in front of the notation for the type of $\g$, the first and the second
index in the notation being the growth and the size of $J$ (see Section 6
for their definition).

A posteriori it turns out that the even and odd parts of a simple 
commutative or anti-commutative rigid linearly compact superalgebra have
equal growth (resp.\ size), unless its odd part is zero.

In conclusion, we use the classification of odd type simple rigid commutative
linearly
compact superalgebras to derive the classification of simple linearly
compact odd generalized Poisson superalgebras: all of them turn out to
be gauge equivalent to $PO(n,n)$ or $PO(n,n+1)$.

Throughout the paper the base field $\F$ is algebraically closed of 
characteristic 0.

\section{Preliminaries on superalgebras}\label{pre}
\begin{definition}
A superalgebra (resp.\ odd type superalgebra) is a $\Z/2\Z$-graded vector space
$A=A_{\bar{0}}\oplus A_{\bar{1}}$, endowed with a product $a \otimes b\mapsto 
ab$, such that $ab\in A_{p(a)+p(b)}$ (resp.\ $A_{p(a)+p(b)+\bar{1}}$).
In this case the product is called even (resp.\ odd).
\end{definition}
\begin{definition}
A superalgebra or an odd type superalgebra, is called commutative
(resp.\ anti-commutative) if 
$ab=(-1)^{p(a)p(b)}ba$ (resp.\ $ab=-(-1)^{p(a)p(b)}ba$). 
\end{definition}

\begin{remark}\em 
Reversing parity exchanges superalgebras with odd
type superalgebras.
\end{remark}
\begin{remark}\em 
\label{rem:1.4}
Let $(R,\mu)$ be a commutative (resp.\ anti-commutative) superalgebra
with the product $\mu$. Define the new product $\bar{\mu}(a \otimes b)=
(-1)^{p(a)} \mu(a \otimes b)$ and let $\bar{R}$ be the vector superspace, 
obtained from $R$ by reversing the parity.
Then $(\bar{R},\bar{\mu})$ is an odd type anti-commutative (resp.\ commutative) 
superalgebra, and viceversa.
\end{remark}

A {\it linearly
compact superalgebra} is a topological superalgebra whose underlying
topological space is linearly compact. For a discussion on
linearly compact vector spaces and algebras we refer to \cite[\S 2]{CantaK2},
restricting ourselves here to a few remarks. Any
finite-dimensional vector space is linearly compact. The basic example of a
commutative associative linearly compact superalgebra is the superalgebra
${\cal O}(m,n)=\Lambda(n)[[x_1,\dots, x_m]]$ with the formal topology, 
where $\Lambda(n)$
denotes the Grassmann algebra over $\F$ on $n$ anti-commuting indeterminates
$\xi_1, \dots, \xi_n$, and the superalgebra parity is defined by
$p(x_i)=\bar{0}$, $p(\xi_j)=\bar{1}$.
The basic example of a
linearly compact Lie superalgebra is $W(m,n)=Der {\cal O}(m,n)$, the
Lie superalgebra of all continuous derivations of the superalgebra
${\cal O}(m,n)$. One has:
$$W(m,n):=\{X=\sum_{i=1}^mP_i(x,\xi)\frac{\partial}{\partial
x_i}+\sum_{j=1}^nQ_j(x,\xi)\frac{\partial}{\partial\xi_j} ~|~ P_i,
Q_j\in{\cal O}(m,n)\}.$$
Letting $a_i=\deg x_i=-\deg\frac{\partial}{\partial x_i}\in\Z_{>0}$ and
$b_j=\deg\xi_j=-\deg\frac{\partial}{\partial \xi_j}\in\Z$ we obtain
a $\Z$-grading of finite depth of the Lie superalgebra $W(m,n)$, called 
the grading of type 
$(a_1,\dots, a_m|b_1,\dots, b_n)$ 
(cf.\ \cite[Example 4.1]{K3}).
 When such a grading induces a $\Z$-grading
on a closed subalgebra $S$ of
$W(m,n)$, the induced grading on $S$ is also  called a grading of type 
$(a_1,\dots, a_m|b_1,\dots, b_n)$. The grading of $W(m,n)$ of type
$(1,\dots,1|1,\dots,1)$ is called the principal grading.

\medskip

The following is a generalization of $W(m,n)$ to the case of an infinite number 
of indeterminates.
Let $J$ be a vector superspace and 
let $W_k(J)$ be the superspace of $(k+1)$-linear supersymmetric
functions on $J$ with values in $J$, i.e., $f\in W_k(J)$ if
$f(\dots,x,y,\dots)=(-1)^{p(x)p(y)}f(\dots,y,x,\dots)$. Let
$W(J)=\prod_{k=-1}^\infty W_k(J)$. 
We endow this vector superspace  with a bracket as follows.
If $f\in W_p(J)$, $g\in W_q(J)$, then we define $f\Box g$ to be
the following element in $W_{p+q}(J)$:
$$f\Box g(x_0, \dots, x_{p+q})=
\sum_{\substack{
i_0<\dots <i_q\\
i_{q+1}<\dots< i_{p+q}}}\epsilon(i_0,\dots, i_q, i_{q+1},\dots,
i_{p+q})f(g(x_{i_0},\dots, x_{i_q}), x_{i_{q+1}},\dots, x_{i_{p+q}}).$$

\medskip

\noindent
Here $\epsilon=(-1)^N$,  where $N$ is the number of interchanges of indeces
of odd $x_i$'s in the permutation. 
\begin{proposition}
The bracket 
\begin{equation}
[f,g]=f\Box g-(-1)^{p(f)p(g)}g\Box f
\label{box}
\end{equation}
defines a Lie superalgebra structure on $W(J)$.
\end{proposition}
{\bf Proof.} It is clear that $f\Box g$ is supersymmetric if both $f$ and $g$
are. The anti-commutativity of bracket (\ref{box}) is immediate. 
Finally, it is easy to see that the
operation $\Box$ is right
associative, i.e.,  $(a,b,c)=(-1)^{p(b)p(c)}(a,c,b)$, where 
$(a,b,c)=(a\Box b)\Box c-a\Box(b\Box c)$. The right associativity
of $\Box$ implies the Jacobi identity for  bracket (\ref{box}).
\hfill$\Box$
\begin{remark}\em If $J$ is a  vector superspace of dimension $(m|n)$,
 then $W(J)$ is isomorphic to the Lie superalgebra $W(m,n)$
 \cite[\S 2.1]{K1}.
\end{remark}
\begin{remark}\label{rem:1.8}\em
According to the above definitions, $W_0(J)=End(J)$, acting on $W_{-1}(J)=J$ 
in the obvious way
and on $W_1(J)$ as follows: for $f\in W_0(J)$, $B\in W_1(J)$ and 
$x,y\in W_{-1}(J)$,
$$
[f,B](x,y)=f(B(x,y))-(-1)^{p(f)p(B)}(B(f(x),y)+(-1)^{p(x)p(y)}B(f(y),x)).
$$
Besides, for $A\in W_k(J)$ and $x,y_i\in W_{-1}(J)$, $i=1,...,k$, we have:
$[A,x](y_1,...,y_k)=A(x,y_1,...,y_k)$.
\end{remark}
\section{Odd generalized Poisson superalgebras}
\begin{definition}\label{oddPoisson} An odd generalized Poisson superalgebra
is a superalgebra $P$ with two operations: a unital commutative associative
superalgebra product $a\otimes b\mapsto ab$, and a Lie superalgebra bracket 
$a\otimes b\mapsto\{a,b\}$ with respect to reversed parity,
satisfying the generalized odd Leibniz rule, namely, for
$a,b,c\in P$ one has:
\begin{equation}
\{a,bc\}=\{a,b\}c+(-1)^{(p(a)+1)p(b)}b\{a,c\}+(-1)^{p(a)+1}D(a)bc
\label{oddPbracket}
\end{equation}
where $D(a)=\{e,a\}$, $e$ being the unit element.
If $D=0$, then relation (\ref{oddPbracket}) becomes the 
odd Leibniz rule; in this case $P$ is called an odd Poisson superalgebra.
(Note that $D$ is an odd derivation of the associative product and of the Lie 
superalgebra bracket.)
\end{definition}

\begin{example}\label{PO(n,n)}\em
Consider the associative superalgebra ${\cal O}(n,n)$
with the following bracket, known as the Buttin
bracket
($f,g\in {\cal O}(n,n)$):
\begin{equation}\{f,g\}_H=\sum_{i=1}^n(\frac{\partial f}{\partial x_i}\frac{\partial
g}{\partial \xi_i}+(-1)^{p(f)}\frac{\partial f}{\partial \xi_i}\frac{\partial
g}{\partial x_i}).
\label{buttinbracket}
\end{equation}
Then ${\cal O}(n,n)$ with this bracket is an odd Poisson superalgebra,
 which we denote  by $PO(n,n)$.
\end{example}

\begin{example}\label{PO(n,n+1)}\em
Consider the associative superalgebra ${\cal O}(n,n+1)$
with even indeterminates
$x_1, \dots, x_n$
 and odd indeterminates $\xi_1,\dots, \xi_n, \xi_{n+1}=\tau$.
Define on ${\cal O}(n, n+1)$ the following bracket
($f,g\in {\cal O}(n,n+1)$):
\begin{equation}
\{f,g\}_K=\{f,g\}_H+(E-2)(f)\frac{\partial g}{\partial \tau}
+(-1)^{p(f)}\frac{\partial
f}{\partial \tau}(E-2)(g),
\label{Leibniz!!}
\end{equation}
where $\{\cdot,\cdot\}_H$ is the Buttin bracket (\ref{buttinbracket})
and $E=\sum_{i=1}^n(x_i\frac{\partial}{\partial
x_i}+\xi_i
\frac{\partial}{\partial \xi_i})$ is the Euler operator.
Then ${\cal O}(n,n+1)$ with bracket $\{\cdot,\cdot\}_K$
is an odd  generalized Poisson superalgebra
with $D=-2\frac{\partial}{\partial \tau}$ \cite[Remark 4.1]{CantaK},
which we denote  by $PO(n,n+1)$.  

\medskip

Odd generalized Poisson superalgebras give rise to several series of simple 
linearly compact Lie superalgebras. Consider ${\cal O}(n,n)$ with its Lie 
superalgebra structure defined by
bracket (\ref{buttinbracket}). Then
$\F 1$ is a central ideal of this Lie superalgebra, hence
$HO(n,n):={\cal O}(n,n)/\F 1$, with reversed parity, is a linearly compact Lie 
superalgebra \cite[\S 1.3]{CK}. The Lie superalgebra $HO(n,n)$ is simple if 
and only if $n\geq 2$ (cf.\ \cite[Example 4.6]{K3}).

The Lie superalgebra $HO(n,n)$ contains the subalgebra $SHO'(n,n)=
\{f\in {\cal O}(n,n)/\F 1~|~\Delta(f)=0\}$ where $\Delta=\sum_{i=1}^n
\frac{{\partial}^2}{\partial x_i\partial\xi_i}$ is the odd Laplacian. 
For $n\geq 2$ the derived algebra of $SHO'(n,n)$ is an ideal of codimension
1, denoted by $SHO(n,n)$, consisting of elements not containing the
monomial $\xi_1\dots\xi_n$ \cite[\S 1.3]{CK}. $SHO(n,n)$ is simple
if and only if $n\geq 3$ \cite[Example 4.7]{K3}.

We shall denote by $KO(n,n+1)$ the superalgebra
$PO(n,n+1)$ with its Lie superalgebra structure and reversed parity.
This is a simple linearly compact Lie superalgebra for every $n\geq 1$
(cf.\ \cite[Example 4.8]{K3}). Since $KO(1,2)\cong W(1,1)$ (cf.\
\cite[Remark 6.6]{K3}), when dealing with the Lie superalgebra $KO(n,n+1)$ we shall
always assume $n\geq 2$.

For $\beta\in \F$, let $div_{\beta}=\Delta+(E-n\beta)\frac{\partial}{\partial
\tau}$. The Lie superalgebra $KO(n,n+1)$ contains the subalgebra
$SKO'(n,n+1;\beta)=\{f\in{\cal O}(n,n+1)~|~ div_{\beta}f=0\}$. Let 
$SKO(n,n+1;\beta)$ be the derived algebra of $SKO'(n,n+1;\beta)$. Then
$SKO(n,n+1;\beta)$ is simple for $n\geq 2$ and coincides with
$SKO'(n,n+1;\beta)$ unless $\beta=1$ or $\beta=\frac{n-2}{n}$.
The Lie superalgebra $SKO(n,n+1;1)$ (resp.\ $SKO(n,n+1;\frac{n-2}{n})$)
consists of the elements of $SKO'(n,n+1;1)$ (resp.\
$SKO'(n,n+1;\frac{n-2}{n})$) not containing the momomial
$\tau\xi_1\dots\xi_n$ (resp.\ $\xi_1\dots\xi_n$) \cite[\S 1.4]{CK}.  
\end{example}

\begin{example}\label{gauge}\em
Let $P$ be a unital odd generalized Poisson superalgebra with bracket 
$\{\cdot,\cdot\}$ and derivation $D$. Then for any even invertible $\varphi
\in P$, formula 
$$\{a,b\}^{\varphi}=\varphi^{-1}\{\varphi a,\varphi b\}$$
defines another odd generalized Poisson bracket on $P$ if and only if $\{\varphi,
\varphi\}=0$. In this case the derivation corresponding to
$\{\cdot,\cdot\}^{\varphi}$ is $D^{\varphi}(a)=-D(\varphi)a+\{\varphi, a\}$.
We denote this odd generalized Poisson superalgebra by $P^{\varphi}$. The
odd generalized Poisson superalgebras $P$ and $P^{\varphi}$ are called
{\em gauge equivalent}. 
\end{example}

\section{Rigid superalgebras}\label{jordansection}
Let $J$ be a vector superspace. Then products (resp.\ odd type products)
on $J$ that make it a superalgebra (resp.\ odd type superalgebra)
correspond bijectively to even (resp.\ odd) elements $\mu$
in the vector superspace 
$Hom(J\otimes J, J)$, 
so that $a\circ b =\mu(a\otimes b)$ 
is the product, corresponding to $\mu$.
Denote by $\mu_a \in End (J)$ the operator of left multiplication by $a\in J$ 
with respect to this product: $\mu_a(b)=a\circ b$.

\begin{definition}
The structural Lie superalgebra 
$Str(J,\mu)$ 
for the superalgebra $(J,\mu)$ 
is the subalgebra of the Lie superalgebra $End(J)_-$, 
generated by all $\mu_a$ with $a\in J$.
\end{definition}

Note that the vector superspace $Hom(J\otimes J, J)$ is naturally a module 
over the Lie superalgebra $(End J)_-$, hence over its subalgebra  
$Str(J,\mu)$. The action is defined as follows ($f\in (End J)_-$,
$B\in Hom(J\otimes J, J)$, $x,y\in J$):
\begin{equation}
[f,B](x,y)=f(B(x,y))-(-1)^{p(f)p(B)}(B(f(x),y)+(-1)^{p(x)p(f)}B(x, f(y))).
\label{action}
\end{equation}

\begin{definition}\label{related}
The space of related products $R(J,\mu)$ for the superalgebra $(J,\mu)$
is the minimal $Str(J,\mu)$-submodule of 
$Hom(J\otimes J, J)$,  
containing $\mu$.
\end{definition}

In the case $J$ is linearly compact 
we define $Str(J,\mu)$ and $R(J,\mu)$ as 
closures of the corresponding subspaces.

\begin{definition}\label{Jordan}
A vector superspace $J$ with a product (resp.\ odd type product)
$\mu$ is called a rigid superalgebra 
if $\mu$ is an even (resp.\ odd) element of $Hom(J\otimes J, J)$, such that  
\begin{displaymath}
R(J,\mu)=Str(J,\mu)(\mu)+\F \mu.
\end{displaymath}
\end{definition}

\begin{example}\em Any Jordan superalgebra is a commutative rigid 
superalgebra (see \cite[Remark 6]{K1}).
\end{example}

\begin{example}\em Any Lie superalgebra $\g$ is an anti-commutative rigid
superalgebra with $Str(\g,\mu)=\g$ and $R(\g,\mu)=\F \mu$.
\end{example}

Other examples of rigid commutative and anti-commutative superalgebras 
will be given in Sections 4 and 5 respectively (see also the introduction).


\begin{remark}\em 
\label{rem:3.6}
If $(J,\mu)$ is a commutative superalgebra, then $\mu$
is an even element of $W_1(J)=Hom(S^2(J),J)$,
and $R(J,\mu)$ is a $Str(J,\mu)$-submodule of $W_1(J)$.
If $(J,\mu)$ is an anti-commutative superalgebra, then
the odd type superalgebra 
$(\bar{J},\bar{\mu})$, defined in Remark
\ref{rem:1.4}, is commutative, $\bar{\mu}$ is an odd element of
$W_1(\bar{J})=Hom(S^2(\bar{J}),\bar{J})$, and $R(\bar{J},\bar{\mu})$ is a 
$Str(\bar{J},\bar{\mu})$-submodule of $W_1(\bar{J})$.
\end{remark}

\begin{proposition} \label{prop:3.7}
A superalgebra $(J,\mu)$ is rigid if and only if 
 $(\bar{J},\bar{\mu})$   is rigid.
\end{proposition}
{\bf Proof.} The map $(End J)_-\rightarrow (End \bar{J})_-$ defined by:
$f\mapsto \bar{f}$, $\bar{f}(\bar{a})=\overline{f(a)}$, $a\in J$,
 is an isomorphism of Lie superalgebras sending $Str(J,\mu)$ to
$Str(\bar{J},\bar{\mu})$. Now
consider the map $Hom(J\otimes J,J)\rightarrow Hom(\bar{J}\otimes \bar{J},\bar{J})$,
$\varphi\mapsto\bar{\varphi}$,
defined by: $\bar{\varphi}(\bar{a}\otimes\bar{b})=(-1)^{p(a)}
\overline{\varphi(a\otimes b)}$, $a,b\in J$.
This is an isomorphism between
the $End(J)_-$-module $Hom(J\otimes J, J)$ and
the $End(\bar{J})_-$-module  $Hom(\bar{J}\otimes \bar{J},\bar{J})$,
sending $R(J,\mu)$ to $R(\bar{J},\bar{\mu})$.
The statement then follows.
\hfill$\Box$

\section{Examples of commutative rigid  superalgebras}\label{evenexamples}
In the introduction we constructed five infinite-dimensional linearly compact
commutative superalgebras: $JS_{1,1}$,
$JSHO_{2,2}$, $JSKO_{1,2}$, and $JS^\alpha_{1,8}$ ($\alpha=0,1$), and a 
2-dimensional one, $JS_{0,2}$. In this section we shall construct the remaining 
four finite-dimensional commutative superalgebras, 
$JW_{0,4}$, $JW_{0,8}$, $JS_{0,8}$, and $JS_{0,16}$, that appear in the statement
of Theorem \ref{01}. We shall prove in
Section \ref{eight} that all these superalgebras 
are indeed simple rigid commutative superalgebras. The proof 
comes from the generalized TKK construction.

\begin{example}\label{is}\em 
Let $A$ be a unital commutative associative algebra (i.e.\ a superalgebra
whose underlying superspace is even) 
with a surjective derivation $D$. We define on $A$ the following two products:
$$(a)~~ f\circ g=D(f)D(g); ~~~~~~~~(b)~~ f\circ g=D(fg).$$
(Note that the map $D$ is a central extension homomorphism of
$(a)$ to $(b)$.)
Then $J=A$, 
with product $(a)$ or $(b)$, is a rigid commutative algebra. 
Indeed, in both cases $Str(J,\mu)=A D$ with the usual Lie
bracket, acting on $J$ in the obvious way in case $(a)$, and
as follows in case $(b)$: $[aD,x]=D(ax)$.
Besides, $R(J,\mu)=A\mu$, where, for $f\in A$, $x,y\in J$,
$(f\mu)(x,y)=fD(x)D(y)$ in case $(a)$, and
$(f\mu)(x,y)=D(fxy)$ in case $(b)$. In both cases
$Str(J,\mu)$ acts on $R(J,\mu)$ as follows:
$[aD,b\mu]=(aD(b)-2bD(a))\mu$.
It follows that $R(J,\mu)=[Str(J,\mu),\mu]$.
Note that, in case $(a)$, $\F 1$ is a proper ideal of $(J,\circ)$.
A special case of this example is $A=\F [[x]]$ and $D=d/dx='$ with the
product $f \circ g=f'g'$, called  the Beltrami product. This algebra is
a central extension of the simple algebra $\F [[x]]$ with the product
$f \circ g=(fg)'$, which we called the Beltrami algebra 
in the introduction.
In \cite{D} an interesting quartic identity in these algebras is derived.
The Beltrami product is the simplest case of one of the operations 
introduced by Beltrami in differential geometry.
\end{example}

\begin{remark}\em Let $(J,\circ)$ be a non-unital 
rigid superalgebra.
  Consider $\tilde{J}=J+\F e$ with the product defined by
  extending the product of $J$ by: $x\circ e=x$, for $x\in
  \tilde{J}$. Differently from what happens with usual non unital
  Jordan superalgebras \cite[\S 5]{CantaK2},
  $\tilde{J}$ is not necessarily a rigid
  superalgebra. Consider, for example, the central extension $J$ of 
the Beltrami  algebra 
  constructed in Example \ref{is}, and extend $'$ from
$J$ to $\tilde{J}$ by letting $e'=0$. Then
  $Str(\tilde{J},\mu)=\langle L_x, \varphi_x, \mu_e | x\in J\rangle$,
  where, for $z\in \tilde{J}$, $L_x(z)=x'z'$,  and $\varphi_x$ is defined as follows: $\varphi_x(z)=0$
  for $z\in J$ and $\varphi_x(e)=x$. It follows that
  $Str(\tilde{J},\mu).\mu+\F \mu=\langle [L_x,\mu], [\varphi_x,\mu],\mu\rangle$
  and one checks that the element $[\varphi_{\omega},[\varphi_x,\mu]]$,
  lying in $R(\tilde{J},\mu)$, for $\omega, x\in J$, does not belong to 
$Str(\tilde{J},\mu). \mu+\F
  \mu$. It follows that  $\tilde{J}$ does not satisfy Definition \ref{Jordan}.
\end{remark}


Given a superspace $A$, we shall denote, as before, by $\overline{A}$ a copy 
of $A$ with reversed parity, and given an element $a\in A$, we shall denote 
the corresponding element of $\overline{A}$  by $\overline{a}$. 


\begin{example}\label{JDD}\em
Let $A$ be a commutative associative superalgebra with 
two odd derivations $D_1$ and $D_2$. 
Let $J(A,D_1,D_2,\mu_1, \mu_2)=AD_1 + AD_2$  
and define, for $f,g\in A$, the following product (see the introduction):
\begin{equation}
\begin{array}{c}
fD_i \circ gD_i = [fD_i,gD_i]_+, \,i=1,2,\\
fD_1 \circ gD_2 = [fD_1,gD_2]_+
+(-1)^{p(g)}(\mu_1(f,g)D_1 + \mu_2(f,g)D_2),
\end{array}
\label{J(A,D_1,D_2)}
\end{equation}
where $\mu_1$ and $\mu_2$ are some odd type products on $A$.
Then $J(A,D_1,D_2, \mu_1, \mu_2)$ is a simple rigid
 superalgebra in the following cases:
\begin{itemize}
\item[a)] $A=\Lambda(\xi)$, $D_1=\frac{\partial}{\partial\xi}$ and $D_2=0$,
$\mu_1(f,g)=fg\xi=\mu_2(f,g)$;
\item[b)] $A=\Lambda(\xi_1, \xi_2)$, 
$D_1=\frac{\partial}{\partial\xi_1}$, $D_2=\frac{\partial}{\partial\xi_2}+\xi_1\xi_2\frac{\partial}{\partial\xi_1}$, 
$\mu_1(f,g)=fg\xi_1=\mu_2(f,g)$;
\item[c)]
$A=\Lambda(\xi_1, \xi_2)$,  $D_1=\frac{\partial}{\partial\xi_1}+\xi_1\xi_2\frac{\partial}{\partial\xi_2}$, 
$D_2=\frac{\partial}{\partial\xi_2}+\xi_1\xi_2\frac{\partial}{\partial\xi_1}-\xi_1\xi_2\frac{\partial}{\partial\xi_2}$,
$\mu_1(f,g)=fg(\xi_1+\xi_2)$, $\mu_2(f,g)=fg\xi_1$;
\item[d)] $A=\Lambda(\xi_1,\xi_2, \xi_3)$,
$D_1=\frac{\partial}{\partial\xi_1}+\xi_1\xi_2\frac{\partial}{\partial\xi_3}$,
$D_2=\frac{\partial}{\partial\xi_2}+\xi_1\xi_2\frac{\partial}{\partial\xi_1}+\xi_2\xi_3\frac{\partial}{\partial\xi_3}+\xi_1\xi_2\frac{\partial}{\partial\xi_3}$,
$\mu_1=0=\mu_2$;
\end{itemize}
We shall denote the superalgebras listed in a)-d) by $JW_{0,4}$,
$JW_{0,8}$, $JS_{0,8}$, and $JS_{0,16}$, respectively. 
\end{example}

\begin{remark}\em For every $\alpha\neq 0$, product (\ref{eq:0.9}) 
is isomorphic to (\ref{J(A,D_1,D_2)}) where
$D_1=\frac{\partial}{\partial\xi_1}+\xi_1\frac{\partial}{\partial x}+
\xi_2\frac{\partial}{\partial x}$, $D_2=\frac{\partial}{\partial\xi_2}$ and
$\mu_1=0=\mu_2$.
\end{remark}

In Section \ref{proofs}  we shall prove 
the following theorem.
\begin{theorem}\label{evenclass}
Any simple linearly compact commutative rigid  
superalgebra, which is not a Jordan superalgebra, 
is isomorphic either to one of the five infinite-dimensional superalgebras
$JS_{1,1}$, $JSHO_{2,2}$, $JSKO_{1,2}$, $JS^\alpha_{1,8}$ for $\alpha =0,1$,
or to one of the five finite-dimensional superalgebras $JS_{0,2}$,
 $JW_{0,4}$, $JW_{0,8}$, $JS_{0,8}$, $JS_{0,16}$. 
\end{theorem}

\section{Examples of anti-commutative rigid superalgebras}\label{oddexamples}
We shall prove in Section \ref{eight} that the series $LP(n,n)$ and
$LP(n,n+1)$ from Example \ref{mainexample}, as well as examples 
\ref{SHO}--\ref{OJW},  \ref{OJS_{1,3}}--\ref{OJS}, \ref{JSKO0}--\ref{OJH},
\ref{HO(3,3)},
 \ref{OJSHO(4,4)}--\ref{OJSKO(3,4;beta)} 
in this section are indeed simple rigid anti-commutative superalgebras.
The proof  comes from the generalized TKK construction.
With the exception of the series $LP(n,n)$ and $LP(n,n+1)$,
we shall denote the anti-commutative simple rigid 
superalgebras  by $LX_{m,n}$, where the first index 
is the growth of the superalgebra 
and the second is its total size (equal to the total 
dimension if the size is 0). The letter $X$ comes
from the connection to the Lie superalgebra $X(m',n')$ via the generalized 
TKK construction. 
\begin{example}\label{mainexample}\em 
%
%
%
%
%
Let $P$ be an odd generalized Poisson superalgebra 
with the bracket $\{\cdot,\cdot \}$ and the derivation $D$ (defined by $D(a)=
\{e,a\}$), 
let $'=d/dx$ on $P[[x]]$ and let
$\eta$ be an odd indeterminate, such that $\eta^2=0$. 
Consider the commutative associative superalgebra $J=P[[x]]+\eta P[[x]]$
and extend on $J$, by commutativity, the following odd product 
($f,g\in P, ~a,b\in \F[[x]]$):

\begin{equation}
\begin{array}{l}
fa\circ gb=(-1)^{p(f)+1}\{f,g\}ab+
\frac{1}{2}(-1)^{p(f)}fD(g)xa'b+
\frac{1}{2}D(f)gxb'a+2\eta fgab\\
\\
fa\circ \eta(gb)=\eta\{f,g\}ab-
((-1)^{p(f)}fg+\frac{1}{2}\eta fD(g)x)a'b
-\frac{1}{2}(-1)^{p(f)}\eta D(f)ga(b+xb')\\
\\
\eta(fa)\circ \eta(gb)=(-1)^{p(f)}\eta fg(a'b-ab').\\
\end{array}
\label{goddP}
\end{equation}
We will denote the resulting odd type superalgebra by $OJP$. 
If $P=PO(n,n)$ (resp.\ $P=PO(n,n+1)$), then we will denote it by $OJP(n,n)$ 
(resp.\ $OJP(n,n+1)$). 

Now we shall show that $OJP$ 
is a rigid odd type superalgebra. 
Denote by $\mu_v$ the operator of left multiplication by $v\in J$ 
with respect to the product $\circ$, defined by (\ref{goddP}), 
and by $l_v$ the operator of left multiplication by $v$ with respect to the
associative product on $J$. Then 
$Str(J)=\langle \mu_v, l_v ~|~ v\in J\rangle$. Indeed,
for $v,w\in J$, we have:
$$
[\mu_v, \mu_w]=(-1)^{p(v)+1}(\mu_{v\circ w-2\eta vw}+2l_{\eta(v\circ w)+(vw)'
+2\eta D(vw)}).
$$ 
If we take $v=fa$ and $w=\eta gb$, for $f, g\in P$ and $a,b\in\F[[x]]$, then
$\eta(v\circ w)+(vw)'
+2\eta D(vw)=\eta fg ab'$, therefore $Str(J)$ contains all $l_z$ for
$z\in \eta P[[x]]$. If we take  $v=fa$ and $w=gb$, then $(vw)'=fg(ab)'$, hence
$Str(J)$ contains also all $l_z$ for $z\in P[[x]]$.
Finally we have:
$$
[\mu_{v}, l_{w}]=l_{v\circ w-2\eta vw}+l_{D(v)w}, \,\, [l_v,l_w]=0.
$$
Next, the following commutation relations hold:
\begin{equation}
[l_v,[l_w,\mu]]=-[l_{vw},\mu]+(-1)^{p(v)}[\mu_{z-2\eta\psi}-2l_{\eta z},\mu]
\label{KV}
\end{equation}
for some $z, \psi\in J$ such that $z'=(-1)^{p(v)+1}(v\circ w-2\eta vw+D(vw))$,
$\psi'=D(z)$;
\begin{equation}
[l_v,[\mu_w,\mu]]=(-1)^{p(v)+1}[\mu_{vw},\mu]+2(-1)^{p(v)+1}[\mu_{z-2\eta\psi}
-2l_{\eta z},\mu],
\label{KVI}
\end{equation}
for some $z, \psi\in J$ such that $z'=v'w+2\eta D(vw)$, $\psi'=D(z)$;
\begin{equation}
[\mu_w,[l_v,\mu]]=[l_{w\circ v-2\eta wv},\mu]+[l_{D(w)v},\mu]+
(-1)^{p(v)(p(w)+1)}[l_v,[\mu_w,\mu]];
\label{KVII}
\end{equation}
\begin{equation}
[\mu_v,[\mu_w,\mu]]=2[l_{v\circ \eta w},\mu]+ (-1)^{p(v)+1}([\mu_{D(v)w},\mu]-
2[\mu_{\eta z},\mu]+[\mu_{\varphi-2\eta\psi}-2l_{\eta\varphi},\mu])
\label{KVIII}
\end{equation}
for some $z,\varphi, \psi\in J$ such that 
$z'=vw'-v'w-\{v,D(w)\}+2(-1)^{p(v)}D(v)D(w)$,  
$\varphi'= 2\eta vw'-v\circ w'-D(v)w'-2(-1)^{p(v)+1}v'D(w)-D(v')w$, 
$\psi'=D(\varphi)$.
Relations (\ref{KV})--(\ref{KVIII}) show that
$[Str(J,\mu),\mu]+\F \mu=R(J,\mu)$, hence $OJP$
is a rigid odd type superalgebra.
Notice that relations (\ref{KV})--(\ref{KVIII}) do not depend on the choice 
of the functions $z$ and $\psi$, since,
for every $c\in \F$,
$\mu_c=2l_{\eta c}-l_c\circ D$, as it can be deduced from
(\ref{goddP}), $[l_c\circ D, \mu]=0$, and $[\mu_{\eta c},\mu]=0$.

Let $LP=\overline{OJP}$ with product $\overline{f}\bullet 
\overline{g}=(-1)^{p(f)}\overline{f\circ g}$. By Proposition \ref{prop:3.7},
$LP$ is an anti-commutative rigid superalgebra. If $P=PO(n,n)$ (resp.\
$P=PO(n,n+1)$), then we will denote  the corresponding anti-commutative rigid 
superalgebra by $LP(n,n)$ (resp.\ $LP(n,n+1)$).
\end{example}

\begin{remark}\label{reconstruct}\em 
The generalized odd Poisson bracket on $P$ 
can be recovered  from
product (\ref{goddP}) as follows ($f,g\in P$):
$\{f,g\}=(-1)^{p(f)+1}f\circ g+2\eta(fx)\circ\eta g$.
Furthermore, $(-1)^{p(f)}\eta fx\circ\eta g=\eta fg$, hence 
also the associative
product on $P$ can be recovered from product (\ref{goddP}).
\end{remark}

\begin{remark}\label{simplier}\em 
The odd Poisson superalgebra $PO(n+1,n+1)$ can be decomposed
as follows: $PO(n+1,n+1)=PO(n,n)[[x_{n+1}]]+\xi_{n+1}PO(n,n)[[x_{n+1}]]$.
Therefore, if $P=PO(n,n)$, then $LP(n,n)=PO(n+1,n+1)$ with reversed parity and
the following product:
$$f\bullet g=-\{f,g\}_H+2(-1)^{p(f)}\xi_{n+1}fg$$
where $\{\cdot,\cdot\}_H$ is the Buttin bracket
on $PO(n+1,n+1)$.

Likewise, $PO(n+1,n+2)=PO(n,n+1)[[x_{n+1}]]+\xi_{n+1}PO(n,n+1)[[x_{n+1}]]$
and if $P=PO(n,n+1)$, then $LP(n,n+1)=PO(n+1,n+2)$ with reversed parity and
 product
$$f\bullet g=-\{f,g\}_K+2(-1)^{p(f)}\xi_{n+1}fg$$
where $\{\cdot,\cdot\}_K$ is the odd generalized Poisson bracket on 
$PO(n+1,n+2)$.
\end{remark}

\begin{example}\label{SHO}\em 
Let us define on
${\cal O}={\cal O}(x_1,\dots,x_n,\xi_1,\dots,\xi_n)$, with $n\geq 2$,
 the following product
($f,g\in {\cal O}$):
\begin{equation}
f\bullet g=\{f,g\}+2(-1)^{p(f)+1}\xi_1fg+2\{x_2\xi_1\xi_2f,g\}
+2(-1)^{p(f)}\{x_2\xi_1\xi_2,f\}g
\label{productS1}
\end{equation}
where $\{\cdot,\cdot\}$ is the usual odd Poisson bracket on 
${\cal O}(x_1,\dots,x_n,\xi_1,\dots,\xi_n)$.

Let
$J'=\{f\in {\cal O} ~|~ \Delta(f)=0\}$,
where $\Delta$ is the odd Laplacian, i.e., $\Delta=\sum_{i=1}^n
\frac{\partial^2}{\partial x_i\partial\xi_i}$. Then $J'$ is closed
under product (\ref{productS1}). Besides,
the span $J$ of all monomials in $J'$ except for 
$\xi_1\dots\xi_n$ is an ideal of $(J',\bullet)$.

$J$ with reversed parity and product (\ref{productS1}) 
is a simple rigid anti-commutative
superalgebra that we will denote by $LSHO_{n,2^{n-1}}$.

\end{example}


\begin{example}\label{OJSKO}\em 
Let us define on
${\cal O}={\cal O}(x_1,\dots,x_n,\xi_1,\dots,\xi_n,\tau)$, for $n\geq 1$,
 the following product
($f,g\in {\cal O}$):
\begin{equation}
f\bullet g=\{f,g\}_K+\{\xi_1\tau f,g\}_K+(-1)^{p(f)+1}
(\xi_1(2E-(n+1)\beta)(f)g+\{\xi_1\tau,f\}_K g-2\tau\frac{\partial f}{\partial x_1}g)
\label{productSK1}
\end{equation}
where $\{\cdot,\cdot\}_K$ is the usual generalized odd Poisson bracket on 
${\cal O}(x_1,\dots,x_n,\xi_1,\dots,\xi_n,\tau)$,
$E$ is the Euler operator, and $\beta\in\F$.
Let
$J'=\{f\in {\cal O} ~|~ div_\beta(f)+(1-\beta)\frac{\partial f}{\partial\tau}=0\}$,
where $div_\beta$ is the $\beta$-divergence introduced in Example
\ref{PO(n,n+1)}. Then $J'$ is closed
under product (\ref{productSK1}). Besides, if $\beta=1$ (resp.\ 
$\beta=\frac{n-1}{n+1}$), then
the span $J_\beta$ of all monomials in $J'$ except for 
$\xi_1\dots\xi_n\tau$ (resp.\ $\xi_1\dots\xi_n$) is an ideal of $(J',\bullet)$.
For $\beta\neq 1,\frac{n-1}{n+1}$, let us set $J_{\beta}=J'$. Then,
for every $\beta\neq \frac{4}{n+1}$,
$J_{\beta}$ with reversed parity and product $\bullet$ is a
simple rigid anti-commutative superalgebra that we shall denote by
$LSKO_{n,2^n}(\beta)$.
\end{example}

\begin{example}\label{LSKO'}\em
Let $n=2$, $\beta=1$ and $J_{\beta}$ be as in Example \ref{OJSKO}. Consider
the following product on $J_1$:
$$f\bullet g=\{f,g\}_K+\{\xi_2(\tau-\xi_1)f,g\}_K-(-1)^{p(f)+1}(\xi_2fg-
\{\xi_2(\tau-\xi_1),f\}g-\frac{\partial f}{\partial \tau}\xi_1\xi_2g).$$
Then $J_1$ with reversed parity and product $\bullet$ is a simple rigid 
anti-commutative superalgebra that we will denote
by $LSKO'_{2,4}$.
\end{example}

%
%
%
%
%

\begin{example}\em 
Let $LW_{0,2}=\F a+\F\bar{a}$, where $a$ is even and $\bar{a}$ is odd,
 with product
$$a\bullet a=0, ~~ a\bullet\bar{a}=-\bar{a}\bullet a=\bar{a}, ~~\bar{a}\bullet\bar{a}=a.$$
\end{example}

\begin{example}\em
Let $LW_{1,2}=\F[[x]]+\overline{\F[[x]]}$ with product
$$f\bullet g=0, ~~ f\bullet\bar{g}=-\bar{g}\bullet f=\overline{fg}, ~~\bar{f}\bullet\bar{g}=
2fg-\frac{d(fg)}{dx}.$$
\end{example}

\begin{example}\label{OJHO_{1,2}}\em 
Let $LHO_{1,2}=\F[[x]]/\F 1+\overline{\F[[x]]}$ with product ($f,g\in \F[[x]]$):
\begin{equation}
f\bullet g=0, ~~~~ f\bullet \overline{g}=-\overline{g}\bullet f=\overline{\frac{df}{dx}g}, ~~~~
\overline{f}\bullet\overline{g}=2fg.
\end{equation}
\end{example}

\begin{example}\label{SHO(3,3)}\em
Let 
$LSHO'_{2,2}=\F[[x_1, x_2]]/\F 1+\overline{\F[[x_1, x_2]]}$
with product ($f,g\in\F[[x_1, x_2]]$):
$$f\bullet g=-\{f,g\}, ~~ f\bullet\bar{g}=-\bar{g}\bullet f=\overline{\{f,g\}}+
\overline{gD_2(f)}, ~~ \bar{f}\bullet\bar{g}=-2fg$$
where $\{f,g\}=D_1(f)D_2(g)-D_2(f)D_1(g)$,
$D_1=(1+x_1)\frac{\partial}{\partial x_1}$
and $D_2=\frac{\partial}{\partial x_2}$.
\end{example}

\begin{example}\label{OJW}\em
Let 
$LW^\alpha_{1,2}=\F[[x]]\frac{d}{dx}+\F[[x]]$ with
product ($f,g\in \F[[x]]$, $\alpha\in\F$):
\begin{equation}
f\frac{d}{dx}\bullet g\frac{d}{dx}=\left[f\frac{d}{dx}, g\frac{d}{dx}\right]_-,
~~~~ f\frac{d}{dx}\bullet g=-g\bullet f\frac{d}{dx}
=-(\alpha+x)(fg)\frac{d}{dx}+f\frac{dg}{dx},~~~~
f\bullet g=0.
\label{W(1,2)}
\end{equation}
\end{example}

\begin{remark}\em For every $\alpha\neq 0$, product (\ref{W(1,2)}) is
isomorphic to the following ($f,g\in \F[[x]]$, $\alpha\in\F$):
\begin{equation}
f\frac{d}{dx}\bullet g\frac{d}{dx}=\left[f\frac{d}{dx}, g\frac{d}{dx}\right]_-,
~~~ f\frac{d}{dx}\bullet g=-g\bullet f\frac{d}{dx}
=-(fg)\frac{d}{dx}+f\frac{dg}{dx},~~~~
f\bullet g=0.
\end{equation}
\end{remark}

\begin{example}\label{OJS_{1,3}}\em
Let $LS_{1,3}=\F[[x]]\frac{d}{dx}+\F[[x]]+\widetilde{\F[[x]]}$ with
the following product ($f,g\in\F[[x]]$):
$$f\frac{d}{dx}\bullet g\frac{d}{dx}=\left[f\frac{d}{dx}, g\frac{d}{dx}
\right]_-,~~~
f\bullet g=0,~~~\tilde{f}\bullet\tilde{g}=0,$$
$$f\frac{d}{dx}\bullet g=f\frac{dg}{dx}, ~~~~f\frac{d}{dx}\bullet\tilde{g}=
\widetilde{f\frac{dg}{dx}},~~~~
f\bullet\tilde{g}=fg\frac{d}{dx},$$
extended to all other pairs of elements of $LS_{1,3}$ by anti-commutativity.
\end{example}

\begin{example}\em
Let $LW^\alpha_{2,2}=\F[[x_1,x_2]]\frac{\partial}{\partial x_1}
+\F[[x_1,x_2]]\frac{\partial}{\partial x_2}$ with
the following product ($f,g\in \F[[x_1,x_2]]$, $\alpha\in\F$):
\begin{equation}
\begin{array}{c}
f\frac{\partial}{\partial x_i}\bullet g\frac{\partial}{\partial x_i}=
\left[f\frac{\partial}{\partial x_i}, g\frac{\partial}{\partial x_i}\right]_-,~~~i=1,2\\
f\frac{\partial}{\partial x_1}\bullet g\frac{\partial}{\partial x_2}=
-g\frac{\partial}{\partial x_2}\bullet f\frac{\partial}{\partial x_1}=
[f\frac{\partial }{\partial x_1}, g\frac{\partial }{\partial x_2}]_- 
+fg(1+x_1)\frac{\partial}{\partial x_1}
-\alpha x_2fg\frac{\partial}{\partial x_2}.
\end{array}
\end{equation}
\end{example}

\begin{example}\label{OJS}\em
Let $LS^\alpha_{2,2}=\F[[x_1,x_2]]D_1+\F[[x_1,x_2]]D_2$, where
$D_1=\frac{\partial}{\partial x_1}$ and $D_2=\frac{\partial}{\partial x_2}+
\alpha x_1\frac{\partial}{\partial x_2}+x_1x_2\frac{\partial}{\partial x_2}$,
$\alpha\in\F$, with
the following product ($f,g\in \F[[x_1,x_2]]$):
\begin{equation}
\begin{array}{l}
fD_i\bullet gD_i=\left[fD_i, gD_i\right]_-, ~~i=1,2\\
fD_1\bullet gD_2=-gD_2\bullet fD_1=\left[fD_1, gD_2\right]_-+x_1fgD_1.
\end{array}
\label{OJS(2,2)}
\end{equation}
\end{example}

\begin{remark}\em
For every $\alpha\neq 0$, $LS^\alpha_{2,2}$  is isomorphic to the
anti-commutative superalgebra $J'=\F[[x_1,x_2]]D'_1+\F[[x_1,x_2]]D'_2$,
where $D'_1=D_1$, $D'_2=\frac{\partial}{\partial x_2}+
x_1\frac{\partial}{\partial x_2}$ and the product is defined as follows
($f,g\in \F[[x_1,x_2]]$):
\begin{equation}
\begin{array}{l}
fD'_i\bullet gD'_i=\left[fD'_i, gD'_i\right]_-, ~~i=1,2\\
fD'_1\bullet gD'_2=-gD'_2\bullet fD'_1=\left[fD'_1, gD'_2\right]_-.
\end{array}
\end{equation}
\end{remark}

\begin{example}\label{JSKO0}\em
Let $J=\F[[x]]\frac{d}{dx}+\overline{\F[[x]]}$. 
Consider the following product ($f,g\in \F[[x]]$, $\beta\in\F$):
$$f\frac{d}{dx}\bullet g\frac{d}{dx}=-\beta\left[f\frac{d}{dx}, g\frac{d}{dx}\right]_-, 
~~f\frac{d}{dx}\bullet\overline{g}=-\overline{g}\bullet f\frac{d}{dx}=
-\overline{g\frac{df}{dx}}+\beta\overline{f\frac{dg}{dx}},
~~\overline{f}\bullet\overline{g}=xfg\frac{d}{dx}.$$
Then, for every $\beta\neq 0,1,2+\frac{1}{b},2+\frac{2}{b}$, $b\in\Z_{>0}$,
 $(J,\bullet)$ is a simple rigid anti-commutative superalgebra that we shall denote by 
$LSKO'_{1,2}(\beta)$.
\end{example}

\begin{example}\label{OJH}\em 
Let $LH^\alpha_{2,2}=\F[[x]]\frac{d}{dx}+\overline{\F[[x]]}/\F\bar{1}$
with the following product
($f,g\in \F[[x]]$, $\alpha\in\F$):
\begin{equation}
f\frac{d}{dx}\bullet g\frac{d}{dx}=\left[f\frac{d}{dx}, g\frac{d}{dx}\right]_-,
~~~
f\frac{d}{dx}\bullet\bar{g}=-\bar{g}\bullet f\frac{d}{dx}=\overline{f\frac{dg}{dx}},~~~
\bar{f}\bullet\bar{g}=-2(\alpha+x)\frac{df}{dx}\frac{dg}{dx}.
\label{OJH(2,2)}
\end{equation}
\end{example}

\begin{remark}\em
For every $\alpha\neq 0$, product (\ref{OJH(2,2)}) is isomorphic 
to the following ($f,g\in \F[[x]]$, $\alpha\in\F$):
$$
f\frac{d}{dx}\bullet g\frac{d}{dx}=\left[f\frac{d}{dx}, g\frac{d}{dx}\right]_-,~~~
f\frac{d}{dx}\bullet\bar{g}=-\bar{g}\bullet f\frac{d}{dx}=\overline{f\frac{dg}{dx}}, ~~~
\bar{f}\bullet\bar{g}=-2\frac{df}{dx}\frac{dg}{dx}.
$$
\end{remark}

\medskip

Let $A$ be a commutative associative unital algebra. We call a generalized 
bivector field any finite sum of the form $Z+\sum_iX_i\wedge Y_i$, 
where $Z, X_i, Y_i$ are derivations of $A$. The corresponding quasi Poisson 
bracket on $A$ is:
$$\{f,g\}=Z(f)g - fZ(g)+\sum_iX_i(f)Y_i(g)-Y_i(f)X_i(g).$$
This is a skew-symmetric product which in general does not satisfy the Jacobi 
identity.
The following examples of anti-commutative rigid algebras come from 
quasi Poisson brackets.

\begin{example}\label{HO(3,3)}\em
$LHO^\alpha_{3,1}=\F[[x_1,x_2, x_3]]/\F 1$ with the quasi Poisson bracket
associated with the bivector field
$\frac{\partial}{\partial x_1}\wedge\frac{\partial}{\partial x_2}+
(\alpha x_1+x_1^2)
\frac{\partial}{\partial x_1}\wedge\frac{\partial}{\partial x_3}$,
$\alpha\in\F$. We call this the Jacobi-Mayer algebra.
\end{example}

\begin{remark}\em If $\alpha\neq 0$ in Example \ref{HO(3,3)},
we can replace $\alpha x_1+x_1^2$ by $x_1$, getting an isomorphic 
algebra. This algebra appears in the theory of PDEs under the name 
Jacobi-Mayer bracket (see \cite {D2}).
\end{remark}

\begin{example}\label{OJSHO(4,4)}\em
$LSHO^\alpha_{4,1}=\F[[x_1, x_2, x_3, x_4]]/\F 1$
with the quasi Poisson bracket associated with the bivector field 
$\frac{\partial}{\partial x_1}\wedge\frac{\partial}{\partial x_2}+
(\alpha + x_1)
\frac{\partial}{\partial x_3}\wedge\frac{\partial}{\partial x_4}$, 
$\alpha\in\F$.
\end{example}

\begin{example}\em $LKO_{2,1}=\F[[x_1,x_2]]$ with the quasi Poisson bracket
associated with the generalized bivector field
$2\frac{\partial}{\partial x_1}+(x_1+x_2)\frac{\partial}{\partial x_1}
\wedge \frac{\partial}{\partial x_2}$.
\end{example}

\begin{example}\label{OJSKO(3,4;beta)}\em
Let $LSKO^\alpha_{3,1}(\beta)=\F[[x_1, x_2, x_3]]$ with the quasi Poisson
bracket associated with the generalized bivector field
$\frac{\partial}{\partial x_1}
\wedge \frac{\partial}{\partial x_2}-3\beta(\alpha+ x_1)
x_1\frac{\partial}{\partial x_1}
\wedge \frac{\partial}{\partial x_3}-(\alpha+2x_1)(x_2\frac{\partial}{\partial 
x_2}\wedge\frac{\partial}
{\partial x_3}-2\frac{\partial}
{\partial x_3})$, where $\alpha, \beta\in\F$.
Then, for every $\alpha\neq 0$ and $\beta\neq\frac{1}{3}$ 
and for $\alpha=0$, $LSKO^\alpha_{3,1}(\beta)$ is a simple rigid 
anti-commutative superalgebra.
\end{example}



In Section \ref{proofs}  we shall prove the following theorem.

\begin{theorem}\label{oddclass} Any simple linearly compact anti-commutative
rigid superalgebra, which is not a Lie superalgebra,
is isomorphic to one of the following superalgebras:
\begin{itemize}
\item[-] $LP(n,n)$ or $LP(n,n+1)$ for $n\geq 0$;
\item[-] $LSKO'_{2,4}$, 
$LSHO_{n,2^{n-1}}$ or $LSKO_{n-1,2^{n-1}}(\beta)$  for $n\geq 2$
and for $\beta\neq \frac{4}{n}$;
\item[-] $LW_{k,2}$ for $k=0,1$, $LHO_{1,2}$, $LSHO'_{2,2}$; 
\item[-] $LW^\alpha_{1,2}$ for $\alpha =0,1$;
\item[-] $LW^\alpha_{2,2}$, $LS^\alpha_{2,2}$ for $\alpha =0,1$;
\item[-] $LS_{1,3}$;
\item[-] $LH^\alpha_{1,2}$ for $\alpha =0,1$, $LSKO'_{1,2}(\beta)$ for $\beta\neq 0,1,2+\frac{1}{b},
2+\frac{2}{b}, b\in\Z_{>0}$;
\item[-] $LKO_{2,1}$, $LHO^\alpha_{3,1}$, $LSHO^\alpha_{4,1}$,
 $LSKO^\alpha_{3,1}(\beta)$ for $\alpha =0$ or $\alpha=1$ and $\beta\neq \frac{1}{3}$.
\end{itemize}
\end{theorem}

\section{Admissible gradings}
\begin{definition}\label{admissible}
An even admissible (resp.\ odd admissible) grading   of a Lie superalgebra
$\g$ is a transitive $\Z$-grading $\g=\prod_{j\geq -1}\g_j$ such that there 
exists an even (resp.\ odd)
element $\mu\in\g_1$ satisfying the following properties:
\begin{equation}
[\g_0,\mu]+\F \mu=\g_1;
\label{a}
\end{equation}
\begin{equation}
\mbox{the subspace}~ [\g_{-1},\mu] ~\mbox{generates the Lie superalgebra}~ 
\g_0;
\label{b}
\end{equation}
\begin{equation}
\g_k=\g_1^k ~\mbox{for every}~ k\geq 1.
\label{c}
\end{equation}
\end{definition}

In \cite[\S 1]{CantaK} the notions of growth and size of an artinian linearly
compact Lie superalgebra $\g$, related
to a Weisfeiler filtration of $\g$, were defined. We define the
growth and size of a subspace of $\g$ as those 
for the induced filtration from a Weisfeiler filtration of $\g$. 
This definition is independent of the choice of the filtration due to the
equivalence of Weisfeiler filtrations for any artinian $\g$. 

\begin{remark}\label{necessarycond}\em
If $\g=\prod_{j\geq -1}\g_j$ is an even admissible or odd admissible grading
of a Lie superalgebra $\g$ with $\dim \g_0 < \infty$, then, by 
condition (\ref{a}), 
$\dim\g_1\leq\dim\g_0+1$.
If $\g_0$ is infinite-dimensional and $\g$ is artinian linearly compact,
then, by condition (\ref{a}), 
$\growth \g_1 \leq \growth\g_0$ and $\size \g_1 \leq \size \g_0$.
\end{remark}

\begin{proposition}\label{oddJordan}
Let $\g=\prod_{j\geq{-1}}\g_j$ be an even admissible (resp.\ odd admissible)
 grading of a
Lie superalgebra $\g$, and let $\mu\in(\g_1)_{\bar{0}}$
(resp.\ $\mu\in(\g_1)_{\bar{1}}$) be as in Definition
\ref{admissible}. Define the following product on $J=\g_{-1}$:
\begin{equation}
x\circ y=[[\mu,x],y] ~~~~~~(x,y\in J).
\label{basicdef}
\end{equation}
Then $(J,\circ)$ is a 
 rigid (resp.\ rigid odd type) superalgebra 
that we will denote by $J(\g,\mu)$.
\end{proposition}
{\bf Proof.} By definition, for $x\in\g_{-1}$, 
$\mu_x=ad\, [\mu,x]$, hence,
by property (\ref{b}), $Str(J,\mu)=\g_0$ acting on $J=\g_{-1}$
via the adjoint action.
Therefore, by property
(\ref{a}), $[Str(J,\mu),\mu]+\F \mu=\g_1$, hence $R(J,\mu)=\g_1$
and Definition \ref{Jordan}
is satisfied. \hfill$\Box$

\begin{remark}\label{evenshort}\em
Let $\g=\prod_{i\geq -1} \g_i$ be an even admissible grading with
$\mu\in(\g_1)_{\0}$ as in Definition \ref{admissible}. The grading is 
{\em short},
i.e.\ $\g_i = 0$ for $i> 1$, if and only if $J(\g,\mu)$ is a Jordan 
superalgebra
\cite[Proposition 5.1, Lemma 5.13]{CantaK2}.
\end{remark}

\begin{proposition}\label{[p,p]=0} Let $\g=\prod_{j\geq{-1}}\g_j$ be an 
odd admissible grading of a
Lie superalgebra $\g$, and let 
$\mu\in(\g_1)_{\bar{1}}$ be an element as in Definition
\ref{admissible}, such that $[\mu,\mu]=0$. Then $ \g=\g_{-1}+\g_0+\g_1\cong\xi\g_0+\g_0+\F \frac{d}{d\xi}$, for some odd indeterminate $\xi$.
\end{proposition}
{\bf Proof.}  Since  $[\mu,\mu]=0$,  by the Jacobi identity one has
$(ad \,\mu)^2=0$.
Now consider  product (\ref{basicdef}) on $\g_{-1}$. Condition
$(ad \,\mu)^2=0$ implies the following identity
$$(a\circ b)\circ c=(-1)^{p(a)+1}a\circ (b\circ c)-(-1)^{p(b)(p(a)+1)}
b\circ (a\circ c).$$
It follows that, if we set $[a,b]'=(-1)^{p(a)}a\circ b$, then
$[\cdot,\cdot]'$ 
defines
on $\g_{-1}$ with reversed parity a Lie superalgebra bracket.
By condition (\ref{b}), $\g_0$ is generated by the left translations
$\mu_x=[\mu,x]$ with $x\in\g_{-1}$, and the map $x\mapsto-\mu_x$,
from $(\g_{-1},[\cdot,\cdot]')$ with reversed parity to $\g_0$,
defines an isomorphism of Lie superalgebras since $(ad\,\mu)^2=0$. 
Besides, 
$(ad\,\mu)^2=0$ implies $[\mu,\mu_x]=0$ for every $x\in \g_{-1}$,
hence, by condition 
(\ref{a}), $\g_1=\F \mu$.
\hfill$\Box$

\begin{remark}\label{short}\em
Let $\g=\oplus_{j= -1}^d\g_j$ be an even admissible grading
of a Lie superalgebra $\g$, with $\mu\in (\g_1)_{\bar{0}}$
as in Definition \ref{admissible}. Suppose that $J=\g_{-1}$ with
product (\ref{basicdef}) is a rigid superalgebra with unit element
$e$. Set $h=\mu_e$, then $\{e, h, \mu\}$ is an $sl_2$ triple,
whose adjoint action on $\g$ exponentiates to the action of $SL_2$
by continuous
automorphisms on $\g$,
since the grading on $\g$ has finite height.  Therefore we have an
automorphism of $\g$ which exchanges $\g_{-1}$ with $\g_1$, 
hence the grading is short, i.e., $\g_i =0$ for $i>1$, and therefore $J$ is a 
Jordan superalgebra.
\end{remark}

\section{Generalized TKK construction}
\begin{definition} Let $J$ be a commutative superalgebra with
product $\mu$ (even or odd). 
We denote by $Lie(J,\mu)$ the $\Z$-graded Lie subalgebra
of $W(J)$ generated by $J=W_{-1}(J)$ and $\mu$. 
\end{definition}

\begin{proposition}\label{irr} Let $J$ be a commutative superalgebra
with product $\mu$.
Then:
\begin{itemize}
\item[$(a)$] The $\Z$-grading of $W(J)$
induces on $Lie(J,\mu)$ a $\Z$-grading 
$Lie(J,\mu)=\prod_{k\geq -1}Lie_k(J)$,
where
$Lie_{-1}(J)=J$, $Lie_0(J)=Str(J,\mu)$ and $Lie_1(J)=R(J,\mu)$, which
is admissible if $J$ is rigid.
\item[$(b)$]
$J$ is simple if and only if $Lie(J,\mu)$ is an irreducible
$\Z$-graded Lie superalgebra;
\item[$(c)$] If $J$ is linearly compact then $Lie(J,\mu)$ is linearly
compact.
\end{itemize}
\end{proposition}
{\bf Proof.}
The transitivity of the grading and properties (\ref{b}) and (\ref{c})
follow from the definitions.
Moreover, if $J$ is a rigid superalgebra, then also 
property (\ref{a}) holds, i.e., 
the $\Z$-grading defined on $Lie(J,\mu)$ is admissible, proving $(a)$.
 Denote by $\mu$ the commutative product of $J$. By construction, if $x,y$ lie 
in $J$,
then $x\circ y=[[\mu,x],y]$, hence a proper ideal of $J$ is
a proper $Str(J,\mu)$-submodule of $Lie_{-1}(J)$, proving $(b)$.
Finally, $(c)$ follows from $(a)$ and \cite[Lemma 2.1]{CantaK2}.
\hfill$\Box$

\begin{theorem}\label{LieJ}
(A) If $J$ is a simple linearly compact rigid superalgebra with 
an even product $\mu$, 
then one of the following two possibilities holds:
\begin{enumerate}
\item[$(a)$] $Lie(J,\mu)$ is a simple linearly compact Lie superalgebra;
\item[$(b)$] $Lie(J,\mu)=S\rtimes\F \mu$, where $S$ is a simple
linearly compact Lie superalgebra
and $\mu$ is an even outer derivation of $S$.
\end{enumerate}
($A^\prime$) If $J$ is a simple linearly compact rigid odd type 
superalgebra with
an odd product $\mu$,  
then one of the following four possibilities holds:
\begin{enumerate}
\item[$(a')$] $Lie(J,\mu)$ is a simple linearly compact Lie superalgebra;
\item[$(b')$] $Lie(J,\mu)=S+\F\mu+\F [\mu,\mu]$, where $S$ is a simple
linearly compact Lie superalgebra and $\mu$ is an odd outer derivation of $S$
such that $[\mu,\mu]\neq 0$;
\item[$(c')$] $Lie(J,\mu)\cong\xi\a+\a+d/d\xi$, where $\a$ is a simple
linearly compact Lie superalgebra and $\xi$ is an 
odd indeterminate;
\item[$(d')$] $Lie(J,\mu)=S\otimes\Lambda(1)+\F \mu+\F d$ where $S$ is a
  simple linearly compact Lie superalgebra, 
$d$ is an even outer derivation of $S$ and $\mu=d\otimes\xi+1\otimes d/d\xi$.
\end{enumerate}
\end{theorem}
{\bf Proof.} By Proposition
\ref{irr}, $Lie(J,\mu)$ is a transitive irreducible $\Z$-graded Lie 
superalgebra. 
Let $I'$ be a non-zero closed ideal of 
$Lie(J,\mu)$. 
Then, by transitivity, 
$I'\cap Lie_{-1}(J,\mu)\neq\emptyset$,
hence, by the irreducibility of the grading, $I'\cap Lie_{-1}(J,\mu)=J$. 
Now, let $I$ be the intersection of all non-zero closed ideals of
$Lie(J,\mu)$. Then, by the above remark, $I$ is a minimal closed ideal of
$Lie(J,\mu)$, containing $J$. therefore, by construction,  
$Lie(J,\mu)=I+\F \mu+\F [\mu,\mu]$. Next, by the super-analogue of 
Cartan-Guillemin's theorem
\cite{B,G1}, established in \cite{FK}, 
$I=S\hat{\otimes}\Lambda(m,n)$, for some simple linearly compact
Lie superalgebra $S$ and some
$m,n\in\Z_{\geq 0}$, and $\mu$ lies in $Der(S\hat{\otimes}\O(m,n))$. Since
$Der(S\hat{\otimes}\O(m,n))=Der S\hat{\otimes}\O(m,n)+1\otimes W(m,n)$
\cite{FK}, we have: 
$\mu=\sum_i(d_i\otimes a_i) +1\otimes \mu'$ for some $d_i\in Der S$,
$a_i\in\O(m,n)$ and $\mu'\in W(m,n)$.

First consider the case when $J$ is a simple rigid superalgebra.
Then $\mu'$ is an even element of $W(m,n)$ hence, by the minimality of
the ideal $S\hat{\otimes} \O(m,n)$, $n=0$. Now suppose 
$m\geq 1$. If $\mu'$ lies
in the non-negative part of $W(m,0)$ with the principal grading, then the
ideal generated by $Sx_1$ is a proper $\mu$-invariant ideal of 
$S\hat{\otimes} \O(m,0)$,
contradicting its minimality. Therefore we may assume, up to a linear change 
of indeterminates, that $\mu'=\frac{\partial}{\partial x_1}+D$, for some
derivation $D$ lying in the non-negative part of
$W(m,0)$. Since $\mu$ lies in $Lie_1(J)$, we have $\deg(x_1)=-1$, but this
is a contradiction since the $\Z$-grading of $Lie(J,\mu)$ has depth 1.
It follows that $m=0$. Therefore, either $\mu$ is an inner derivation of $S$
and $Lie(J,\mu)=S$, or $\mu$ is an outer derivation of $S$ and 
$Lie(J,\mu)=S\rtimes\F \mu$ ($[\mu,\mu]=0$ since $\mu$ is even). 

Now consider the case when $J$ is a simple rigid odd type
superalgebra.
Consider the principal grading of $W(m,n)$ and
denote by $W(m,n)_{\geq 0}$ its non-negative part. If $\mu'\in
W(m,n)_{\geq 0}$, then the minimality of the ideal
$S\hat{\otimes}{\cal O}(m,n)$ implies $m=n=0$, hence we get $(a')$, in
case 
$\mu$ is an inner derivation of $S$, or we get $(b')$, in case $\mu$ is
an outer derivation of $S$ such that $[\mu,\mu]\neq 0$, or 
we get $(c')$ by Proposition \ref{[p,p]=0}, in case $\mu$ is
an outer derivation of $S$ such that $[\mu,\mu]=0$. Now let $n\geq 1$
and suppose that $\mu'$ has a non-zero projection on $W(m,n)_{-1}$. Then,  up
to a linear change of indeterminates, 
$\mu'=\frac{\partial}{\partial\xi_1}+D$ for some odd derivation $D\in
W(m,n)_{\geq 0}$.   Write $D=\xi_1 D_0+D_1$ where $D_0\in W(m,n)_{\0}
\cap (\O(x_1,\dots,x_m,
\xi_2,\dots,\xi_n)\otimes W(m,n)_{-1})$
and $D_1\in W(m,n)_{\1}\cap (\O(x_1,\dots,x_m,
\xi_2,\dots,\xi_n)\otimes W(m,n)_{-1})$. If $n\geq 2$, 
then  the ideal of $S\hat{\otimes}\O(m,n)$ generated
by $S\xi_2$ is a proper $\mu$-stable ideal of
$S\hat{\otimes}\O(m,n)$, contradicting the  minimality of 
$S\hat{\otimes}\O(m,n)$. It follows that
$n=1$. We shall now prove that $m=0$. Indeed, if $D_0$ lies in 
$W(m,1)_{\geq 0}$, then $m=0$,
otherwise the ideal generated by $Sx_1$ would be a proper $\mu$-stable
ideal of $S\hat{\otimes}\O(m,1)$.
Now suppose that $m\geq 1$ and that $D_0$ has non-zero projection on 
$W(m,1)_{-1}$,
i.e.,  up to a linear change of indeterminates,
$D_0=\frac{\partial}{\partial x_1}+\delta_0$ for some $\delta_0\in W(m,0)_{\geq
  0}$ and $\mu'=\frac{\partial}{\partial
  \xi_1}+\xi_1\frac{\partial}{\partial x_1}+\xi_1\delta_0+D_1$. 
Since $\mu$ lies in $Lie_1(J)$, we have: $\deg(\xi_1)=-1$ and
$\deg(x_1)=-2$, but this is a contradiction since the $\Z$-grading of
$Lie(J,\mu)$ has depth 1. We thus have proved that $m=0$, i.e., 
$Lie(J,\mu)=S\otimes\Lambda(1)+\F \mu+\F[\mu,\mu]$.
Notice that either $\mu=\delta\otimes 1+\alpha\otimes\frac{\partial}{\partial
  \xi_1}$  for some $\delta\in (Der S)_{\1}$ and some $\alpha\in\F$
(resp.\ $\mu=d\otimes\xi_1$ for some $d\in (Der S)_{\0}$),
 hence $[\mu,\mu]=0$ and, by
Proposition \ref{[p,p]=0}, $(c')$ holds, or
$\mu=d\otimes \xi_1+1\otimes\frac{\partial}{\partial
  \xi_1}$ for some even derivation $d\in Der S$, hence $[\mu,\mu]=2d$,
 and $Lie(J,\mu)=S\otimes\Lambda(1)+\F\mu+\F d$. If $d$ is an inner derivation
of $S$ we are in case $(c')$, otherwise $(d')$ holds.
\hfill$\Box$

\section{Classification of admissible gradings}\label{eight}
The following proposition lists the $\Z$-gradings of depth 1 of all simple
linearly compact Lie superalgebras. Following \cite{K2} we denote
a finite-dimensional contragredient Lie superalgebra by 
$G(A,\tau)/C$
(where $C$ is a central ideal of $G(A,\tau)$),
we denote
its standard Chevalley generators by $e_i$ and $f_i$, and let  
$\sum_ia_i\alpha_i$ be the highest root. 
Likewise, we will denote by $e_i$, $f_i$ the standard
Chevalley generators of $q(n)_{\bar{0}}$ and by $\bar{e}_i$, $\bar{f}_i$
the corresponding elements in $q(n)_{\bar{1}}$.
Besides, we will identify $p(n)$ with the subalgebra of the Lie superalgebra 
$SHO(n,n)$ spanned by the following elements:
$\{x_ix_j, \xi_i\xi_j: i,j=1,\dots,n;~ x_i\xi_j: i\neq j=1,\dots,n;~
x_i\xi_{i}-x_{i+1}\xi_{i+1}: i=1, \dots,n-1\}$ (cf.\ \cite[\S 1.3]{CK}),
and thus describe the $\Z$-gradings of $p(n)$ as induced by $\Z$-gradings of 
$SHO(n,n)$.

\begin{proposition}\label{listdepth1} All $\Z$-gradings of depth 1 of all 
simple linearly compact Lie
superalgebras $S$ are, up to isomorphism, the following:
\begin{enumerate}
\item $S=G(A,\tau)/C$:
$\deg e_i=-\deg f_i=k_i$
with $k_s=1$ for $s$ such that $a_s=1$ and $k_i=0$ for every $i\neq s$;
\item $S=q(n)$: $\deg(e_i)=\deg(\bar{e}_i)=-\deg(f_i)=-\deg(\bar{f}_i)=k_i$
with $k_s=1$ for some $s$ and $k_i=0$ for every $i\neq s$;
\item $S=p(n)$, $n\geq 2$: 
\begin{enumerate}
\item $(1,0,\dots,0|-1,0,\dots,0)$;
\item $(\underbrace{1,\dots,1}_h,0,\dots,0|\underbrace{0,\dots,0}_h,1,\dots,1)$ with $h=0,\dots,n$;
\end{enumerate}
\item $S=W(m,n)$, $(m,n)\neq (0,1)$; $S(m,n)$, $m>1$, or $m=0$ and $n\geq 3$,
or $m=1$ and $n\geq 2$:
\begin{enumerate}
\item $(\underbrace{1,\dots,1}_h,0,\dots,0|
\underbrace{1,\dots,1}_k,0,\dots,0)$ 
with $h=0,\dots,m$, $k=0,\dots,n$;
\item $(0,\dots,0|-1,0,\dots,0)$;
\end{enumerate}
\item $S=H(2k,n)$, $k\geq 1$, or $k=0$ and $n\geq 4$:
\begin{enumerate}
\item $(1,\dots,1|1,\dots,1)$;
\item $n\geq 2$, $(0,\dots,0|1,0,\dots,0,-1)$;
\item $n=2t$, $(\underbrace{1,\dots,1}_h,
\underbrace{0,\dots,0}_{k-h},
\underbrace{0,\dots,0}_h,1,\dots,1|\underbrace{1,\dots,1}_t,0,\dots,0)$
with $h=0,\dots,k$;
\end{enumerate}
\item $S=K(2k+1,n)$:
\begin{enumerate}
\item $n\geq 2$, $(0,\dots,0|1,0,\dots,0,-1)$;
\item $n=2t$, $(1,\underbrace{1,\dots,1}_h,\underbrace{0,\dots,0}_{k-h},
\underbrace{0,\dots,0}_h,1,\dots,1|\underbrace{1,\dots,1}_t,0,\dots,0)$
with $h=0,\dots,k$;
\end{enumerate}
\item $S=HO(n,n)$, $n\geq 2$, $SHO(n,n)$, $n\geq 3$:
\begin{enumerate}
\item $(1,\dots,1|1,\dots,1)$;
\item $(0,\dots,0,1|0,\dots,0,-1)$;
\item $(\underbrace{1,\dots,1}_h,0,\dots,0|\underbrace{0,\dots,0}_h,1,\dots,1)$
with $h=0,\dots,n$;
\end{enumerate}
\item 
$S=KO(n,n+1)$, $n\geq 2$, $SKO(n,n+1;\beta)$, $n\geq 2$:
\begin{enumerate}
\item $(0,\dots,0,1|0,\dots,0,-1,0)$;
\item $(\underbrace{1,\dots,1}_h,0,\dots,0|
\underbrace{0,\dots,0}_h,1,\dots,1,1)$ with $h=0,\dots,n$;
\end{enumerate}
\item $S=SKO(2,3;0)$: 
\begin{enumerate}
\item $(1,1|-1,-1,0)$;
\end{enumerate}
\item $S=SHO^\sim(n,n)$, $n\geq 2$, even: 
\begin{enumerate}
\item $(1,\dots,1|0,\dots,0)$;
\end{enumerate}
\item $S=E(5,10)$:
\begin{enumerate}
\item $(1,1,1,1,0)$;
\item $(1,1,0,0,0)$;
\end{enumerate}
\item $S=E(4,4)$:
\begin{enumerate}
\item $(1,1,1,1)$;
\item $(1,1,0,0)$;
\end{enumerate}
\item $S=E(3,6)$:
\begin{enumerate}
\item $(1,0,0;\frac{1}{2})$;
\item $(1,1,0;0)$;
\item $(1,1,1;\frac{1}{2})$;
\end{enumerate}
\item $S=E(3,8)$:
\begin{enumerate}
\item $(1,0,0;0)$;
\item $(1,0,0;-1)$;
\item $(1,1,0;-1)$;
\end{enumerate}
\item $S=E(1,6)$:
\begin{enumerate}
\item $(0|1,0,0,-1,0,0)$;
\item $(1|1,1,1,0,0,0)$;
\item $(1|1,1,0,0,0,1)$.
\end{enumerate}
\end{enumerate}
\end{proposition}
{\bf Proof.} If $S$ is a finite-dimensional contragredient Lie superalgebra,
any $\Z$-grading of $S$ is defined by setting 
$\deg(e_i)=-\deg(f_i)=k_i\in\Z_{\geq 0}$.
Then it is clear that such a grading has depth one if and only if all
$k_i$'s are 0 except for $k_s=1$ for some $s$ such that $a_s=1$.
Likewise, one gets the statement for $S=q(n)$, using that $q(n)_{\bar{0}}\cong
A_n$.

Recall that any $\Z$-grading of finite depth of $W(m,n)$ is isomorphic to
that of type $(a_1, \dots, a_m|$ $b_1, \dots, b_n)$ with $a_i\in\Z_{\geq 0}$
\cite{CK}, \cite{CantaK}.
Suppose that it has depth 1. Then 
$\deg \frac{\partial}{\partial x_i}=-a_i\geq -1$ and 
$\deg \frac{\partial}{\partial\xi_j}=-b_j\geq -1$, i.e., the grading is 
isomorphic to a grading of type $(1, \dots, 1,0,\dots,0|b_1, \dots, b_n)$
with $h$ 1's and $m-h$ 0's, for some $h=0, \dots, m$ and $b_j\leq 1$.
If $n=0$, then the statement for $S=W(m,n)$ is proved.
Let $n>0$ and assume $m\geq 1$. If $h\geq 1$, then $b_j\geq 0$, since
$\deg\xi_j\frac{\partial}{\partial x_1}=b_j-1$, i.e., the grading is isomorphic
to the grading of type $(1, \dots, 1,0,\dots,0|1,\dots,1,0,\dots,0)$ with
$h+k$ 1's, for some $h=1, \dots, m$ and some $k=0,\dots,n$.
If $m\geq 1$ and $h=0$, then $b_j\geq -1$, since 
$\deg\xi_j\frac{\partial}{\partial x_1}=b_j$. In this case, either $b_j\geq 0$
for every $j=1, \dots,n$, and we are in the same situation as above, or
there exists some $s$ such that $b_s=-1$. Since 
$\deg\xi_s\frac{\partial}{\partial \xi_t}=b_s-b_t$, it follows that
$b_t=0$ for every $t\neq s$, i.e.\ the grading is isomorphic to
the grading of type $(0,\dots,0|-1,0,\dots,0)$.
Finally, if $m=0$ and $n\geq 2$, then $b_i\geq -1$
and 
$b_i-b_j\geq -1$, since
$\deg\xi_i\xi_j\frac{\partial}{\partial\xi_j}=b_i$ and
$\deg\xi_i\frac{\partial}{\partial \xi_j}=b_i-b_j$. It follows
that the grading is isomorphic either to the grading of type
$(1,\dots,1,0,\dots,0)$ with $h$ 1's, for some $h=1, \dots,n$, or
to the grading of type $(-1,0,\dots,0)$. Note that, if $n=2$, then the grading
of type $(-1,-1)$ is isomorphic to the grading of type $(1,1)$,
and if $n\geq 3$, then $b_i+b_j\geq -1$, since $\deg\xi_i\xi_j\xi_k
\frac{\partial}{\partial\xi_k}=b_i+b_j$.

The arguments for the closed subalgebras of $W(m,n)$ are similar.

If $S$ is a simple infinite-dimensional exceptional Lie superalgebra, then one
uses the description of the $\Z$-gradings of $S$ given in \cite{CK} and
\cite{CantaK}. For example, let $S=E(5,10)$. Then $E(5,10)_{\bar{0}}\cong S_5$
and $E(5,10)_{\bar{1}}\cong d\Omega^1(5)$ as an $S_5$-module.
A $\Z$-grading of finite depth of $E(5,10)$ is defined by setting $\deg x_i=
-\deg\frac{\partial}{\partial x_i}=a_i\in\Z_{\geq 0}$ and 
$\deg d=-1/4\sum_{i=1}^5a_i
\in \frac{1}{2}\Z$. Therefore, if such a grading has depth one, then 
$a_i\leq 1$, hence the statement holds, since $\sum_{i=1}^5a_i\in 2\Z$.

Similar arguments apply to $S=E(4,4), E(3,6), E(3,8)$ and $E(1,6)$.\hfill$\Box$

\begin{remark}\label{isoSKO(2,3;beta)}\em
Let us consider the Lie superalgebra $SKO(2,3;\beta)$. For every $\beta\neq 0, -1$
we have:
$SKO(2,3;\beta)_{\0}\cong W(2,0)$,
$SKO(2,3;\beta)_{\1}\cong \Omega^0(2)^{-\frac{1}{\beta+1}}
\oplus\Omega^0(2)^{-\frac{\beta}{\beta+1}}$  \cite[Remark 4.15]{CantaK}.
Let us consider $W(2,0)$ and $\Omega^0(2)$ with the grading of type $(0,1|)$, 
so that
$W(2,0)=\prod_{j\geq -1}W(2,0)_j$ and
$\Omega^0(2)=\prod_{j\geq 0}\Omega^0(2)_j$.
Denote by $\g=\prod_{j\geq -1}\g_j$ and 
$\g=\prod_{j\geq -1}\h_j$ the Lie superalgebra $SKO(2,3;\beta)$ with
the gradings of type 
$(0,1|1,0,1)$ and $(0,1|0,-1,0)$, respectively. Then, for every $\beta\neq 0,-1$ we have:
$(\g_j)_{\0}=(\h_j)_{\0}\cong W(2,0)_j$,
$(\g_j)_{\1}\cong 
\Omega^0(2)_j^{-\beta/(\beta+1)}\oplus\Omega^0(2)_{j+1}^{-1/(\beta+1)}$ and
$(\h_j)_{\1}\cong 
\Omega^0(2)_j^{-1/(\beta+1)}\oplus\Omega^0(2)_{j+1}^{-\beta/(\beta+1)}$
($\Omega^0(2)_{-1}=0$).
It follows that for every
$\beta\neq 0, -1$, $SKO(2,3;\beta)$ with the grading of type $(0,1|0,-1,0)$ 
is isomorphic
to $SKO(2,3;1/\beta)$ with the grading of type $(0,1|1,0,1)$.

\medskip

If $\beta=-1$, then $SKO(2,3;-1)_{\0}$ is not simple: it has a commutative 
ideal isomorphic to $\Omega^0(2)$. We have: 
$SKO(2,3;-1)_{\0}\cong \Omega^0(2)\rtimes S(2,0)$ and
$SKO(2,3;-1)_{\1}\cong \Omega^0(2)_+\oplus \Omega^0(2)_-$,
where $\Omega^0(2)_{+}$ and $\Omega^0(2)_{-}$ are two odd copies of
$\Omega^0(2)$ on which $S(2,0)$ acts 
in the natural way. The even functions in $\Omega^0(2)$
act by multiplication on $\Omega^0(2)_{+}$ and by $-$multiplication on
$\Omega^0(2)_{-}$
\cite[Remark 4.19]{CantaK}. Besides, if  
$f_+\in \Omega^0(2)_{+}$ and $g_-\in\Omega^0(2)_{-}$, then
$[f_+, g_-]=\frac{\partial f_+}{\partial x_1}\frac{\partial g_-}{\partial x_2}
-\frac{\partial f_+}{\partial x_2}\frac{\partial g_-}{\partial x_1}-
(f_+dg_-+g_-df_+)$. Notice that if  $f_-$ and $g_+$ are the corresponding elements
in $\Omega^0(2)_{-}$ and $\Omega^0(2)_{+}$, respectively, then
$[f_-, g_+]=-\frac{\partial f_-}{\partial x_1}\frac{\partial g_+}{\partial x_2}
+\frac{\partial f_-}{\partial x_2}\frac{\partial g_+}{\partial x_1}-
(f_-dg_++g_+df_-)$.

As above, let us denote by $\g=\prod_{j\geq -1}\g_j$
and $\g=\prod_{j\geq -1}\h_j$ the Lie superalgebra
$SKO(2,3;-1)$ with respect to the gradings of type $(0,1|1,0,1)$ and
$(0,1|0,-1,0)$, respectively.
We have:
$(\g_j)_{\0}=(\h_j)_{\0}\cong S(2,0)_j+\Omega^0(2)_j$,
$(\g_j)_{\1}\cong 
(\Omega^0(2)_-)_j\oplus(\Omega^0(2)_+)_{j+1}$
and $(\h_j)_{\1}\cong 
(\Omega^0(2)_+)_j\oplus(\Omega^0(2)_-)_{j+1}$.
It follows that the map $\Phi: \g=\prod\g_j  \longrightarrow \g=\prod\h_j$,
$\Phi(a)=-a$ for $a\in \Omega^0(2)$, $\Phi(X)=X$ for $X\in S(2,0)$,
$\Phi(f_+)=f_-$ for $f_+\in(\Omega^0(2)_+)$,
$\Phi(f_-)=f_+$ for $f_-\in(\Omega^0(2)_-)$, is an isomorphism
of $\Z$-graded Lie superalgebras.
\end{remark}

\begin{theorem}\label{eag}\begin{itemize}
\item[i)] An even admissible grading of a
simple linearly compact Lie superalgebra is either short or isomorphic to
one of the following:
\begin{enumerate}
\item $S=W(0,3)$, $(|1,1,0)$;
\item $S=W(0,4)$, $(|1,1,0,0)$;
\item $S=S(0,4)$, $(|1,1,0,0)$;
\item $S=S(0,5)$, $(|1,1,0,0,0)$;
\item $S=S(1,4)$, $(0|1,1,0,0)$;
\item $S=S(2,0)$, $(1,1|)$.
\item $S=S(2,0)$, $(1,0|)$.
\end{enumerate}
\item[ii)] An odd admissible grading of a
simple linearly compact Lie superalgebra is either short or isomorphic to
one of the following:
\begin{enumerate}
\item $S=W(1,1)$, $(1|0)$;
\item $S=W(1,2)$, $(0|1,1)$;
\item $S=W(2,1)$, $(1,0|0)$;
\item $S=W(2,2)$, $(0,0|1,1)$;
\item $S=S(1,3)$, $(0|1,1,1)$;
\item $S=S(2,2)$, $(0,0|1,1)$;
\item $S=H(2,2)$, $(1,0|1,0)$;
\item $S=HO(2,2)$, $(1,0|0,1)$;
\item $S=HO(3,3)$, $(0,0,0|1,1,1)$;
\item $S=HO(n+1,n+1)$ with $n\geq 1$, $(0,\dots,0,1|0,\dots,0,-1)$;
\item $S=SHO(3,3)$, $(1,0,0|0,1,1)$;
\item $S=SHO(4,4)$, $(0,0,0,0|1,1,1,1)$;
\item $S=SHO(n+1,n+1)$ with $n\geq 2$, $(0,\dots,0,1|0,\dots,0,-1)$;
\item $S=KO(n+1,n+2)$ with $n\geq 1$, $(0,\dots,0,1|0,\dots,0,-1,0)$;
\item $S=KO(2,3)$, $(0,0|1,1,1)$;
\item $S=SKO(2,3;\beta)$, $(1,0|0,1,1)$;
\item $S=SKO(2,3;\beta)$, $(0,1|0,-1,0)$;
\item $S=SKO(3,4;\beta)$, $(0,0,0|1,1,1,1)$;
\item $S=SKO(n,n+1;\beta)$ with $n>2$ and $\beta\neq\frac{4}{n}$, 
$(0,\dots,0,1|0,\dots,0,-1,0)$.
\end{enumerate}
\end{itemize}
\end{theorem}
{\bf Proof.} By Remark \ref{necessarycond}, among the $\Z$-gradings of 
depth 1 of all
simple linearly compact Lie superalgebras $S=\prod_{j\geq -1}S_j$ listed in
Proposition \ref{listdepth1}, we select the non-short ones such that either 
the dimension of $S_0$
is infinite and  
$\growth S_1 \leq \growth S_0$ and $\size S_1 \leq \size S_0$
, or the dimension of $S_0$ is finite and
$\dim S_1\leq\dim S_0+1$. 
Up to isomorphism, we get 
a grading in the following list:
\begin{enumerate}
\item[A)] $S=p(n)$, $(1,0,\dots,0|-1,0,\dots,0)$;
\item[B)] $S=W(0,3)$, $(|1,1,1)$;
\item[C)] $S=W(0,3)$ or $S=S(0,3)$, $(|1,1,0)$;
\item[D)] $S=W(0,4)$ or $S=S(0,4)$, $(|1,1,0,0)$;
\item[E)] $S=W(1,2)$, $(0|1,1)$;
\item[F)] $S=W(1,3)$ or $S=S(1,3)$, $(0|1,1,0)$;
\item[G)] $S=W(2,2)$ or $S=S(2,2)$, $(0,0|1,1)$;
\item[H)] $S=W(m,n)$ or $S=S(m,n)$ with $m\geq 1$, $(1,0,\dots,0|0,\dots,0)$;
\item[I)] $S=S(0,4)$, $(|1,1,1,0)$;
\item[J)] $S=S(0,5)$, $(|1,1,0,0,0)$;
\item[K)] $S=S(1,2)$, $(1|1,1)$;
\item[L)] $S=S(1,3)$, $(0|1,1,1)$;
\item[M)] $S=S(1,4)$, $(0|1,1,0,0)$;
\item[N)] $S=S(2,3)$, $(0,0|1,1,0)$;
\item[O)] $S=S(3,2)$, $(0,0,0|1,1)$;
\item[P)] $S=H(0,5)$, $(|1,1,1,1,1)$;
\item[Q)]$S=S(2,0)$, $(1,1|)$;
\item[R)] $S=H(2,2)$, $(1,0|1,0)$;
\item[S)] $S=H(0,6)$, $(|1,1,1,0,0,0)$;
\item[T)] $S=HO(n,n)$ with $n\geq 2$ and $S=SHO(n,n)$ with $n\geq 3$, $(0,\dots,0,1|0,\dots,0,-1)$;
\item[U)] $S=HO(2,2)$, $(1,0|0,1)$;
\item[V)] $S=HO(3,3)$, $(0,0,0|1,1,1)$;
\item[W)] $S=SHO(3,3)$, $(1,0,0|0,1,1)$;
\item[X)] $S=SHO(4,4)$, $(0,0,0,0|1,1,1,1)$;
\item[Y)] $S=KO(n,n+1)$ and $S=SKO(n,n+1,\beta)$ with $n\geq 2$, 
$(0,\dots,0,1|0,\dots,0,-1,0)$;
\item[Z)] $S=KO(2,3)$, $(0,0|1,1,1)$;
\item[Z$^\prime$)] $S=SKO(3,4;\beta)$, $(0,0,0|1,1,1,1)$.
\end{enumerate}

Note that the grading of type $(1|1)$ and the
grading of type $(1|0)$ of $W(1,1)$ are isomorphic. 
Moreover,  in the above list we omitted  
the Lie superalgebra $H(2,0)$, since it is isomorphic to $S(2,0)$,
the Lie superalgebras $KO(1,2)$ and $K(1,2)$, 
since they are isomorphic to $W(1,1)$, and the Lie
superalgebra $S(2,1)$  since it is isomorphic to 
$HO(2,2)$. Besides, $SKO(2,3;0)\cong HO(2,2)$
and, due to Remark \ref{isoSKO(2,3;beta)}, for every $\beta\neq 0$
 the Lie superalgebra $SKO(2,3;\beta)$ with the grading of type $(0,1|1,0,1)$
is isomorphic to the Lie superalgebra $SKO(2,3;\frac{1}{\beta})$ with the
grading of type $(0,1|0,-1,0)$.

All gradings A)--Z$^\prime$) are transitive and satisfy property
(\ref{c}), with the exception of the grading of type
$(1|)$ of $W(1,0)$, since in this case $S_1$ is one-dimensional 
hence $S_2\neq S_1^2$. In order to prove that one of these gradings 
is even admissible (resp.\ odd admissible)
 it is therefore sufficient to show that there exists an
element $\mu\in (S_1)_{\bar{0}}$ (resp.\ $\mu\in (S_1)_{\bar{1}}$)
satisfying properties (\ref{a}) and (\ref{b}).
In Tables 1 and 2 we list, up to isomorphism,
all even admissible and odd admissible gradings and indicate
such an element $\mu$ for each of them. One shows that properties
(\ref{a}) and (\ref{b}) hold by direct computation,
as in Example \ref{HO}.
For notation concerning finite-dimensional and infinite-dimensional Lie superalgebras,  
not introduced in this paper, we refer to \cite{K2} and
\cite{CantaK}, respectively.

In all cases not appearing in Tables 1 and 2, either property (\ref{a}) or 
property
(\ref{b}) fails for every even or odd element in $S_1$. For
instance, let $S=W(m,n)$ with $m\geq 1$ and $n\geq 1$,
and consider the grading of $S$ of type $(1,0,\dots,0|0,\dots,0)$.
Then $S_{-1}=\langle \frac{\partial}{\partial x_1}\rangle\otimes A$,
where $A=\F[[x_2,\dots,x_m]]\otimes\Lambda(n)$,
$S_0=\langle x_1\frac{\partial}{\partial x_1},
\frac{\partial}{\partial x_i}, \frac{\partial}{\partial \xi_j} | i=2,\dots, m,
j=1,\dots,n\rangle\otimes A$ and 
$S_1=\langle x_1^2\frac{\partial}{\partial x_1},
x_1\frac{\partial}{\partial x_i}, x_1\frac{\partial}{\partial \xi_j} | i=2,\dots, m,
j=1,\dots,n\rangle\otimes A$. It follows that for every
$\mu\in (S_1)_{\bar{0}}$, the subalgebra of $S_0$ generated by
 $[S_{-1},\mu]$ does not
contain the vector fields $\frac{\partial}{\partial \xi_j}$,
hence $[S_{-1},\mu]$ does not  generate $S_0$.
The same kind of argument applies to the gradings A), C) for $S=S(0,3)$, 
F) for $S=S(1,3)$, H) for $S=W(m,0)$ with $m\geq 2$, and $S=S(m,n)$ unless
$m=2$, $n=0$, K), I), N), R), S), T), U), W),  
hence these gradings
are not even admissible. Note that the 1st graded component of $S$ with respect to
gradings B), E), G), L), O), P), V), X), Z), Z$^\prime$)
 is completely odd. 

Let $S=W(1,3)$ with the grading of type $(0|1,1,0)$. Then
$S_{-1}=\langle\frac{\partial}{\partial\xi_1},
\frac{\partial}{\partial\xi_2}\rangle\otimes A$ with $A=\F[[x]]\otimes
\Lambda(\xi_3)$,
$S_{0}=\langle\frac{\partial}{\partial x},
\frac{\partial}{\partial\xi_3},
\xi_i\frac{\partial}{\partial\xi_j} | i,j=1,2
\rangle\otimes A$,
$S_{1}=\langle\xi_i\frac{\partial}{\partial x},
\xi_i\frac{\partial}{\partial\xi_3},
\xi_1\xi_2\frac{\partial}{\partial\xi_i} | i=1,2
\rangle\otimes A$. $S_0$ and $S_1$ are therefore
infinite-dimensional subspaces of growth 1 and size 12.
One shows that if $\mu$ lies in $(S_1)_{\0}$, then
the centralizer of $\mu$ in $S_0$ has size strictly greater than 0,
therefore condition (\ref{a}) of Definition \ref{admissible}
cannot hold and the grading is not even admissible.
We thus get the list of cases in $i)$.

Let $S=W(0,3)$ with the grading of type $(|1,1,0)$. Then
$S_{-1}=\langle\frac{\partial}{\partial\xi_1},
\frac{\partial}{\partial\xi_2}\rangle\otimes A$ with $A=\Lambda(\xi_3)$,
$S_0=\langle\xi_i\frac{\partial}{\partial\xi_j},
\frac{\partial}{\partial\xi_3} | i,j=1,2\rangle\otimes A$ and
$S_1=\langle\xi_1\xi_2\frac{\partial}{\partial\xi_i},
\xi_i\frac{\partial}{\partial\xi_3} | i=1,2\rangle\otimes A$.
It follows that, for every $\mu\in (S_1)_{\bar{1}}$,
the subalgebra of $S_0$ generated by $[S_{-1},\mu]$
does not contain the vector field $\frac{\partial}{\partial\xi_3}$
lying in $S_0$, therefore the grading is not odd admissible.
The same kind of argument applies to the gradings A),
C), D), F), H) unless $(m,n)=(1,1)$ or $(m,n)=(2,1)$, I), J), M), N), O), S), 
hence
these gradings are not odd admissible.  Notice that the superalgebra
$S=S(2,0)$ 
in Q) is completely even.

Let $S=H(0,5)$ with the principal grading, i.e., the grading
of type $(|1,1,1,1,1)$. Then $S_{-1}=\langle\xi_i
| i=1,\dots, 5\rangle$,
$S_0=\langle\xi_i\xi_j | i,j=1,\dots, 5\rangle$ and 
$S_1=\langle\xi_i\xi_j\xi_k | i,j,k=1,\dots, 5\rangle$. It follows
that $\dim(S_0)=\dim(S_1)=10$. One shows that, for every choice of
$\mu\in(S_1)_{\1}$, $\mu$ lies in $[S_0,\mu]$ and the
centralizer of $\mu$ in $S_0$ has positive dimension. It follows that
condition (\ref{a}) of Definition \ref{admissible} cannot be satisfied,
hence the grading is not odd admissible. Similarly, 
consider $S=W(0,3)$ with the grading of type $(|1,1,1)$. Then
$S_0\cong gl_3$ and $S_1\cong \F^3\oplus S^2(\F^3)^*$, where
$\F^3$ denotes the standard $gl_3$-module. The stabilizer of a generic
point in $S^2(\F^3)^*$ is isomorphic to $so_3$ and the stabilizer
of $so_3$ in $\F^3$ is isomorphic to $\F$. It follows that the
codimension of $[S_0,\mu]$ in $S_1$, for a generic point $\mu\in S_1$, is
equal to one. Since $\mu\in[S_0,\mu]$, condition (\ref{a}) cannot
be satisfied, hence the grading is not odd admissible.
The same kind of argument applies to K) and Y) for $S=SKO(n,n+1;4/n)$ 
when $n\geq 2$,
therefore these gradings are not odd admissible.
We thus get the list
in $ii)$.
\hfill$\Box$

\medskip

In Theorem \ref{orbits} we will prove that if $\g=\prod_{j\geq -1}\g_j$ is a 
simple linearly compact Lie superalgebra with an admissible grading, then 
Tables 1 and 2 below provide a complete list, up to an automorphism from 
$exp(ad(\g_{int}))$, where $\g_{int}$ is the subalgebra of $(\g_0)_{\bar{0}}$, 
generated by its elements which are $ad$-exponentiable in $\g$, of even and 
odd elements $\mu\in \g_1$, satisfying properties (\ref{a}) and 
(\ref{b}).

\bigskip

\begin{center}
\begin{tabular}{|c|c|c|}
\hline
$S$ & type of grading & $\mu$\\
\hline
\vspace{0.001cm} & \vspace{0.001cm} & \vspace{0.001cm}\\
$W(0,3)$ & $(|1,1,0)$ & $\xi_1\frac{\partial}{\partial\xi_3}+
\xi_1\xi_2\xi_3(\frac{\partial}{\partial\xi_1}+
\frac{\partial}{\partial\xi_2})$ \\
$W(0,4)$ & $(|1,1,0,0)$ & $P_W=\xi_1(\frac{\partial}{\partial\xi_3}+
\xi_2\xi_3\frac{\partial}{\partial\xi_1}+
\xi_2\xi_3\frac{\partial}{\partial\xi_2})+
\xi_2(\frac{\partial}{\partial\xi_4}+
\xi_3\xi_4\frac{\partial}{\partial\xi_3})$ \\
$S(0,4)$ & $(|1,1,0,0)$ & $P_W+\xi_1\xi_3\xi_4
\frac{\partial}{\partial\xi_4}+\xi_1\xi_2\xi_4
\frac{\partial}{\partial\xi_1}-
\xi_2\xi_3\xi_4\frac{\partial}{\partial\xi_4}$ \\
$S(0,5)$ & $(|1,1,0,0,0)$ & $\xi_1(\frac{\partial}{\partial\xi_3}+
\xi_3\xi_4\frac{\partial}{\partial\xi_5})+
\xi_2(\frac{\partial}{\partial\xi_4}+
\xi_3\xi_4\frac{\partial}{\partial\xi_3}+
\xi_4\xi_5\frac{\partial}{\partial\xi_5}+
\xi_3\xi_4\frac{\partial}{\partial\xi_5}
)$ \\
$S(1,4)$ & $(0|1,1,0,0)$ & $\xi_1
\frac{\partial}{\partial\xi_3}+
\xi_2\frac{\partial}{\partial\xi_4}+
\xi_1\xi_3\frac{\partial}{\partial x}
+\alpha\xi_1\xi_4\frac{\partial}{\partial x}+
\xi_2\xi_4(x\frac{\partial}{\partial x}+\xi_3\frac{\partial}{\partial \xi_3}), ~\alpha\in\F$\\
$S(2,0)$ & $(1,1|)$ & $x_1^2\frac{\partial}{\partial x_2}+
x_2^2\frac{\partial}{\partial x_1}$ \\
$S(2,0)$ & $(1,0|)$ & $x_1\frac{\partial}{\partial x_2}$ \\

\hline 
\end{tabular}
\end{center}
\begin{center} Table 1: even admissible gradings
\end{center}


\begin{center}
\begin{tabular}{|c|c|c|}
\hline
$S$ & type of grading & $\mu$~~($\alpha\in\F)$\\
\hline
\vspace{0.001cm} & \vspace{0.001cm} & \vspace{0.001cm}\\
$W(1,1)$ & $(1|0)$ & $x^2\xi\frac{\partial}{\partial x}+x\frac{\partial}
{\partial\xi}$ \\
$W(1,2)$ & $(0|1,1)$ & $\xi_1\frac{\partial}{\partial x}+
(\alpha+x)\xi_1\xi_2
\frac{\partial}{\partial\xi_1}$ \\
$W(2,1)$ & $(1,0|0)$ & $x_1^2\xi\frac{\partial}{\partial x_1}+
x_1\frac{\partial}{\partial\xi}+x_1\xi\frac{\partial}{\partial x_2}$ \\
$W(2,2)$ & $(0,0|1,1)$ & $\xi_1\frac{\partial}{\partial x_1}
+\xi_2\frac{\partial}{\partial x_2}+\xi_1\xi_2((1+x_1)
\frac{\partial}{\partial\xi_1}
+\alpha x_2\frac{\partial}{\partial\xi_2})$ \\
$S(2,2)$ & $(0,0|1,1)$ & $\xi_1\frac{\partial}{\partial x_1}+(1+\alpha x_1+x_1x_2)\xi_2\frac{\partial}{\partial x_2}
-x_1\xi_1\xi_2\frac{\partial}{\partial \xi_1}$ \\
$S(1,3)$ & $(0|1,1,1)$ & $\xi_3\frac{\partial}{\partial x}+
\xi_1\xi_2\frac{\partial}{\partial\xi_3}$ \\
$H(2,2)$ & $(1,0|1,0)$ & $(\alpha+q)p^2\xi_2+p\xi_1$ \\
$HO(2,2)$ & $(1,0|0,1)$ & $x_1^2+x_1\xi_1\xi_2$ \\
$HO(3,3)$ & $(0,0,0|1,1,1)$ & $\xi_1\xi_2+(\alpha+x_1)x_1\xi_1\xi_3$
 \\
$HO(n,n)$ ${}_{(n\geq 2)}$ & $(0,\dots,0,1|0,\dots,0,-1)$ & $x_n+x_n^2\xi_{n}\xi_{n-1}$ \\
$SHO(3,3)$ & $(1,0,0|0,1,1)$ & $x_1^2+\xi_2\xi_3+x_1\xi_3\xi_1-x_2\xi_3\xi_2$ \\
$SHO(4,4)$ & $(0,0,0,0|1,1,1,1)$ & $\xi_1\xi_2+(\alpha+x_1)\xi_3\xi_4$ \\
$SHO(n,n)$ ${}_{(n\geq 3)}$ &  $(0,\dots,0,1|0,\dots,0,-1)$ & $x_n+x_n^2\xi_n\xi_{n-1}+
2x_nx_{n-2}\xi_{n-1}\xi_{n-2}$ \\
$KO(n,n+1)$ ${}_{(n\geq 2)}$ & $(0,\dots,0,1|0,\dots,0,-1,0)$ & $x_n+x_n^2\xi_n\xi_{n-1}$ \\
$KO(2,3)$ & $(0,0|1,1,1)$ & $\xi_1\tau+x_1\xi_1\xi_2$\\
${\begin{array}{c}
SKO(n,n+1;\beta)\\ 
~~~~{}^{(n\geq 2, ~~\beta\neq 4/n)}
\end{array}}$ & $(0,\dots,0,1|0,\dots,0,-1,0)$ &
$x_n+x_n\xi_{n-1}(\tau-\frac{2-n\beta}{2}x_n\xi_{n})$ \\
${\begin{array}{c}
SKO(2,3;\beta)\\ 
~~~~{}^{(\beta\neq 0,1,2+\frac{1}{b}, 2+\frac{2}{b}, \forall b\in\Z_{>0})}
\end{array}}$ & $(0,1|0,-1,0)$ &
$x_2(1+x_2\xi_2\tau-x_{1}\xi_{1}(2\tau-(3-2\beta)x_2
\xi_2))$ \\
$SKO(3,4;1)$
 & $(0,0,1|0,0,-1,0)$ &
$x_3(1+\xi_1\xi_2+\xi_2\tau+\frac{1}{2}x_3\xi_2\xi_3)$\\
$SKO(3,4;\beta)$ & $(0,0,0|1,1,1,1)$ & $\xi_1\xi_2+
(\alpha+2x_1)(x_1\xi_1\xi_3-
\tau\xi_3)-3\beta(\alpha+x_1)x_1\xi_1\xi_3$\\
\hline 
\end{tabular}
\end{center}
\begin{center} Table 2: odd admissible gradings
\end{center}

\medskip


\begin{lemma}\label{uniquenessofp} Let $S=S_0\supset 
S_1\supset\dots$ be a filtered linearly compact Lie algebra
acting on a linearly compact filtered vector space 
$V=V_{-k}\supset V_{-k+1}\supset\dots\supset V_0\supset V_1\supset\dots$,
such that $S_j$ and $V_j$ have finite codimension and 
$S_i (V_j) \subset V_{i+j}$. Suppose that
\begin{equation}
\mbox{there exists}~ \mu\in V, ~\mbox{such that}~ S(\mu) +\F \mu=V.
\label{*}
\end{equation}
Then for any other such element $\mu'$
there exists an element  $a\in\exp(S)\times\F^{\times}$, such that 
$a(\mu)=\mu'$.
\end{lemma}
{\bf Proof.}  Condition (\ref{*}) is equivalent to the existence, for every 
$j$, of an element $\mu_j$
in $V/V_j$ such that $S(\mu_j)+\F \mu_j=V/V_j$, or, equivalently, such that
\begin{equation}
\s_j(\mu_j)+\F \mu_j=V/V_j,
\label{fdred}
\end{equation}
where $\s_j=S/S_{j+k}$.
If $\mu'$ is another element satisfying condition (\ref{*}) then there exist
elements $\mu'_j\in V/V_j$ such that $\s_j(\mu'_j)+\F \mu'_j=V/V_j$ . 
From (\ref{fdred}) it follows
that that the orbit of $\F\mu_j$ is Zariski open in $V/V_j$. Likewise, the 
orbit of $\F\mu'_j$ is Zariski open in $V/V_j$. 
Since any two Zariski open subsets in $V/V_j$ have a 
non-empty intersection, we conclude that 
for every $j$ there exists an element $a_j\in\exp(\s_j)\times\F^\times$ such 
that $a_j(\mu_j)=\mu'_j$ (see also the proof of \cite[Proposition 2]{K1}). 
Then the inverse limit $a$ of the $a_j$'s sends $\mu$ to
$\mu'$. \hfill$\Box$

\medskip

Recall that  the {\em canonical filtration} of a simple infinite-dimensional 
linearly compact Lie superalgebra $\g$ is the Weisfeiler filtration of $\g$, 
associated to the canonical subalgebra $L_0$ of $\g$, i.e., the  open maximal 
subalgebra of $\g$ equal to the intersection of all open  subalgebras of 
minimal codimension in $\g$ \cite[\S 11]{CantaK}. The subalgebra $L_0$ 
contains all $ad$-exponentiable elements of $\g$.

\begin{proposition}\label{indofp}  
Let $\g=\prod_{j\geq -1}\g_j$ be an even (resp.\ odd) admissible grading of a
simple linearly compact Lie superalgebra
$\g$ and let $\mu\in(\g_1)_{\0}$
(resp.\ $\mu\in(\g_1)_{\1}$) be as in Definition \ref{admissible}.  
If \begin{equation}
[\g_{int},\mu]+\F \mu=(\g_1)_{\0} ~~(\mbox{resp.}~ [\g_{int},\mu]+
\F \mu=(\g_1)_{\1}),
\label{a'}
\end{equation} then the  superalgebra $J(\g,\mu)$ 
does not depend on the choice of the element $\mu$. 
\end{proposition} 
{\bf Proof.} It is sufficient to use Lemma \ref{uniquenessofp} where
 $S=\g_{int}$, $V=(\g_1)_{\0}$ if $\mu$ is even or $V=(\g_1)_{\1}$ if
 $\mu$ is odd, and the filtrations of $S$ and $V$ are induced by the 
canonical filtration of $\g$. 
Since the canonical subalgebra $L_0$ of $\g$ 
contains all its  exponentiable elements,  $\g_{int}$ is
 contained in $L_0$, hence the filtration of $L_0$ induces a filtration of 
$\g_{int}$
satisfying the hypotheses of Lemma \ref{uniquenessofp}.
  \hfill$\Box$ 

\begin{lemma}\label{lemma} Let $\g=\prod_{j\geq -1}\g_j$
be a depth 1 $\Z$-graded linearly compact Lie superalgebra, and let
$\g=\prod_{j\geq -k} S_j$ be a $\Z$-grading of $\g$, such that 
$\dim(S_j)<\infty$  for all $j$, and $\g_{int}\subset \prod_{j\geq 0} S_j$.
Let $\mu=\mu_0+h$, where $\mu_0, h\in\g_1$, $\mu_0\in S_{-k}+\dots+S_t$, 
$h\in\prod_{j> t}S_j$, and
$[\g_{int},\mu_0]+(\g_1\cap\oplus_{j\leq t}S_j)=\g_1$ for some $t \in \Z_+$.
Then $[\g_{int},\mu]=[\g_{int},\mu_0]$.  
Consequently, by Lemma \ref{uniquenessofp}, $\mu$ and
$\mu_0$ are conjugate  by an element of $exp(ad(\g_{int}))$.
\end{lemma}
{\bf Proof.}
We have: $[\g_{int},\mu_0]= \sum_{i \in \Z_+}[\g_{int}\cap S_i,\mu_0]$
and $[\g_{int},\mu]=\sum_{i \in \Z_+}([\g_{int}\cap S_i,\mu_0]+
[\g_{int}\cap S_i,h])$,
where
$[\g_{int}\cap S_i,h]$ lies in $\oplus_{j>i+t}S_j$.
It follows that  $[\g_{int},\mu]=[\g_{int},\mu_0]$. 
\hfill$\Box$

\begin{example}\label{HO}\em  In Example \ref{PO(n,n)} we introduced
 the Buttin bracket $\{\cdot,\cdot\}$ on 
${\cal O}(n,n)$. Recall that it induces on the 
superspace $HO(n,n)={\cal O}(n,n)/\F 1$ with reversed parity a simple Lie
superalgebra bracket, denoted by $[\cdot, \cdot]$. Consider the Lie superalgebra
$\mathfrak{g}=HO(n+1,n+1)$ for $n\geq 1$, 
with the grading 
$\mathfrak{g}=\prod_{j\geq -1}\g_j$ of type $(0,\dots,0,1|0,\dots,0,-1)$ (see \S \ref{pre}).
Then we have: $\g_{-1}=\langle\xi_{n+1}\rangle\otimes
{\cal O}(n,n)$, $\g_{0}=(\langle1, x_{n+1}\xi_{n+1}\rangle\otimes
{\cal O}(n,n))/\F 1$, $\g_{1}=\langle x_{n+1}, x_{n+1}^2\xi_{n+1}\rangle\otimes
{\cal O}(n,n)$, with reversed parity, hence we can identify $\g_{-1}$, $\g_0$ and $\g_1$
with the following subspaces:
\begin{itemize}
\item[] $\g_{-1}\equiv {\cal O}(n,n)$,
\item[] $\g_0\equiv {\cal O}(n,n)/\F 1\oplus \eta {\cal O}(n,n)$ with reversed parity,
\item[] $\g_1\equiv {\cal O}(n,n)\oplus \xi {\cal O}(n,n)$ with reversed parity,
\end{itemize}

\noindent
where $\xi$ and $\eta$ are odd indeterminates. Under this identification,
the bracket of $\g_0$ and $\g_1$ can be written 
as follows ($f\in {\cal O}(n,n)/\F 1$, $g, s_1, s_2\in {\cal O}(n,n)$):
\begin{equation}
[f+\eta g, s_1+\xi s_2]=[f,s_1]+(-1)^{p(g)+1}gs_1+
\xi((-1)^{p(f)+1}[f,s_2]+[g,s_1]-gs_2),
\label{g0g1}
\end{equation}
and the bracket of $\g_{-1}$ and $\g_1$ as follows: if  ($t\in \g_{-1}$
has zero constant term, then 
\begin{equation}
[t, s_1+\xi s_2]=(-1)^{p(t)+1}ts_1+
\eta([t,s_1]-2ts_2);
\label{g-1g1}
\end{equation}
and 
\begin{equation}
[1, s_1+\xi s_2]=-2\eta s_2.
\label{g-1g1'}
\end{equation}

In order to show that the grading of type $(0,\dots,0,1|0,\dots,0,-1)$ is odd 
admissible, consider the element
 $\mu_0=1+\xi\xi_n\in(\g_1)_{\bar{1}}$. Then  (\ref{a}) follows from 
(\ref{g0g1}). 
Besides, (\ref{g-1g1}) and (\ref{g-1g1'}) show that
the Lie subalgebra generated by $[\mu_0,\g_{-1}]$ contains
the element $\eta\xi_n$ and all elements 
$t+2\eta\xi_nt$ with $t\in {\cal O}(n,n)$, $t\neq 1$, hence it contains the 
subalgebra $\eta {\cal O}(n,n)$,
since $[\eta\xi_n, t+2\eta\xi_n t]=-\eta\frac{\partial t}{\partial x_n}$,
hence it contains ${\cal O}(n,n)$.
Therefore $\mu_0$ and $\g_{-1}$ generate $\g_0$, hence
the grading of type 
$(0,\dots,0,1|0,\dots,0,-1)$ is odd admissible.

We have $J(\g,\mu_0)={\cal O}(n,n)$
with product ($f,g\in {\cal O}(n,n)$):
\begin{equation}
f\bullet g=(-1)^{p(f)+1}\{f,g\}+2\xi_nfg.
\label{product1}
\end{equation}
By Remark \ref{simplier}, $J(\g,\mu_0)\cong OJP(n-1,n-1)$. 


\medskip

We now want to show that the choice of $\mu\in(\g_1)_{\bar{1}}$ satisfying
(\ref{a}) and (\ref{b}) is unique up to automorphisms in $\exp(ad(\g_{int}))$.
Let $\mu=\mu_1+\xi \mu_2$, with $\mu_1,\mu_2\in {\cal O}(n,n)$, be an odd element in
$\g_1$ satisfying properties (\ref{a}) and (\ref{b}).
We shall first show that if $\mu$ satisfies (\ref{b}), then
$\mu_1$ has a non-zero constant term.
In order to construct the Lie subalgebra of $\g_0$
generated by $[\g_{-1},\mu]$, 
notice that, if $r\in {\cal O}(n,n)$, then, by (\ref{g-1g1}),
$[r, \mu]=(-1)^{p(r)+1}r\mu_1+\eta t'$ for some $t'\in {\cal O}(n,n)$.
Besides, if $r',r'',t',t''\in {\cal O}(n,n)$, then:
\begin{equation}
[r'\mu_1+\eta t', r''\mu_1+\eta t'']=[r'\mu_1,r''\mu_1]+\eta\omega
\label{ideal}
\end{equation}
for some $\omega\in {\cal O}(n,n)$. Let us denote by $I$ the principal ideal 
of ${\cal O}(n,n)$, with respect to
the usual associative product, generated by $\mu_1$.
Since $[\g_{-1},\mu]$ generate the whole Lie superalgebra $\g_0$, 
(\ref{ideal}) implies that the Lie superalgebra ${\cal O}(n,n)/\F 1$ is 
generated by the image $\bar{I}$ of $I$.
Let us now denote by $K$ the principal ideal of 
${\cal O}(n,n)$, 
with respect to the usual associative
product, generated by $[\mu_1,\mu_1]$, and denote by $\bar{K}$ its image in
${\cal O}(n,n)/\F 1$. We have: $[I,I]\subset I+K$ and
$[I,K]\subset I+K$ since $[\mu_1,[\mu_1,\mu_1]]=0$ ($\mu_1$ is odd). 
Finally, $[K,K]\subset K$,
hence $I+K$ is a subalgebra of the Lie superalgebra ${\cal O}(n,n)$,
 containing the Lie subalgebra,
generated by $I$, hence $\bar{I}+\bar{K}={\cal O}(n,n)/\F 1$.
Now suppose that $\mu_1$ has zero constant term, 
hence, up to a linear change of indeterminates, $\mu_1=ax_1+q_1$ where
$a\in\F$ and $q_1$ is a sum of monomials of degree greater than or equal to 2.
Then $[\mu_1,\mu_1]=2a\frac{\partial q_1}{\partial\xi_1}+[q_1,q_1]$, hence
$I+K=(x_1+q_1){\cal O}(n,n)+(2a\frac{\partial q_1}{\partial\xi_1}+[q_1,q_1]){\cal O}(n,n)$, therefore it
does not contain any even indeterminate
different from $x_1$, since $q_1$ is an even element of ${\cal O}(n,n)$.
This contradicts the equality $\bar{I}+\bar{K}={\cal O}(n,n)/\F 1$,
therefore $\mu_1$ is an invertible element in ${\cal O}(n,n)$. 

Up to an automorphism, we may thus assume
$\mu=1+\xi \mu_2$. Indeed, if $\mu_1=1+\mu'_1$, where $\mu'$ has zero constant 
term,then the automorphism 
$\exp(ad \eta \mu'_1)$ maps $\mu$ to an element of the form $1+\xi q_2$,
for some $q_2$ in ${\cal O}(n,n)$, which still satisfies both conditions 
(\ref{a}) and (\ref{b}).
By (\ref{g0g1}) and condition (\ref{b}), we thus have:
\begin{equation}
\langle \xi [f,\mu_2], ~g-(-1)^{p(g)+1}\xi g\mu_2~|~ f\in {\cal O}(n,n)/\F 1,
g\in {\cal O}(n,n)\rangle+\F \mu=\g_1, 
\label{g1spanned}
\end{equation}
hence $\langle [f,\mu_2]~|~ f\in {\cal O}(n,n)/\F1\rangle={\cal O}(n,n)$.
It follows that $\mu_2$ has a non-zero linear term, hence, up to a linear
change of indeterminates, since $\mu_2$ is
an odd element of ${\cal O}(n,n)$, we may assume $\mu_2=\xi_n+\rho$,
 where $\rho$ is a sum of monomials of degree greater than or equal to 2,
i.e.\ $\mu=1+\xi\xi_n+\xi \rho$.
In other words, if $\g=\prod_{j\geq -1}S_j$ denotes the principal grading of 
$\g$, then $\mu_0\in S_{-1}+S_2$ and $\mu=\mu_0+h$ with $h\in \prod_{j>2}S_j$.
 By Lemma \ref{lemma}, $\mu$ is conjugate to $\mu_0$.
\end{example}

\medskip

In what follows, if $\g$ is one of the
simple infinite-dimensional linearly compact Lie
superalgebras appearing in $i)$ and $ii)$ of Theorem \ref{eag},
we will denote by $\g=\prod_{j}S_j$ the principal grading,
defined in \cite{K3}.

\begin{theorem}\label{orbits} Let 
$\g=\prod_{j\geq -1}\g_j$ be a simple, linearly compact, 
infinite-dimensional Lie superalgebra with an admissible grading
of height strictly greater than one.
Then, up to conjugation by $exp(ad(\g_{int}))$, a complete list of
elements $\mu\in\g_1$ satisfying properties (\ref{a}) and (\ref{b}), is as follows:
\begin{enumerate}
\item $\g=W(1,1)$ with the grading of type $(1|0)$: $\mu=x^2\xi\frac{\partial}
{\partial x}+x\frac{\partial}{\partial \xi}$ and $[\g_{int},\mu]+\F\mu=\g_1$.
\item $\g=W(1,2)$ with the grading $(0|1,1)$: $\mu=\xi_1\frac{\partial}{\partial x}+(\alpha+x)\xi_1\xi_2\frac{\partial}
{\partial\xi_1}$, where either $\alpha\neq 0$
 and $[\g_{int},\mu]+\F\mu=\g_1$, or $\alpha=0$
 and $[\g_{int},\mu]+\F\mu$ has codimension one in $\g_1$. 
\item $\g=W(2,1)$ with the grading $(1,0|0)$: $\mu=x_1\frac{\partial}{\partial \xi}+
x_1\xi\frac{\partial}{\partial x_2}+
x_1^2\xi\frac{\partial}
{\partial x_1}$
and $[\g_{int},\mu]=\g_1$. 
\item $\g=W(2,2)$ with the grading $(0,0|1,1)$: 
$\mu=\xi_1\frac{\partial}{\partial x_1}+\xi_2\frac{\partial}{\partial x_2}+
(1+x_1)\xi_1\xi_2\frac{\partial}{\partial\xi_1}+\alpha x_2\xi_1\xi_2\frac{\partial}{\partial\xi_2}$, where either $\alpha\neq 0$ and $[\g_{int},\mu]+\F\mu=\g_1$, or $\alpha=0$
 and $[\g_{int},\mu]+\F\mu$ has
codimension one in $\g_1$.
\item $\g=S(1,4)$ with the grading $(0|1,1,0,0)$:
$\mu=\xi_1\frac{\partial}{\partial \xi_3}+
\xi_2\frac{\partial}{\partial\xi_4}+
\alpha\xi_1\xi_4\frac{\partial}
{\partial x}+
\xi_1\xi_3\frac{\partial}
{\partial x}+
x\xi_2\xi_4\frac{\partial}
{\partial x}-
\xi_2\xi_3\xi_4\frac{\partial}{\partial\xi_3}$,
where either $\alpha\neq 0$ 
and $[\g_{int},\mu]+\F\mu=\g_1$, or $\alpha=0$
and $[\g_{int},\mu]+\F\mu$ has codimension one in $\g_1$. 
\item $\g=S(2,0)$ with the grading $(1,1|)$: $\mu=x_1^2\frac{\partial}{\partial x_2}+x_2^2\frac{\partial}{\partial x_1}$ and $[\g_{int},\mu]+\F\mu=\g_1$.
\item $\g=S(2,0)$ with the grading $(1,0|)$: $\mu=x_1\frac{\partial}{\partial x_2}$ and $[\g_{int},\mu]+\F\mu=\g_1$.
\item $\g=S(2,2)$ with the grading $(0,0|1,1)$:
$\mu=\xi_1\frac{\partial}{\partial x_1}+\xi_2\frac{\partial}{\partial x_2}+
(\alpha x_1+x_1x_2)\xi_2\frac{\partial}{\partial x_2}-x_1\xi_1\xi_2\frac{\partial}{\partial\xi_1}$, where either $\alpha\neq 0$ and $[\g_{int},\mu]+\F\mu=\g_1$, or $\alpha=0$ and $[\g_{int},\mu]+\F\mu$ has
codimension two in $\g_1$.
\item $\g=S(1,3)$ with the grading $(0|1,1,1)$: $\mu=\xi_3\frac{\partial}
{\partial x}+\xi_1\xi_2\frac{\partial}{\partial\xi_3}$ and
$[\g_{int},\mu]+\F\mu=\g_1$. 
\item $\g=H(2,2)$ with the grading $(1,0|1,0)$: $\mu=p\xi_1+(\alpha+q)p^2\xi_2$, where either $\alpha\neq 0$
and $[\g_{int},\mu]+\F\mu=\g_1$, or $\alpha=0$
and $[\g_{int},\mu]+\F\mu$ has codimension one in $\g_1$.
\item $\g=HO(2,2)$ with the grading $(1,0|0,1)$:
$\mu=x_1^2+x_1\xi_1\xi_2$ and $[\g_{int},\mu]+\F\mu=\g_1$.
\item $\g=HO(3,3)$ with the grading $(0,0,0|1,1,1)$: 
$\mu=\xi_1\xi_2+(\alpha+x_1)x_1\xi_1\xi_3$, where $\alpha\neq 0$ and $[\g_{int},\mu]+\F\mu=\g_1$,
or $\alpha=0$ and $[\g_{int},\mu]+\F\mu$ has
codimension one in $\g_1$.
\item $\g=HO(n,n)$, $n\geq 2$, with the grading $(0,\dots,0,1|0,\dots,0,-1)$:
$\mu=x_n+x_n^2\xi_n\xi_{n-1}$ and $[\g_{int},\mu]+\F\mu=\g_1$.
\item $\g=SHO(3,3)$ with the grading $(1,0,0|0,1,1)$: $\mu=x_1^2+
\xi_2\xi_3+x_1\xi_3\xi_1-x_2\xi_3\xi_2$ and $[\g_{int},\mu]+\F\mu=\g_1$. 
\item $\g=SHO(4,4)$ with the grading $(0,0,0,0|1,1,1,1)$:
$\mu=\xi_1\xi_2+(\alpha+x_1)\xi_3\xi_4$, where either $\alpha\neq 0$ and
$[\g_{int},\mu]+\F\mu=\g_1$, or $\alpha=0$ and
$[\g_{int},\mu]+\F\mu$ has codimension 1 in $\g_1$.
\item $\g=SHO(n,n)$, $n\geq 3$, with the grading $(0,\dots,0,1|0,\dots,0,-1)$:
$\mu=x_n+x_n^2\xi_n\xi_{n-1}+2x_nx_{n-2}\xi_{n-1}\xi_{n-2}$
and $[\g_{int},\mu]+\F\mu=\g_1$. 
\item $\g=KO(n,n+1)$, $n\geq 2$, with the grading $(0,\dots,0,1|0,\dots,0,-1,0)$:
$\mu=x_n+x_n^2\xi_n\xi_{n-1}$ and $[\g_{int},\mu]+\F\mu=\g_1$. 
\item $\g=KO(2,3)$ with the grading $(0,0|1,1,1)$: 
$\mu=\xi_1\tau+x_1\xi_1\xi_2$ and $[\g_{int},\mu]+\F\mu$ has codimension one 
in $\g_1$.
\item $\g=SKO(n,n+1;\beta)$, $n\geq 2$, $\beta\neq \frac{4}{n}$,
with the grading $(0,\dots,0,1|0,\dots,0,-1,0)$: 
$\mu=x_n+x_n\xi_{n-1}\tau+\frac{2-n\beta}{2}x_n^2\xi_n\xi_{n-1}$
and $[\g_{int},\mu]+\F\mu=\g_1$ if $(n,\beta)\neq (3,1)$,
$[\g_{int},\mu]+\F\mu$ has codimension one in $\g_1$ if $n=3$ and $\beta=1$.
\item $\g=SKO(2,3;\beta)$, $\beta\neq 0,1, 2+\frac{1}{b}, 2+\frac{2}{b}$
$\forall b\in\Z_{>0}$, with the grading $(0,1|0,-1,0)$:
$\mu=x_2+
x_2^2\xi_{2}\tau-2x_{1}x_2\xi_{1}\tau+(3-2\beta)x_{1}x_2^2
\xi_{1}\xi_2$  and $[\g_{int},\mu]+\F\mu$ has codimension one in $\g_1$.
\item $\g=SKO(3,4;1)$, with the grading $(0,0,1|0,0,-1,0)$:
$\mu=x_3+x_3\xi_1\xi_2+x_3\xi_2\tau+\frac{1}{2}x_3^2\xi_2\xi_3$
and $[\g_{int},\mu]+\F\mu=\g_1$.
\item $\g=SKO(3,4;\beta)$, with the grading $(0,0,0|1,1,1,1)$:
$\mu=\xi_1\xi_2+\alpha((1-3\beta)x_1\xi_1\xi_3-\tau\xi_3)+
(2-3\beta)x_1^2\xi_1\xi_3-2x_1\tau\xi_3$, where either $\alpha\neq 0$ and 
$\beta\neq 1/3$ and $[\g_{int},\mu]=\g_1$, or $\alpha=0$ and 
$[\g_{int},\mu]$ has codimension one in $\g_1$.
\end{enumerate}

\end{theorem}
{\bf Proof.}  If $\g=W(1,1)$ with the grading of type $(1|0)$, then
one checks by direct computation that $\mu=x\frac{\partial}{\partial \xi}+
x^2\xi\frac{\partial}{\partial x}$ satisfies properties (\ref{a})
and (\ref{b}). In this case
$(\g_0)_{\0}=\g_{int}$, hence 1.\ follows from Proposition \ref{indofp}.
The same argument holds in case 6.
In all other cases we use the same kind of argument as
in Example \ref{HO}, where we prove 12.
For instance,
let $\g=S(2,0)$ with the grading of type $(1,0|)$. Then 
$\g_{-1}=\langle\frac{\partial}{\partial x_1}\rangle\otimes\F[[x_2]]$,
$\g_0=\langle x_2^r\frac{\partial}{\partial x_2}-rx_1x_2^{r-1}
\frac{\partial}{\partial x_1}~|~ r\geq 0\rangle$ and
$\g_1=\langle 2x_1x_2^r\frac{\partial}{\partial x_2}-rx_1^2x_2^{r-1}
\frac{\partial}{\partial x_1}~|~ r\geq 0\rangle$.
Let $\mu=\sum_{r\geq 0}\alpha_r
(2x_1x_2^r\frac{\partial}{\partial x_2}-rx_1^2x_2^{r-1}
\frac{\partial}{\partial x_1})\in\g_1$ satisfy properties
(\ref{a}) and (\ref{b}). If $\alpha_0=0$, then
the subalgebra of $\g_0$ generated by $[\g_{-1},\mu]$ is contained
in $S_0$ hence it is properly contained in $\g_0$, contradicting
 property (\ref{b}). It follows that $\alpha_0\neq 0$, hence, by Lemma
\ref{lemma}, $\mu$ is conjugate to $x_1\frac{\partial}{\partial x_2}$,
hence $[\g_{int},\mu]+\F\mu=\g_1$.

Let us now consider $\g=W(1,2)=\prod_{j\geq -1}\g_j$ with the grading of type $(0|1,1)$.
We have:
$\g_{-1}=\langle\frac{\partial}{\partial\xi_i}~|~i=1,2\rangle\otimes\F[[x]]$,
$\g_0=\langle\frac{\partial}{\partial x},
\xi_j\frac{\partial}{\partial\xi_i}~|~j,i=1,2\rangle\otimes\F[[x]]$,
$\g_1=\langle\xi_i\frac{\partial}{\partial x},
\xi_1\xi_2\frac{\partial}{\partial\xi_i}~|~i=1,2\rangle\otimes\F[[x]]$.
Let $\mu\in\g_1$ satisfy properties (\ref{a}) and (\ref{b}).
If $\mu\in\g_1\cap \prod_{j\geq 1}S_j$, then the subalgebra of $\g_0$
generated by $[\g_{-1},\mu]$ is contained in $\prod_{j\geq 0}S_j$, hence it
is properly contained in $\g_0$. It follows that, 
 up to a linear change of indeterminates,
we may assume $\mu=\xi_1\frac{\partial}{\partial x}+h$ for some
$h\in\g_1\cap\prod_{j\geq 1}S_j$.
We will show that, up to conjugation by elements in $\exp(ad(\g_{int}))$,
we may assume  $h=\alpha\xi_1\xi_2\frac{\partial}{\partial \xi_1}+\varphi$,
for some $\alpha\in\F$ and some $\varphi\in\prod_{j\geq 2}S_j$.
Indeed, let 
$\mu=\xi_1\frac{\partial}{\partial x}+a\xi_1\xi_2\frac{\partial}{\partial\xi_1}
+b\xi_1\xi_2\frac{\partial}{\partial\xi_2}+cx\xi_1\frac{\partial}{\partial x}
+dx\xi_2\frac{\partial}{\partial x}+h_1$ for some
$a, b, c, d\in\F$ and some $h_1\in\g_1\cap\prod_{j\geq 2}S_j$.
Then $\exp(ad(-dx\xi_2\frac{\partial}{\partial\xi_1}))(\mu)=
\xi_1\frac{\partial}{\partial x}+(a+d)\xi_1\xi_2\frac{\partial}{\partial\xi_1}
+b\xi_1\xi_2\frac{\partial}{\partial\xi_2}+cx\xi_1\frac{\partial}{\partial x}
+h'_1$ for some
$h'_1\in\g_1\cap\prod_{j\geq 2}S_j$.
It follows that we may assume
$\mu=\xi_1\frac{\partial}{\partial x}+a\xi_1\xi_2\frac{\partial}{\partial\xi_1}
+b\xi_1\xi_2\frac{\partial}{\partial\xi_2}+cx\xi_1\frac{\partial}{\partial x}
+h_1$ for some
$h_1\in\g_1\cap\prod_{j\geq 2}S_j$.
Then $\exp(ad(bx\xi_2\frac{\partial}{\partial\xi_2}))(\mu)=
\xi_1\frac{\partial}{\partial x}+a\xi_1\xi_2\frac{\partial}{\partial\xi_1}
+cx\xi_1\frac{\partial}{\partial x}
+h''_1$ for some
$h''_1\in\g_1\cap\prod_{j\geq 2}S_j$.
It follows that we may assume
$\mu=\xi_1\frac{\partial}{\partial x}+a\xi_1\xi_2\frac{\partial}{\partial\xi_1}
+cx\xi_1\frac{\partial}{\partial x}
+h_1$ for some
$h_1\in\g_1\cap\prod_{j\geq 2}S_j$.
Then
$\exp(ad(\frac{c}{2}x^2\frac{\partial}{\partial x}))(\mu)=
\xi_1\frac{\partial}{\partial x}+a\xi_1\xi_2\frac{\partial}{\partial\xi_1}
+\tilde{h}$ for some
$\tilde{h}\in\g_1\cap\prod_{j\geq 2}S_j$.
It follows that we may assume
 $\mu=\xi_1\frac{\partial}{\partial x}+a\xi_1\xi_2\frac{\partial}
{\partial\xi_1}+\tilde{h}$ for some $a\in\F$ and
$\tilde{h}\in\g_1\cap\prod_{j\geq 2}S_j$. If $a\neq 0$, then we may assume, up
to a linear change of indeterminates, $a=1$. 
By direct computation one shows that   
$\mu_0=\xi_1\frac{\partial}{\partial x}+\xi_1\xi_2\frac{\partial}
{\partial\xi_1}$ satisfies conditions (\ref{a}) and (\ref{b}), and that
 $[\g_0,\mu_0]=[\g_{int},\mu_0]$, since $[\frac{\partial}{\partial x},
\mu_0]=0$. Therefore, by Lemma \ref{lemma}, for every $a\neq 0$,  
$[\g_{int},\mu]=[\g_{int},\mu_0]$, hence 
$\mu$ is conjugate to $\mu_0$ by an element in $\exp(ad(\g_{int}))$.
Let us now suppose $a=0$. Then,
 by the same arguments as above, we can assume
either $(i)$
$\mu=\xi_1\frac{\partial}{\partial x}+
x\xi_1\xi_2\frac{\partial}{\partial\xi_1}+
\tilde{h}$ for some
$\tilde{h}\in\g_1\cap\prod_{j\geq 3}S_j$, or
$(ii)$
$\mu=\xi_1\frac{\partial}{\partial x}+
\tilde{h}$ for some
$\tilde{h}\in\g_1\cap\prod_{j\geq 3}S_j$.
Notice that $\mu_1=\xi_1\frac{\partial}{\partial x}+x\xi_1\xi_2\frac{\partial}
{\partial\xi_1}$ satisfies
conditions (\ref{a}) and (\ref{b}) and
$[\g_{int},\mu]$ has codimension one in $[\g_0,\mu]$. 
Therefore, by Lemma \ref{lemma},
in case $(i)$  $[\g_{int},\mu]=[\g_{int},\mu_1]$, hence
$\mu$ is conjugate to $\mu_1$.
Finally, in case $(ii)$, $[\g_{int},\mu]+\F\mu$ has
at least codimension 2 in $\g_1$, and this contradicts
condition $[\g_0,\mu]+\F\mu=\g_1$, since $\g_{int}$ has codimension 1
in $\g_0$. 

\medskip

Now let $\g=S(1,4)=\prod_{j\geq -1}\g_j$ with the grading of type $(0|1,1,0,0)$. We have: $\g_{-1}=\langle \frac{\partial}{\partial\xi_i}~|~i=1,2
\rangle\otimes\F[[x]]\otimes\Lambda(\xi_3,\xi_4)$,
$\g_0=\langle \frac{\partial}{\partial x}, 
\frac{\partial}{\partial\xi_k},
\xi_j\frac{\partial}{\partial\xi_i}~|~$
$k=3,4, i,j=1,2
\rangle\otimes\F[[x]]\otimes\Lambda(\xi_3,\xi_4)$,
$\g_1=\langle \xi_i\frac{\partial}{\partial x}, 
\xi_i\frac{\partial}{\partial\xi_k},
\xi_1\xi_2\frac{\partial}{\partial\xi_i}~|~k=3,4, i=1,2
\rangle\otimes\F[[x]]\otimes\Lambda(\xi_3,\xi_4)$.

Let $\mu\in(\g_1)_{\0}$ satisfy properties (\ref{a}) and (\ref{b}).
If $\mu\in\g_1\cap \prod_{j\geq 1}S_j$, then the subalgebra of $\g_0$
generated by $[\g_{-1},\mu]$ is contained in $\prod_{j\geq 0}S_j$, hence it
is properly contained in $\g_0$. It follows that, 
 up to a linear change of indeterminates,
we may assume $\mu=\xi_1\frac{\partial}{\partial\xi_3}+\xi_2\frac{\partial}
{\partial\xi_4}+h$ for some
$h\in\g_1\cap\prod_{j\geq 1}S_j$. 
We will show that, up to conjugation by elements in $\exp(ad(\g_{int}))$,
we may assume either $h=\alpha\xi_1\xi_4\frac{\partial}{\partial x}+
\xi_1\xi_3\frac{\partial}{\partial x}+
\varphi$ for some $\alpha\in \F$,
$\alpha\neq 0$
and some $\varphi\in\prod_{j\geq 2}S_j$, or
$h=\xi_1\xi_3\frac{\partial}{\partial x}+
x\xi_2\xi_4\frac{\partial}{\partial x}
-\xi_2\xi_3\xi_4\frac{\partial}{\partial\xi_3}+\psi$ for some
$\psi\in\prod_{j\geq 3}S_j$.
Indeed, let $\mu=\xi_1\frac{\partial}{\partial\xi_3}+\xi_2
\frac{\partial}{\partial\xi_4}+\alpha\xi_1\xi_4\frac{\partial}{\partial x}+
\beta\xi_1\xi_3\frac{\partial}{\partial x}+\gamma\xi_2\xi_4
\frac{\partial}{\partial x}+\delta\xi_2\xi_3\frac{\partial}{\partial x}+
Ax\xi_1\frac{\partial}{\partial \xi_3}+
Bx\xi_2\frac{\partial}{\partial \xi_3}+
Cx\xi_1\frac{\partial}{\partial \xi_4}+
Dx\xi_2\frac{\partial}{\partial \xi_4}+h_1$ for some
$\alpha, \beta, \gamma, \delta, A, B, C, D\in\F$ and some $h_1\in\g_1\cap\prod_{j\geq 2}S_j$.
Then $\exp(ad(-Ax\xi_1\frac{\partial}{\partial\xi_1}))(\mu)=
\mu-Ax\xi_1\frac{\partial}{\partial\xi_3}+h'_1$
for some
$h'_1\in\g_1\cap\prod_{j\geq 2}S_j$. It follows that we may assume $A=0$.
Then $\exp(ad(-Bx\xi_2\frac{\partial}{\partial\xi_1}))(\mu)=
\mu-Bx\xi_2\frac{\partial}{\partial\xi_3}+h''_1$
for some
$h''_1\in\g_1\cap\prod_{j\geq 2}S_j$. It follows that we may assume $B=0$.
By similar arguments we may assume $C=D=0$.
Then $\exp(ad(-\delta\xi_3\xi_4\frac{\partial}{\partial x}))(\mu)=
\mu+\delta\xi_1\xi_4\frac{\partial}{\partial x}-\delta\xi_2\xi_3\frac{\partial}
{\partial x}$ hence we may assume $\delta=0$, i.e.,
$\mu=\xi_1\frac{\partial}{\partial\xi_3}+\xi_2\frac{\partial}
{\partial\xi_4}+\alpha\xi_1\xi_4\frac{\partial}{\partial x}+
\beta\xi_1\xi_3\frac{\partial}{\partial x}+\gamma\xi_2\xi_4
\frac{\partial}{\partial x}+h_1$.
Notice that, since $[\g_{-1},\mu]$ generates $\g_0$ and 
$\frac{\partial}{\partial x}$ lies in $\g_0$, we necessarily have
$(\alpha,\beta,\gamma)\neq (0,0,0)$.

By direct computation one shows that   
$\mu_{\alpha,\beta,\gamma}=\xi_1\frac{\partial}{\partial \xi_3}+
\xi_2\frac{\partial}{\partial\xi_4}+
\alpha\xi_1\xi_4\frac{\partial}
{\partial x}+
\beta\xi_1\xi_3\frac{\partial}
{\partial x}+
\gamma\xi_2\xi_4\frac{\partial}
{\partial x}$ satisfies conditions (\ref{a}) and (\ref{b})
 for every $\alpha, \beta, \gamma\in \F$ such that
$\alpha\neq 0$ or $\alpha=0$ and $\beta\neq 0\neq\gamma$.
Besides, $[\g_0,\mu_{\alpha,\beta,\gamma}]=[\g_{int},\mu_{\alpha, \beta,
\gamma}]$, since $[\frac{\partial}{\partial x},
\mu_{\alpha,\beta, \gamma}]=0$. It follows that the elements 
$\mu_{\alpha,\beta,\gamma}$,  with $\alpha, \beta, \gamma\in \F$ such that
$\alpha\neq 0$ or $\alpha=0$ and $\beta\neq 0\neq\gamma$, are conjugate 
to each other by automorphisms in $exp(ad(\g_{int}))$.
Therefore,
if $\alpha\neq 0$ or $\alpha=0$ and $\beta\neq 0\neq \gamma$,
 by Lemma \ref{lemma}, $[\g_{int},\mu]=[\g_{int},\mu_{\alpha,
\beta,\gamma}]=[\g_{int}, \mu_{a,1,0}]$ for $a\neq 0$, hence, under these
hypotheses, 
$\mu$ is conjugate to $\mu_{a,1,0}$ by an element in $\exp(ad(\g_{int}))$.

Now suppose $\alpha=0=\gamma$ and $\beta\neq 0$
(the case $\alpha=0=\beta$ and $\gamma\neq 0$ is equivalent
to this up to a linear change of indeterminates). Then we may assume
$\mu=\xi_1\frac{\partial}{\partial\xi_3}+
\xi_2\frac{\partial}{\partial\xi_4}+
\xi_1\xi_3\frac{\partial}{\partial x}+q$ for some $q\in\prod_{j\geq 2}S_j$.

Notice that $\tilde{\mu}_{a,b,c}=
\xi_1\frac{\partial}{\partial \xi_3}+
\xi_2\frac{\partial}{\partial\xi_4}+
\xi_1\xi_3\frac{\partial}
{\partial x}+
ax\xi_2\xi_4\frac{\partial}
{\partial x}+b\xi_1\xi_2\xi_4\frac{\partial}{\partial\xi_1}+
c\xi_2\xi_3\xi_4\frac{\partial}{\partial\xi_3}$
satisfies
conditions (\ref{a}) and (\ref{b})
 for every
$a, b, c\in\F$ such that $a\neq 0$ and $a-b+c=0$  but in this case
 $[\g_{int},\tilde{\mu}_{a,b,c}]+\F\tilde{\mu}_{a,b,c}$ has codimension one in $[\g_0,\tilde{\mu}_{a,b,c}]+\F\tilde{\mu}_{a,b,c}$. 
More precisely, $[\g_{int},\tilde{\mu}_{a,b,c}]\oplus\F\xi_2\xi_4\frac{\partial}{\partial x}=
[\g_0,\tilde{\mu}_{a,b,c}]$. It follows that,  for 
$a, b, c\in\F$ such that $a\neq 0$ and $a-b+c=0$, the elements
$\tilde{\mu}_{a,b,c}$ are conjugate to each other by automorphisms 
in $exp(ad(\g_{int}))$.

One verifies that the vector field $\xi_2\xi_4\frac{\partial}{\partial x}$
lies in $[\g_0,\mu]+\F\mu$ if and only if 
$q=ax\xi_2\xi_4\frac{\partial}{\partial x}+b\xi_1\xi_2\xi_4\frac{\partial}
{\partial\xi_1}+c\xi_2\xi_3\xi_4\frac{\partial}{\partial\xi_3}+h+t$
for some $a,b,c\in\F$, with $a\neq 0$ and $a-b+c=0$, some 
$t\in\prod_{j\geq 3}S_j$ and some
$h\in\langle x\xi_1\xi_3\frac{\partial}{\partial x},
x\xi_1\xi_4\frac{\partial}{\partial x}, x\xi_2\xi_3\frac{\partial}{\partial x},
\xi_1\xi_3\xi_4\frac{\partial}{\partial x}, x^2\xi_i\frac{\partial}
{\partial\xi_4}, \xi_1\xi_2\xi_k\frac{\partial}{\partial\xi_i}~|~i=1,2,
k=3,4\rangle$.
It follows that $[\g_{int},\mu]+\F\mu$ consists of vector fields
in $\g_1$ not containing the element $\xi_2\xi_4\frac{\partial}{\partial x}$.
Therefore $\mu$ is conjugate to $\tilde{\mu}_{a,b,c}$ for
$a\neq 0$ and $a-b+c=0$, i.e., to $\tilde{\mu}_{1,0,-1}$.

\medskip

Let us now consider
$\g=W(2,1)=\prod_{j\geq -1}\g_j$ with the grading of type $(1,0|0)$.
We have: $\g_{-1}=\langle \frac{\partial}{\partial x_1}
\rangle\otimes A$,
$\g_0=\langle x_1\frac{\partial}{\partial x_1}, 
\frac{\partial}{\partial x_2},
\frac{\partial}{\partial\xi}\rangle\otimes A$,
$\g_1=\langle x_1^2\frac{\partial}{\partial x_1}, 
x_1\frac{\partial}{\partial x_2}$,
$x_1\frac{\partial}{\partial\xi}\rangle\otimes A$,
where $A=\F[[x_2]]\otimes\Lambda(\xi)$.

Let $\mu\in(\g_{1})_{\bar{1}}$ satisfy properties
(\ref{a}) and (\ref{b}). If $\mu\in\g_1\cap\prod_{j\geq 1}S_j$, then the
subalgebra of $\g_0$ generated by  $[\g_{-1},\mu]$ is contained in $S_0$,
hence it is properly contained in $\g_0$. It follows that, up to
a linear change of indeterminates, we may assume $\mu=x_1\frac{\partial}
{\partial\xi}+h$ for some $h\in\prod_{j\geq 1}S_j$.
Let $\mu=x_1\frac{\partial}{\partial\xi}+ax_1\xi\frac{\partial}{\partial x_2}+
bx_1x_2\frac{\partial}{\partial \xi}+h'$,
for some $a, b\in\F$ and some $h'\in\prod_{j\geq 2}S_j$.
Since the subalgebra generated by $[\g_{-1},\mu]$ contains $\frac{\partial}
{\partial x_2}$, we have $a\neq 0$. Besides,
$\exp(ad(-bx_1x_2\frac{\partial}{\partial x_1}))(\mu)=
\mu-bx_1x_2\frac{\partial}{\partial\xi}+h''$, for some
$h''\in\prod_{j\geq 2}S_j$. Therefore we can assume
$\mu=x_1\frac{\partial}{\partial\xi}+x_1\xi\frac{\partial}{\partial x_2}+
\tilde{h}$ for some $\tilde{h}\in\prod_{j\geq 2}S_j$.
By direct computation one shows that 
if $\mu_0=x_1\frac{\partial}{\partial \xi}+
x_1\xi\frac{\partial}{\partial x_2}$, then $[\g_0,\mu_0]=
[\g_{int},\mu_0]=\g_1$. By Lemma \ref{lemma}, $\mu$ is conjugate to
$\mu_0$.
Note that the subalgebra of $\g_0$ generated by
 $[\g_{-1},\mu_0]$  is properly contained in $\g_0$ since $\mu_0$ and all
vector fields in $\g_{-1}$ have zero divergence.
Let $\mu_1=\mu_0+x_1^2\xi\frac{\partial}{\partial x_1}$. Then
$[\g_{int},\mu_1]=[\g_{int},\mu_0]$ by Lemma \ref{lemma}, and a direct
computation shows that $[\g_{-1},\mu_0]$ generates $\g_0$.
It follows that $\mu$ is conjugate to
$\mu_1$.

By similar arguments one proves cases 6.--21.
 \hfill$\Box$

\bigskip

\begin{remark}\em We recall that the Lie superalgebras $S(1,2)$ and
$SKO(2,3;1)$ have an $sl_2$-copy $\mathfrak{a}$ of outer derivations.
We shall denote by $e$, $f$, $h$ the standard
basis of $\mathfrak{a}$ defined in \cite[Lemma 5.9]{FK} and
\cite[Remark 4.15]{CantaK}, respectively.
Besides, $Der SHO(3,3)=SHO(3,3)+\mathfrak{a}$ with $\mathfrak{a}\cong gl_2$
\cite[Proposition 1.8]{CantaK}.
The subalgebra $\mathfrak{a}$ of outer derivations of $SHO(3,3)$ is
generated by the Euler operator $E$ and by a copy of $sl_2$ with Chevalley
basis $\{e,h.f\}$ where
$$e=ad(\xi_1\xi_3\frac{\partial}{\partial x_2}
-\xi_2\xi_3\frac{\partial}{\partial x_1}-\xi_1\xi_2\frac{\partial}{\partial x_3}) ~~~\mbox{and}~~~ h=-ad\sum_{i=1}^3\xi_i\frac{\partial}{\partial\xi_i}.$$
(Note that the formula for $h$, given in \cite[Remark 2.37]{CantaK},
is incorrect.) In order to define $f$, consider $SHO(3,3)$ with its principal grading, i.e., the grading of type $(1,1,1|1,1,1)$. In this grading,
 for every $j>1$, $SHO(3,3)_j=SHO(3,3)_1^j$, therefore
it is sufficient to define $f$ on the local part 
$SHO(3,3)_{-1}\oplus SHO(3,3)_0\oplus SHO(3,3)_1$ of $SHO(3,3)$. One has:
$f(\xi_1\xi_2)=x_3$, $f(\xi_1\xi_3)=-x_2$, 
$f(\xi_2\xi_3)=x_1$;
$f(x_1\xi_2\xi_3)=\frac{1}{2}x_1^2$,
$f(x_2\xi_1\xi_3)=-\frac{1}{2}x_2^2$,
$f(x_3\xi_1\xi_2)=\frac{1}{2}x_3^2$,
$f(x_1\xi_1\xi_2-x_3\xi_3\xi_2)=x_1x_3$,
$f(x_2\xi_2\xi_1-x_3\xi_3\xi_1)=-x_2x_3$,
$f(x_1\xi_1\xi_3-x_2\xi_2\xi_3)=-x_1x_2$, 
and $f=0$ elsewhere on $SHO(3,3)_{-1}\oplus SHO(3,3)_0\oplus SHO(3,3)_1$.
\end{remark}

\begin{theorem}\label{derivations} Let
$L=S\rtimes\F d$, where $S$ is a simple
linearly compact Lie superalgebra and $d$ is an even outer derivation
of $S$. Then 
an even admissible grading of $L$ such that $d$ has degree 1, is either 
short or isomorphic to
one of the following:
\begin{enumerate}
\item $S=SHO(3,3)$ , $(1,0,0|0,1,1)$, $d=e$;
\item $S=SKO(2,3;1)$ , $(1,0|0,1,1)$ , $d=e$.
\end{enumerate}
\end{theorem}
{\bf Proof.} Let $L=\prod_{j\geq -1} L_j$ be a $\Z$-grading of $L$ of depth 1.
According to \cite[Proposition 5.1.2]{K2}, \cite{K3}
and \cite[Proposition 1.8]{CantaK}, 
one of the following possibilities may occur:
\begin{enumerate}
\item[i)] $S=A(1,1)$, $d=D_1$ (see \cite[Proposition 5.1.2(e)]{K2});
\item[ii)] $S=H(0,n)$ for some even $n> 5$ and $d=\xi_1\dots\xi_n$;
\item[iii)] $S=S(1,n)$ for some even $n\geq 4$ and 
$d=\xi_1\dots\xi_n\frac{\partial}{\partial x}$;
\item[iv)] $S=SKO(m,m+1;(m-2)/m)$ for some odd 
$m\geq 3$ and $d=\xi_1\dots\xi_m$;
\item[v)] $S=SKO(m,m+1;1)$ for some even $m\geq 4$ and $d=\tau\xi_1\dots\xi_m$;
\item[vi)] $S=SHO(m,m)$ for some odd $m\geq 5$ and $d=\xi_1\dots\xi_m$;
\item[vii)] $S=S(1,2)$ and $d=e$; 
\item[viii)] $S=SHO(3,3)$ and $d=e$;
\item[ix)] $S=SKO(2,3;1)$ and $d=e$.
\end{enumerate}
Let us analyze all cases $i)-ix)$.
From \cite[Proposition 5.1.2(e)]{K2} one deduces that if 
$S=A(1,1)$ and $L=S\rtimes\F D_1$, then any $\Z$-grading of $L$ of depth 1,
such that $D_1$ has degree 1, is short.

Let $S=H(0,n)$ for some even $n>5$, say $n=2t$. By
Proposition \ref{listdepth1},
the gradings of depth 1 of $L$ are,
up to isomorphism, 
the gradings of type $(|1,\dots,1)$ and the grading of type  
$(|1,\dots,1,0,\dots,0)$ with $t$ 1's and $t$ 0's. 
In the grading of type $(|1,\dots,1)$, $d$ has degree $n-2>1$  
(since $n\geq 5$) in the grading of type  
$(|1,\dots,1,0,\dots,0)$ it has degree $t-1>1$. Therefore if $L=H(0,n)\rtimes
\F\xi_1\dots\xi_n$ for some even $n>5$, then $L$ has no $\Z$-grading
of depth one such that 
$\xi_1\dots\xi_n$ has degree 1.

Let us consider case iii). By Proposition \ref{listdepth1},
the gradings of depth 1 of $L$ such that $d$ has degree 1, are,
up to isomorphism, 
the gradings of type $(0|1,0,\dots,0)$ and
$(1|1,1,0,\dots,0)$. The first one is short, therefore
let us consider the second one. The 0th and 1st graded components of $S$ in this
grading have dimension $(n+6)2^{n-2}$ and $(3+3n)2^{n-2}$, respectively.
Hence, by Remark \ref{necessarycond}, the grading of type $(1|1,1,0,\dots,0)$
 is not even admissible.

Let us now consider case iv). By Proposition \ref{listdepth1},
$L$ has, up to isomorphism, only one grading of depth 1 such that $d$ 
has degree 1, i.e.,  
the grading of type $(1,\dots,1,0,0|0,\dots,0,1,1,1)$. 
The 0th and 1st graded components of $L$ in this
grading are infinite-dimensional vector spaces of growth 2,
and size $m2^{m-2}$ and $\frac{m^2+m-4}{2}2^{m-2}$, respectively.
Hence, by Remark \ref{necessarycond}, the grading of type 
$(1,\dots,1,0,0|0,\dots,0,1,1,1)$ 
 is not even admissible.

By similar arguments one can rule out cases v) and vi).

Let $L=S(1,2)+\F d$. Then, up to isomorphism, one can assume that
the grading of $L$ is either  of
type $(0|1,0)$ or of type $(1|1,1)$.
The first grading is short, thus let us consider the second one.
The 0th and 1st graded components of $L$ in this grading are finite-dimensional
vector spaces of dimension 8 and 9, respectively.
One shows that for every choice of $\mu\in (L_1)_{\0}$,
the centralizer of $\mu$ in $L_0$ has positive dimension, hence
$\mu$ cannot satisfy property (\ref{a}) of Definition \ref{admissible}.

Let $S=SHO(3,3)$. 
By Proposition \ref{listdepth1},
the gradings of depth 1 of $L$ such that $e$ has degree 1, are,
up to isomorphism, 
the gradings of type $(1,1,1|1,1,1)$ and
$(1,0,0|0,1,1)$. The 0th and 1st graded component of $L$ in the grading of type
 $(1,1,1|1,1,1)$ have dimension 17 and 32, respectively, hence
this grading is not even admissible. On the contrary
the grading of type $(1,0,0|0,1,1)$ is even admissible since one can
check, by direct calculations, that it satisfies
Definition \ref{admissible}
with $\mu=e+x_1\xi_2$. Notice that the gradings of type
$(1,1,1|1,1,1)$ and $(1,0,0|0,1,1)$ of $S+\F e$ are isomorphic
to the gradings of type 
$(1,1,1|0,0,0)$  and $(1,0,0|-1,0,0)$, respectively,
via the map $\exp(e)\exp(-f)\exp(e)$ (see also
\cite[Remark 2.38]{CantaK}).

Similarly, if $S=SKO(2,3;1)$, 
then the Lie superalgebra
$L=SKO(2,3;1)+\F e$ has one even admissible grading, namely
 the grading of type $(1,0|0,1,1)$.
In this case $\mu=e+x_1\xi_2$.
\hfill$\Box$

\begin{theorem}\label{moreorbits}
Let $L=S\rtimes\F d$ as in the statement of Theorem \ref{derivations} and let
$L=\prod_{j\geq -1}L_j$ be an even admissible grading of $L$,
of height strictly greater than one, such that $d$ has degree one.
Then a complete list, up to conjugation by $\exp(ad(L_{int}))$,
of elements $\mu\in L_1$ satisfying properties (\ref{a}) and
(\ref{b}), is as follows:
\begin{enumerate}
\item $S=SHO(3,3)$ with the grading of type $(1,0,0|0,1,1)$:
$\mu=e+x_1\xi_2$ and $[L_0,\mu]+\F\mu=L_1$;
\item $S=SKO(2,3;1)$ with the grading of type $(1,0|0,1,1)$: 
$\mu=e+x_1\xi_2$ and $[L_0,\mu]+\F\mu=L_1$.
\end{enumerate}
\end{theorem}
{\bf Proof.} We will denote by $S=\prod_j S_j$ the principal grading of $S$
and argue as in the proof of Theorem \ref{orbits}.
Let us first consider the case $S=SHO(3,3)$ with the grading of type
$(1,0,0|0,1,1)$. We have: $L_{-1}=A/\F 1$, with $A=\F[[x_2,x_3]]\otimes\Lambda(\xi_1)$;
$L_0=\{f\in \langle x_1, \xi_2, \xi_3\rangle\otimes A ~|~ \Delta(f)=0\}$;
$L_1=\{f\in \langle x_1^2, x_1\xi_2, x_1\xi_3, \xi_2\xi_3\rangle\otimes A ~|~ \Delta(f)=0\}$.
Let $\mu\in(L_1)_{\0}$ satisfy properties (\ref{a}) and (\ref{b}).
Then we may assume, up to a linear change of indeterminates, that
$\mu=x_1\xi_2+h$ for some $h\in\prod_{j\geq 1}S_j$. Indeed, if 
$\mu$ lies in $L_1\cap\prod_{j\geq 1}S_j$, then $[L_{-1},\mu]$ is
contained in $\prod_{j\geq 0}S_j$ hence it cannot generate $L_0$, 
contradicting property (\ref{b}). Up to automorphism in $\exp(ad(L_{int}))$
we may thus assume $\mu=x_1\xi_2+\alpha e+h'$ for some
$h'\in\prod_{j\geq 2}S_j$ and some $\alpha\in\F$. By Theorem \ref{eag}$i)$,
we have $\alpha\neq 0$, hence we may assume $\alpha=1$. One checks, by direct computations,
that the element $\mu_0=x_1\xi_2+e$ satisfies
properties (\ref{a}) and (\ref{b}). Besides $[L_0, \mu_0]=[L_{int},\mu_0]$,
since $[\xi_i, \mu_0]=0$ for $i=1,2$. Then statement 1. follows from 
Lemma \ref{lemma}.

Now let $S=SKO(2,3;1)$ with the grading of type $(1,0|0,1,1)$. We have:
$L_{-1}=\F[[x_2]]\otimes\Lambda(\xi_1)=:A$, $L_0=\{f\in\langle
x_1,\xi_2, \tau\rangle\otimes A ~|~ div_1(f)=0\}$,
$L_1=\{f\in\langle
x_1^2, x_1\xi_2, x_1\tau, \xi_2\tau\rangle\otimes A ~|~ div_1(f)=0\}$.
Let $\mu\in(L_0)_{\bar{0}}$ satisfy properties (\ref{a}), (\ref{b}).
Then $\mu=z_0+z_1+h$ for  some 
$h\in\prod_{j\geq 2}S_j$ and some $z_i\in S_i$, $i=1,2$, such that
$z_0+z_1\neq 0$, since if $\mu$ lies in $\prod_{j\geq 2}S_j$,
then $[L_{-1},\mu]\subset L_0\cap \prod_{j\geq 0}S_j$ contradicting
(\ref{b}). Suppose that $\mu=x_1\xi_2+h'$ for some $h'\in\prod_{j\geq 1}S_j$.
Then, we may assume, up to automorphisms,
$\mu=x_1\xi_2+\alpha e+h$ for some 
$h\in\prod_{j\geq 3}S_j$ and some $\alpha\in\F$.
Besides, by Theorem \ref{eag}$i)$, we have $\alpha\neq 0$, hence we may suppose
$\alpha=1$. A direct computation shows that the element 
$\mu_0=x_1\xi_2+e$ satisfies (\ref{a}) and (\ref{b}); besides,
$[L_0,\mu]=[L_{int},\mu]$ since $[\xi_2, \mu]=0$. By Lemma \ref{lemma}
$\mu$ is conjugate to $\mu_0$.

Now suppose that $z_0=0$, i.e., $\mu=Ax_1^2\xi_1+Bx_1\tau+Cx_1x_2\xi_2+h$
for some $h\in\prod_{j\geq 2}S_j$ and some $A,B,C\in\F$ such that
$2A-B+C=0$. Notice that if $B=0$, then $[L_{-1},\mu]\subset L_0\cap \prod_{j\geq 0}S_j$ contradicting (\ref{b}), hence we may assume $B=1$.
It follows that, up to automorphisms, we may assume
$\mu=x_1\tau+Ax_1^2\xi_1+Cx_1x_2\xi_2+\alpha e+h'$ for some
$h'\in\prod_{j\geq 3}S_j$ and some $\alpha, A,C\in\F$ such that $2A+C-1=0$.
By Theorem \ref{eag}$i)$ we have $\alpha\neq 0$.
Then one checks that $L_1\cap S_1$ is not contained in
$[L_0,\mu]+\F\mu$, hence contradicting property (\ref{a}).
\hfill$\Box$

\begin{proposition}\label{outerforodd}\em Let $S$ be a simple linearly compact 
Lie superalgebra. 
\begin{enumerate}
\item[a)] Let
$L=S+\F \mu+\F[\mu,\mu]$, where 
$\mu$ is an odd outer derivation
of $S$ such that $[\mu,\mu]\neq 0$. Then there is no 
admissible grading of $L$ such that $\mu$ has degree 1.
\item[b)] If $L=S\otimes\Lambda(1)+\F\mu+\F d$, where $d$ is an 
even outer derivation of $S$ and $\mu=d\otimes\xi+1\otimes d/d\xi$, 
then there exists no
$\Z$-grading of $L$ of depth 1 such that $\mu$ has degree 1.
\end{enumerate}
\end{proposition}
{\bf Proof.} a)
From  the description of outer
derivations of simple linearly compact Lie superalgebras given in
\cite[Proposition 5.1.2]{K2}, \cite{K3}
and \cite[Proposition 1.8]{CantaK}, 
we deduce that
$L=S\rtimes\F d$, where the following possibilities for $S$ and $d$ may 
occur, up to isomorphism:
\begin{itemize}
\item[i)] $S=q(n)$ and $d=D$ (see \cite[Proposition 5.1.2(c)]{K2} for
the definition of $D$); 
\item[ii)] $S=H(0,n)$ for some odd $n\geq 5$ and $d=\xi_1\dots\xi_n$;
\item[iii)] $S=S(1,n)$
for some odd $n\geq 3$ and $d=\xi_1\dots\xi_n\frac{\partial}{\partial x}$;
\item[iv)] $S=SKO(m,m+1;(m-2)/m)$ for some even $m\geq 2$ and $d=\xi_1\dots\xi_m$;
\item[v)] $S=SKO(m,m+1;1)$ for some odd $m\geq 3$ and $d=\tau\xi_1\dots\xi_m$;
\item[vi)] $S=SHO(m,m)$ for some even $m\geq 4$ and $d=\xi_1\dots\xi_m$.
\end{itemize}
In all cases $[d,d]=0$, hence $\mu=d+z$ for some odd inner derivation 
$z\neq 0$, 
and $d$ and $z$ have degree 1.
In case i) it follows immediately from the definition that
every $\Z$-grading of depth 1 of $L$ is short, 
hence there is no such $\mu$.
Let us consider case iii). By Proposition \ref{listdepth1},
the gradings of depth 1 of $L$ such that $d$ has degree 1, are,
up to isomorphism, 
the gradings of type $(0|1,0,\dots,0)$ and
$(1|1,1,0,\dots,0)$. The first one is short, hence there is no odd derivation
$\mu$ satisfying the hypotheses.
Let us now consider the grading of type $(1|1,1,0,\dots,0)$. The 0th and 1st graded components of $S$ in this
grading have dimension $(n+6)2^{n-2}$ and $(3+3n)2^{n-2}$, respectively.
Hence, by Remark \ref{necessarycond}, the grading of type $(1|1,1,0,\dots,0)$
 is not admissible. Similarly one rules out cases ii), v), vi) and
iv) unless $S=SKO(2,3;0)$ with the grading of type $(0,0|1,1,1)$.
Since this grading is short we get statement a). 

Now let $L$, $\mu$ and $d$ be as in b) and consider a $\Z$-grading of depth 1
of $L$ such that $\mu$ has degree 1.
 Then $\xi$
has degree $-1$, hence $S=\prod_{j\geq 0}S_j$, and $d$ has degree 2. 
Since $S$ is simple
it follows that $S=S_0$, therefore there is no outer derivation of $S$
of degree 2. Statement b) follows.
\hfill$\Box$

\begin{remark}\label{J(g,mu)}\em Let $(\g,\mu)$ be one of the pairs listed in Tables 1
and 2, where $\g=\prod_{i\geq -1}\g_i$ 
is a simple linearly compact Lie superalgebra with an admissible grading, and
$\mu\in \g_1$ satisfies Definition \ref{admissible}.
According to Proposition \ref{oddJordan}, $J(\g,\mu)=\g_{-1}$ 
with product $f\circ g=[[\mu,f],g]$.
We hence get the following corresponding list of rigid superalgebras 
(constructed in Sections \ref{evenexamples} and \ref{oddexamples}),
where the bar over $J$ and $\mu$ means parity reversal, described in Remark
\ref{rem:1.4}:
\begin{enumerate}
\item $\g=W(0,3)$ with the grading of type $(|1,1,0)$, $J(\g,\mu)=JW_{0,4}$;
\item $\g=W(0,4)$ with the grading of type $(|1,1,0,0)$, $J(\g,\mu)=JW_{0,8}$;
\item $\g=S(0,4)$ with the grading of type $(|1,1,0,0)$, $J(\g,\mu)=JS_{0,8}$;
\item $\g=S(0,5)$ with the grading of type $(|1,1,0,0,0)$, $J(\g,\mu)=JS_{0,16}$;
\item $\g=S(1,4)$ with the grading of type $(0|1,1,0,0)$, $J(\g,\mu)=JS^\alpha_{1,8}$;
\item $\g=S(2,0)$ with the grading of type $(1,1|)$, $J(\g,\mu)=JS_{0,2}$;
\item $\g=S(2,0)$ with the grading of type $(1,0|)$, $J(\g,\mu)$ is the Beltrami algebra $JS_{1,1}$;
\item $\g=W(1,1)$ with the grading of type $(1|0)$, $J(\bar{\g},\bar{\mu})=
LW_{0,2}$;
\item $\g=W(1,2)$ with the grading of type $(0|1,1)$, $J(\bar{\g},\bar{\mu})=LW^\alpha_{1,2}$;
\item $\g=W(2,1)$ with the grading of type $(1,0|0)$, $J(\bar{\g},\bar{\mu})=LW_{1,2}$;
\item $\g=W(2,2)$ with the grading of type $(0,0|1,1)$, $J(\bar{\g},\bar{\mu})=LW^\alpha_{2,2}$;
\item $\g=S(2,2)$ with the grading of type $(0,0|1,1)$, $J(\bar{\g},\bar{\mu})=LS^\alpha_{2,2}$;
\item $\g=S(1,3)$ with the grading of type $(0|1,1,1)$, $J(\bar{\g},\bar{\mu})=LS_{1,3}$;
\item $\g=H(2,2)$ with the grading of type $(1,0|1,0)$, $J(\bar{\g},\bar{\mu})=LH^\alpha_{1,2}$;
\item $\g=HO(2,2)$ with the grading of type
$(1,0|0,1)$, $J(\bar{\g},\bar{\mu})=LHO_{1,2}$;
\item $\g=HO(3,3)$ with the grading of type
$(0,0,0|1,1,1)$, $J(\bar{\g},\bar{\mu})=LHO^\alpha_{3,1}$;
\item $\g=HO(n+1,n+1)$, $n\geq 1$, with the grading of type
$(0,\dots,0,1|0,\dots,0,-1)$, $J(\bar{\g},\bar{\mu})=LP(n-1,n-1)$;
\item $\g=SHO(3,3)$ with the grading of type $(1,0,0|0,1,1)$,
$J(\bar{\g},\bar{\mu})=LSHO'_{2,2}$;
\item $\g=SHO(4,4)$ with the grading of type $(0,0,0,0|1,1,1,1)$,
$J(\bar{\g},\bar{\mu})=LSHO^\alpha_{4,1}$;
\item $\g=SHO(n+1,n+1)$, $n\geq 2$, with the grading of type
$(0,\dots,0,1|0,\dots,0,-1)$, $J(\bar{\g},\bar{\mu})=LSHO_{n,2^{n-1}}$;
\item $\g=KO(n+1,n+2)$, $n\geq 1$, with the grading of type
$(0,\dots,0,1|0,\dots,0,-1,0)$, $J(\bar{\g},\bar{\mu})=LP(n-1,n)$;
\item $\g=KO(2,3)$ with the grading of type $(0,0|1,1,1)$, $J(\bar{\g},\bar{\mu})=LKO_{2,1}$;
\item $\g=SKO(n,n+1;\beta)$, $n\geq 2$, $\beta\neq \frac{4}{n}$,
 with the grading of type
$(0,\dots,0,1|0,\dots,0,$ $-1,0)$, $\mu=x_n+x_n\xi_{n-1}\tau-\frac{2-n\beta}
{2}x_n^2\xi_{n-1}\xi_n$,
 $J(\bar{\g},\bar{\mu})=LSKO_{n,2^n}$;
\item $\g=SKO(2,3;\beta)$, $\beta\neq 0,1, 2+\frac{1}{b}, 2+\frac{2}{b}$,
$\forall b\in\Z_{>0}$, with the grading of type $(0,1|0,-1,0)$, $\mu=x_2+x_2^2\xi_2\tau-2x_1x_2\xi_1\tau+(3-2\beta)
x_1x_2^2\xi_1\xi_2$, $J(\bar{\g},\bar{\mu})=LSKO'_{1,2}$;
\item $\g=SKO(3,4;1)$, 
 with the grading of type $(0,0,1|0,0,-1,0)$, $\mu=x_3+x_3\xi_1\xi_2+x_3\xi_2\tau+\frac{1}{2}x_3^2\xi_2\xi_3$,
 $J(\bar{\g},\bar{\mu})=LSKO'_{2,4}$;
\item $\g=SKO(3,4;\beta)$, 
 with the grading of type $(0,0,0|1,1,1,1)$, 
 $J(\bar{\g},\bar{\mu})=LSKO^\alpha_{3,1}$.
\end{enumerate}
\end{remark}

Similarly, we have the following remark.

\begin{remark}\label{J(L,mu)}\em
Let $(L,\mu)$ be one of the pairs listed in Theorem \ref{moreorbits}.
Then the corresponding  rigid superalgebras are
as follows:
\begin{enumerate}
\item $L=SHO(3,3)\rtimes\F e$, $J(L,\mu)=JSHO_{2,2}$;
\item $L=SKO(2,3;1)\rtimes\F e$, $J(L,\mu)=JSKO_{1,2}$.
\end{enumerate}
\end{remark}

\begin{remark}\em Note that, since the Lie superalgebras
$HO(2,2)$ and $SKO(2,3;0)$ are isomorphic \cite[\S 0]{CantaK}, 
the anti-commutative
superalgebras  $LSKO_{1,2}$ and $LHO_{1,2}$ are isomorphic. 
\end{remark}

\section{Proofs of Theorems \ref{evenclass} and \ref{oddclass}
and a corollary}\label{proofs}

{\bf Proof of Theorem \ref{evenclass}}
Let $J$ be a simple linearly compact rigid superalgebra
with product $\mu$.
Then, by Proposition \ref{irr}, $Lie(J,\mu)$ is a linearly compact
Lie superalgebra with an even admissible $\Z$-grading 
$Lie(J,\mu)=\prod_{k\geq -1}Lie_k(J)$,
where
$Lie_{-1}(J)=J$, $Lie_0(J)=Str(J,\mu)$ and $Lie_1(J)=R(J,\mu)$.
By Theorem \ref{LieJ} either $Lie(J,\mu)$ is simple or
$Lie(J,\mu)=S\rtimes\F \mu$, where $S$ is a simple Lie superalgebra
and $\mu$ is an even outer derivation of $S$ lying in $Lie_1(J)$.
The statement then follows from Theorem \ref{eag}i), 
Remark \ref{evenshort}, 
 Theorems \ref{orbits} and \ref{moreorbits}, Proposition \ref{indofp},
 and Remarks \ref{J(g,mu)}, \ref{J(L,mu)}. \hfill$\Box$

\bigskip

\noindent
{\bf Proof of Theorem \ref{oddclass}}
Let $J$ be a simple linearly compact rigid odd type superalgebra
with product $\mu$.
Then, by Proposition \ref{irr}, $Lie(J,\mu)$ is a linearly compact
Lie superalgebra with an odd admissible $\Z$-grading 
$Lie(J,\mu)=\prod_{k\geq -1}Lie_k(J)$,
where
$Lie_{-1}(J)=J$, $Lie_0(J)=Str(J,\mu)$ and $Lie_1(J)=R(J,\mu)$.
By Theorem \ref{LieJ} and Proposition \ref{outerforodd},
either $Lie(J,\mu)\cong\xi\a+\a+d/d\xi$, where $\a$ is a simple Lie 
superalgebra and $\xi$ is an 
odd indeterminate, hence ($\bar{J}, \bar{\mu}$) is a Lie superalgebra,
or else $Lie(J,\mu)$ is simple.
The statement then follows from Theorem \ref{eag}ii),
Proposition \ref{[p,p]=0}, Theorem \ref{orbits},
Proposition \ref{indofp} and Remark \ref{J(g,mu)}.
\hfill$\Box$

\begin{proposition}\label{october} Let $P$ be a linearly compact odd 
generalized Poisson superalgebra.
If $P$ is simple, then the rigid odd type superalgebra $OJP$,
constructed in Example \ref{mainexample}, is simple. 
\end{proposition}
{\bf Proof.} We will divide the proof into three steps:
let $I$ be an ideal of $OJP$, first we will show
that if $I$ contains a non-zero element of the form $\eta fa$ for some 
$f\in P$ and $a\in \F[[x]]$, then $I=OJP$; secondly we will show that 
if $I$ contains a non-zero element of the form $fa$ for some 
$f\in P$ and $a\in \F[[x]]$, then $I=OJP$; finally we will show that 
every non-zero ideal of $OJP$ necessarily 
contains either an element of the form $\eta fa$ or an element of the form $fa$
for some $0\neq f\in P$ and some $0\neq a\in \F[[x]]$.

STEP 1. Suppose that $I$ contains 
a non-zero element $\eta fa$ for some 
$f\in P$ and some $a\in \F[[x]]$. Then $\eta\circ \eta fa=-\eta fa'\in I$,
hence, either $a$ is a polynomial and $\eta f\in I$, or
$a$ is an infinite series and we may assume that it is invertible, i.e.,
$a=\sum_{k\geq 0}\alpha_kx^k$ with $\alpha_0\neq 0$.
We have: $\eta x^r\circ \eta fa=\eta f(rx^{r-1}a-x^ra')=\eta f(r\alpha_0x^{r-1}
+\mbox{higher order terms})$. It follows that in the limit
we can cancel out all terms of $\eta fa$ of degree in $x$ greater than 0,
i.e., $\eta f\in I$. Now, for $g,h\in P$ and $b\in \F[[x]]$, we have:
$\eta gb\circ \eta f=(-1)^{p(g)}\eta gfb'\in I$ and $h\circ \eta f=
\eta\{h,f\}-\frac{1}{2}(-1)^{p(h)}\eta D(h)f\in I$, hence
$\eta\{h,f\}\in I$. Let $J=\{\varphi\in P~|~\eta\varphi\in I\}$. Then
$J\neq 0$ since $f\in J$. Besides, we have just shown that 
$J$ is an ideal of $P$ with respect to both
the associative and the Poisson product. By the simplicity of
$P$, $J=P$, i.e., $\eta P\subset I$. Now, for $a\in \F[[x]]$ and $g\in P$, 
we have: $\eta a\circ \eta g=\eta ga'\in I$, hence $\eta P[[x]]\subset I$.
Likewise, $a\circ\eta g=-ga'+\eta D(g)(a-\frac{1}{2}xa')\in I$, hence
$ga'\in I$, i.e., $P[[x]]\subset I$. It follows that $I=OJP$.

STEP 2. Suppose that $I$ contains 
a non-zero element $fa$ for some 
$f\in P$ and some $a\in \F[[x]]$. Then $1\circ fa=-D(f)a+2\eta fa\in I$, therefore
if $D(f)=0$, then $I=OJP$ by Step 1. Now suppose $D(f)\neq 0$.
We have: \begin{equation}\underbrace
{\eta\circ(\eta \circ \dots(\eta}_{r ~\mbox\small{times}}\circ(1\circ fa)))=
-D(f)a^{(r)}+2(-1)^r\eta fa^{(r)}\in I,
\label{R*}
\end{equation} where we denoted by $a^{(r)}$ the
$r$-th derivative of $a$ with respect to $x$. Likewise,
\begin{equation}
\underbrace
{\eta\circ(\eta \circ \dots(\eta}_{r ~\mbox\small{times}}\circ fa))=(-1)^r
(fa^{(r)}-\frac{r}{2}\eta D(f)a^{(r-1)})\in I.
\label{R**}
\end{equation}
If $a$ is a polynomial then there exists some $s\in\Z_{>0}$ such that 
$a^{(s-1)}\neq 0$ and $a^{(s)}=0$. By (\ref{R**}), $\eta D(f)a^{(s-1)}\in I$, hence $I=OJP$ by Step 1. If $a$ is an infinite series then 
there exists $s\in \Z_{\geq 0}$ such that $a^{(s)}$ is invertible.
For $r\in\Z_{\geq 0}$ we have:
$x^r\circ (-D(f)a^{(s)}+2(-1)^s\eta fa^{(s)})=
-2(-1)^srfx^{r-1}a^{(s)}+(-2-r+2(-1)^s)\eta D(f)x^ra^{(s)}\in I$. 
Since $a^{(s)}$ is invertible,
it follows that $f+\eta D(f)b\in I$ for some invertible $b\in \F[[x]]$.
We have:
$\eta x\circ(1\circ(f+\eta D(f)b))=2\eta f\in I$, hence $I=OJP$ by Step 1. 

STEP 3. Suppose that $I$ contains an element $z=fa+\eta gb$ for some 
non-zero $f,g\in P$ and some non-zero $a,b\in \F[[x]]$. Then
\begin{equation}\underbrace
{\eta\circ(\eta \circ \dots(\eta}_{r ~\mbox\small{times}}\circ(fa+\eta gb)))=
-fa^{(r)}+\eta(\frac{r}{2}D(f)a^{(r-1)}-gb^{(r)})\in I.
\label{R****}
\end{equation}
If $a$ is a polynomial then there exists $s\in\Z_{>0}$ such that 
$a^{(s-1)}\neq 0$ and $a^{(s)}=0$. 
Then $\eta(\frac{s}{2}D(f)a^{(s-1)}-gb^{(s)})\in I$. If 
$\frac{s}{2}D(f)a^{(s-1)}-gb^{(s)}\neq 0$, then
$\eta\circ \eta(\frac{s}{2}D(f)a^{(s-1)}-gb^{(s)})=\eta gb^{(s+1)}\in I$.
Therefore, either $b^{(s+1)}\neq 0$ and $I=OJP$ by Step 1, or
$b^{(s)}=\beta\in \F$. Note that $a^{(s-1)}=\alpha\in \F$ since
$a^{(s)}=0$, i.e., $\eta(\frac{s}{2}\alpha D(f)+\beta g)\in I$ and,
again, $I=OJP$ by Step 1. Now suppose that 
$\frac{s}{2}D(f)a^{(s-1)}-gb^{(s)}=0$. Then, either $D(f)\neq 0$ or $D(f)=0$.
If $D(f)\neq 0$ then $g=\gamma D(f)$ 
and $a^{(s-1)}=\delta b^{(s)}\in\F$, for some $\gamma,
\delta\in \F$, $\gamma\neq 0\neq\delta$.
By (\ref{R****}), $fa^{(s-1)}+\eta(-\frac{s-1}{2}D(f)a^{(s-2)}+gb^{(s-1)})\in I$,
hence $1\circ (fa^{(s-1)}+\eta(-\frac{s-1}{2}D(f)a^{(s-2)}+gb^{(s-1)}))=
-D(f)a^{(s-1)}+2\eta fa^{(s-1)}\in I$, i.e., $-D(f)+2\eta f\in I$. It follows
that $\eta x\circ (-D(f)+2\eta f)=2\eta f\in I$ and $I=OJP$ by Step 1.
If $D(f)=0$, then either $g=0$ or $g\neq 0$ and $b^{(s)}=0$. If $g=0$, then,
by (\ref{R****}), $fa^{(s-1)}\in I$ hence $I=OJP$ by Step 2.
If $g\neq 0$, then $b^{(s)}=0$, i.e., $b^{(s-1)}\in\F$, i.e., we may assume,
by (\ref{R****}) with $r=s-1$, that $f+\eta g\in I$, hence
either $g=0$ and $I=OJP$ by Step 2, or $g\neq 0$ and $\eta x\circ (f+\eta g)=\eta g\in
I$, hence $I=OJP$ by Step 1.

Now suppose that $a$ is an infinite series, then, by (\ref{R****}),
we may assume that $fa+\eta gb\in I$ for some invertible $a\in \F[[x]]$,
some non-zero $f,g\in P$ and some non-zero $b\in \F[[x]]$.
For $k\in\Z_{\geq 0}$, we have:
$\eta x^k\circ (fa+\eta gb)=-fa'x^k+\frac{1-k}{2}\eta D(f)ax^k\in I$.
 It follows that in the limit we
can cancel out all terms of $fa$ of positive degree in $x$, and get
$f+\eta (D(f) c+gd)\in I$ for some  $c, d\in \F[[x]]$. 
We have: $1\circ(f+\eta D(f) c+\eta gd)=-D(f)+2\eta f+\eta D(g)d\in I$ and
$\eta\circ(1\circ(f+\eta D(f) c+\eta gd))=-\eta D(g)d'\in I$. Therefore
 if $D(g)d'\neq 0$, then $I=OJP$ by Step 1, otherwise, either
$D(g)=0$ or $D(g)\neq 0$ and $d'=0$. If $D(g)=0$, then $-D(f)+2\eta f\in I$,
hence $\eta x\circ(-D(f)+2\eta f)=2\eta f\in I$, and $I=OJP$ by Step 1.
If $D(g)\neq 0$ and $d'=0$, i.e., $d\in\F$, then, either $d=0$ and we proceed
as above, or $d\neq 0$ and $\eta x\circ(-D(f)+2\eta f+\eta D(g)d)=
2\eta f+\eta D(g)d\in I$. If $2f+D(g)d\neq 0$, $I=OJP$ by Step 1.
If $2f+D(g)d=0$, then $D(f)\in I$ hence if $D(f)\neq 0$ then $I=OJP$ by
Step 2. Finally, if $2f+D(g)d=0$ and $D(f)=0$, then 
$f+\eta gd\in I$ hence $\eta x\circ (f+\eta gd)=\eta gd$, since $d'=0$,
 and the result follows
from Step 1.
\hfill$\Box$

\begin{corollary} \begin{itemize}
\item[a)] Let $(P,\{\cdot,\cdot\})$ be a simple linearly compact
odd generalized Poisson superalgebra. Then the rigid odd type superalgebra
$OJP$ is isomorphic to one of the superalgebras $OJP(n,n)$, $OJP(n,n+1)$
, constructed in Example \ref{mainexample}.
\item[b)] Any simple linearly compact odd generalized Poisson superalgebra is
gauge equivalent to $PO(n,n)$ or $PO(n,n+1)$.
\item[c)] Any simple linearly compact odd Poisson superalgebra is
isomorphic to $PO(n,n)$.
\end{itemize}
\end{corollary}
{\bf Proof.} 
It will be convenient to talk here about odd type commutative superalgebras 
instead of anti-commutative superalgebras (related to each other by reversal of
parity, as in Remark \ref{rem:1.4}); we shall denote the former superalgebras
by $OJX_{m,n}$ if the latter are denoted by $LX_{m,n}$.
Let $(P,\{\cdot,\cdot\})$ be a simple linearly compact odd generalized Poisson 
superalgebra and suppose that $P$ is non-trivial. We showed in Example
\ref{mainexample} that
 $OJP$ is a  rigid odd type superalgebra, which is linearly compact by
construction. By Proposition \ref{october}, $OJP$ is simple
 hence it is isomorphic to one of the superalgebras listed in Theorem 
\ref{oddclass}.  We have $OJP=P[[x]]+\eta P[[x]]$ where 
$\eta$ is an odd variable. It follows that $OJP$ cannot be 
isomorphic to any of the following rigid superalgebras:
\begin{enumerate}
\item[1)] $OJW_{0,2}$ since it is finite-dimensional;
\item[2)] $OJW^\alpha_{1,2}$, $OJW^\alpha_{2,2}$, $OJS^\alpha_{2,2}$, $OJS_{1,3}$, 
$OJHO^\alpha_{3,1}$, $OJSHO^\alpha_{4,1}$, $OJKO_{2,1}$, $OJSKO^\alpha_{3,1}$,
 since these are
completely odd superalgebras;
\item[3)] $OJW_{1,2}$,  $OJSKO'_{1,2}$, 
 (resp.\ $OJH^\alpha_{1,2}$, $OJHO_{1,2}$) 
since these are equal to $\F[[x]]\otimes\Lambda(1)$ 
(resp.\ to $(\F[[x]]\otimes\Lambda(1))/\F 1$, $(\F[[x]]\otimes\Lambda(1))/\F\xi$);
\item[4)] $OJSHO'_{2,2}$, since $OJSHO'_{2,2}=(\F[[x_1,x_2]]\otimes\Lambda(1))/\F\xi$,
 hence if $OJP$ were isomorphic to $OJSHO'_{2,2}$,
$P$ would be completely even.
\end{enumerate}
Now let us consider the rigid odd type superalgebra $OJSHO_{n,2^{n-1}}$.
We have: $OJSHO_{n,2^{n-1}}=\{f\in {\cal O}(x_1,\dots ,x_n,\xi_1,\dots,\xi_n) ~|~ 
\Delta(f)=0\}$.
Suppose that there exists a super-space $Q$ such that $OJSHO_{n,2^{n-1}}=Q[[x]]+
\eta Q[[x]]$
for some odd element $\eta\in OJSHO_{n,2^{n-1}}$.
Then every element $f\in {\cal O}(x_1,\dots,x_n)$ lies in $Q[[x]]$,
 since there exists no odd element $g$ such that $f =\eta g$. 
In particular $x_i$ lies in $Q[[x]]$ for every $i$. 
It follows that $\eta x_i$ lies in $\eta Q[[x]]$, hence 
$\Delta(\eta x_i)=0$. 
But $\Delta(\eta x_i)= \Delta(\eta)x_i+\frac{\partial \eta}{\partial \xi_i}
=\frac{\partial \eta}{\partial \xi_i}$ 
($\eta\in OJSHO_{n,2^{n-1}}$, hence $\Delta(\eta)=0$). 
It follows that $\eta$ lies in ${\cal O}(x_1,\dots,x_n)$ 
and this is a contradiction since $\eta$ is odd. Therefore $OJP$ cannot be isomorphic to $OJSHO_{n,2^{n-1}}$
for any $n$.

Likewise, we will show that $OJP$ cannot be isomorphic to
the rigid superalgebra $OJSKO_{n,2^n}$ for any $n$. 
Suppose that $OJSKO_{n,2^n}=Q[[x]]+
\eta Q[[x]]$ for some super-space $Q$ and
some odd element $\eta\in OJSKO_{n,2^n}$.
According to Example \ref{OJSKO}, we have $\Delta(\eta)+(E-n\beta)
\frac{\partial\eta}{\partial\tau}+(1-\beta)\frac{\partial\eta}{\partial\tau}=0$. As above, $x_i$ lies in $Q[[x]]$ for every $i$, hence
$\eta x_i\in \eta Q[[x]]$. It follows that
$\Delta(\eta x_i)+(E-n\beta)\frac{\partial(\eta x_i)}{\partial\tau}+(1-\beta)\frac{\partial(\eta x_i)}{\partial\tau}=0$, i.e., 
$\frac{\partial\eta}{\partial\xi_i}=-x_i\frac{\partial\eta}{\partial\tau}$.
Let $\eta=\tau\eta_1+\eta_2$ for some $\eta_1, \eta_2\in
{\cal O}(x_1,\dots,x_n,\xi_1,\dots,\xi_n)$. Then, for $i=1,\dots,n$,
we have:
$-\tau\frac{\partial\eta_1}{\partial\xi_i}+\frac{\partial\eta_2}{\partial\xi_i}
=-x_i\eta_1$, hence $\frac{\partial\eta_1}{\partial\xi_i}=0 ~\forall ~i$, i.e.,
$\eta_1\in{\cal O}(x_1,\dots, x_n)$, and $\frac{\partial\eta_2}{\partial\xi_i}=
-x_i\eta_1$. It follows that $\eta_2=-\eta_1\sum_{k=1}^n x_k\xi_k$,
i.e., $\eta=\tau\eta_1-\eta_1\sum_{k=1}^n x_k\xi_k$.
Therefore, for $i=1,\dots,n$, the elements $\xi_i$'s, lying in $OJSKO_{n,2^n}$,
lie in $Q[[x]]$, since there exists no even element $g$ such that
$\xi_i=\eta g$. It follows that $\eta\xi_i\in \eta Q[[x]]$,
hence $\Delta(\eta\xi_i)+(E-n\beta)(\eta_1\xi_i)+(1-\beta)(\eta_1\xi_i)=0$,
i.e., $-\tau\frac{\partial\eta_1}{\partial x_i}+
\frac{\partial\eta_1}{\partial x_i}(\sum_{k=1}^n x_k\xi_k)+2\eta_1\xi_i=0$.
It follows that $\frac{\partial\eta_1}{\partial x_i}=0 ~\forall ~i$, hence
$\eta_1\in\F$, i.e., $\eta=\tau-\sum_{k=1}^n x_k\xi_k$.
We have $\Delta(\eta)+(E-n\beta)(\frac{\partial\eta}{\partial\tau})+(1-\beta)
(\frac{\partial \eta}{\partial\tau})=
-n+1-(n+1)\beta=0$, which is impossible for every $\beta\neq\frac{1-n}{n+1}$.
Besides, if $\beta=\frac{1-n}{n+1}$, 
$\Delta(\eta x_i)+(E-n\beta)(\frac{\partial (\eta x_i)}{\partial\tau})+(1-\beta)
(\frac{\partial (\eta x_i)}{\partial\tau})=
x_i$ and this is a contradiction since $\eta x_i\in OJSKO_{n,2^n}$.
We conclude that $OJP$ cannot be isomorphic to $OJSKO_{n,2^n}$
for any $n$. Similarly, $OJP$ cannot be isomorphic to $OJSKO'_{2,4}$.
Statement $a)$ follows.

In order to prove $b)$, let us consider a simple odd generalized Poisson 
superalgebra $P$ and the corresponding rigid odd type superalgebra
 $OJP=P[[x]]+\eta P[[x]]$ with product (\ref{goddP}). By $a)$, 
$OJP$ is isomorphic either to $OJP(n,n)$ or to $OJP(n,n+1)$, for some
$n\in\Z_{\geq 0}$. Suppose
$OJP\cong OJP(n,n)$ and let us thus identify
$OJP$ with $OJP(n,n)$. Then also the corresponding Lie superalgebras
$Lie(OJP,\mu')$ and $Lie(OJP(n,n),\mu)$ can be identified, hence there exists
an automorphism $s$ of $HO(n+2, n+2)$, preserving the grading of type
$(0,\dots,0,1|0,\dots,0,-1)$ and the decomposition 
$OJP(n,n)=PO(n,n)[[x_{n+1}]]+
\xi_{n+1} PO(n,n)[[x_{n+1}]]$,
 such that $s(\mu)=\mu'$. From Example \ref{HO} one deduces that
$\mu'=\mu\varphi$ for some invertible $\varphi\in {\cal O}(x_1,\dots,
x_n, \xi_1, \dots, \xi_n)$ such that $\{\varphi,\varphi\}=0$. It follows that, for $f,g\in OJP$,
$f\circ g=[[\mu',f],g]=(-1)^{p(f)+1}\{f,g\}^{\varphi}+2\eta fg$.
By Remark \ref{reconstruct}, $P$ 
is hence gauge equivalent to $PO(n,n)$.
A similar argument shows that if $OJP$ is isomorphic to $OJP(n,n+1)$ then $P$
is gauge equivalent to $PO(n,n+1)$.

Finally, $c)$ follows from $b)$.
Indeed, if $D^{\varphi}=0$, then $D(\varphi)a=\{\varphi,a\}$ for all $a$.
Letting $a=\varphi$, we get $D(\varphi)=0$, hence $\{\varphi,a\}=0$ for all $a$,
hence $\varphi\in\F$.
\hfill$\Box$
 
$$$$ 

\end{document}